\documentclass[11pt]{article}
\usepackage{amssymb}

\usepackage{amsmath}
\usepackage{theorem}
\usepackage{here}
\usepackage[dvipdfmx]{color}

\newtheorem{thm}{Theorem}[section]
\newtheorem{lem}{Lemma}[section]
\newtheorem{cor}{Corollary}[section]
\newtheorem{prp}{Proposition}[section]
\theorembodyfont{\rmfamily}

\newtheorem{remark}{Remark}[section]
\newtheorem{assumption}{Assumption}[section]
\makeatletter

\@addtoreset{equation}{section}
\makeatother

\usepackage{comment} %
\usepackage{bm}
\usepackage[pdftex]{graphicx}
\def\Et{{\widetilde E}}

\def\It{{\widetilde I}}

\def\Xt{{\widetilde X}}
\def\Wt{{\widetilde W}}
\def\Ct{{\widetilde C}}
\def\ta{{\tau}}
\def\xt{{\tilde x}}
\def\at{{\tilde a}}
\def\bt{{\tilde b}}
\def\ct{{\tilde c}}
\def\dt{{\tilde d}}
\def\st{{\tilde s}}
\def\pt{{\tilde p}}

\def\al{{\alpha}}
\def\be{{\beta}}
\def\ga{{\gamma}}
\def\de{{\delta}}

\def\la{{\lambda}}

\def\om{{\omega}}
\def\th{{\theta}}

\def\bga{{\text{\boldmath $\gamma$}}}

\def\bde{{\text{\boldmath $\delta$}}}

\def\mut{{\tilde \mu}}

\def\gat{{\tilde \ga}}
\def\bgat{{\widetilde \bga}}

\def\De{{\Delta}}

\def\Ga{{\Gamma}}

\def\La{{\Lambda}}
\def\a{{\text{\boldmath $a$}}}
\def\b{{\text{\boldmath $b$}}}
\def\c{{\text{\boldmath $c$}}}
\def\d{{\text{\boldmath $d$}}}
\def\e{{\text{\boldmath $e$}}}

\def\i{{\text{\boldmath $i$}}}
\def\j{{\text{\boldmath $j$}}}
\def\k{{\text{\boldmath $k$}}}

\def\p{{\text{\boldmath $p$}}}

\def\w{{\text{\boldmath $w$}}}
\def\x{{\text{\boldmath $x$}}}
\def\y{{\text{\boldmath $y$}}}

\def\W{{\text{\boldmath $W$}}}
\def\X{{\text{\boldmath $X$}}}
\def\Y{{\text{\boldmath $Y$}}}
\def\Z{{\text{\boldmath $Z$}}}

\def\fh{{\hat f}}
\def\gh{{\hat g}}

\def\ph{{\hat p}}

\def\pbh{{\hat \p}}
\def\abt{{\tilde \a}}
\def\bbt{{\tilde \b}}
\def\cbt{{\tilde \c}}

\def\dt{{\tilde d}}
\def\wt{{\tilde w}}

\def\Wbt{{\widetilde \W}}

\def\[{{\text{\boldmath $[$}}}
\def\]{{\text{\boldmath $]$}}}

\def\/{{\Bigr/\!\!}}

\def\1r{{\rm (1)}}
\def\2r{{\rm (2)}}
\def\3r{{\rm (3)}}
\def\4r{{\rm (4)}}
\def\5r{{\rm (5)}}

\def\non{{\nonumber}}
\begin{document}
\title{Bayesian Shrinkage Approaches to Unbalanced Problems of Estimation and Prediction on the Basis of Negative Multinomial Samples\footnote{This preprint has not undergone peer review (when applicable) or any post-submission improvements or corrections. The Version of Record of this article is published in Japanese Journal of Statistics and Data Science, and is available online at https://doi.org/10.1007/s42081-021-00141-z. }
}
\author{
Yasuyuki Hamura\footnote{Graduate School of Economics, University of Tokyo, 
7-3-1 Hongo, Bunkyo-ku, Tokyo 113-0033, JAPAN.\newline{
JSPS Research Fellow. 
E-Mail: yasu.stat@gmail.com}} \
}
%
%
\maketitle
\begin{abstract}
In this paper, we treat estimation and prediction problems where negative multinomial variables are observed and in particular consider unbalanced settings. 
First, the problem of estimating multiple negative multinomial parameter vectors under the standardized squared error loss is treated and a new empirical Bayes estimator which dominates the UMVU estimator under suitable conditions is derived. 
Second, we consider estimation of the joint predictive density of several multinomial tables %
under the Kullback-Leibler divergence and obtain a sufficient condition under which the Bayesian predictive density with respect to a hierarchical shrinkage prior dominates the Bayesian predictive density with respect to the Jeffreys prior. 
Third, our proposed Bayesian estimator and predictive density give risk improvements in simulations. 
Finally, the problem of estimating the joint predictive density of negative multinomial variables is discussed. 
\par\vspace{4mm}
\end{abstract}

\section{Introduction}
\label{sec:introduction}
Properties of shrinkage estimators based on count variables have been extensively investigated within the decision-theoretic framework since the seminal work of Clevenson and Zidek (1975). 
For example, as briefly reviewed in Section 1 of Hamura and Kubokawa (2020b), estimation of Poisson parameters was studied by Ghosh and Parsian (1981), Tsui (1979b), Tsui and Press (1982), and Ghosh and Yang (1988) in various settings while Tsui (1979a), Hwang (1982), and Ghosh, Hwang, and Tsui (1983) showed that similar results hold for discrete exponential families. 
Extending the result of Tsui (1984) and Tsui (1986a), Tsui (1986b) proved that Clevenson--Zidek-type estimators dominate the usual estimator %
in the case of the negagive multinomial distribution, which is a generalization of the negative binomial distribution and is a special case of the general distributions of Chou (1991) and Dey and Chung (1992). 
More recent studies include Chang and Shinozaki (2019), Stoltenberg and Hjort (2019), and Hamura and Kubokawa (2019b, 2020b, 2020c). 
On the other hand, since Komaki (2001), Bayesian predictive densities with respect to shrinkage priors have been shown to dominate those based on noninformative priors and parallels between estimation and prediction have been noted in the literature. 
In particular, Komaki (2004, 2006, 2015) and Hamura and Kubokawa (2019b) obtained dominance conditions in the Poisson case. 

There are still directions in which these results could be generalized further. 
First, although sample sizes will be unbalanced in many practical situations, some of the results are applicable only to the balanced case. 
Weights in loss functions may also be unbalanced in practice (see, for example, Section 7 of Stoltenberg and Hjort (2019)). 
Second, as pointed out by Hamura and Kubokawa (2020b), decision-theoretic properties  of Bayesian procedures have not been fully studied for discrete distributions other than the Poisson distribution. 
Even in the Poisson case, %
it was only after the work of Komaki (2015) that many Bayesian shrinkage estimators were shown to dominate usual estimators in the presence of unbalanced sample sizes (Hamura and Kubokawa (2019b, 2020c)). 
Third, while theoretical properties of Bayesian predictive densities for Poisson models have been investigated in several papers as mentioned earlier, relatively few researchers (Komaki (2012), Hamura and Kubokawa (2019a)) have considered predictive density estimation for other discrete exponential families. 
In this paper, we treat these three issues when considering Bayesian estimators and predictive density estimators based on negative multinomial observations in unbalanced settings. 

In Section \ref{sec:eb}, we consider the problem of estimating negative multinomial parameter vectors under the standardized squared error loss in the general case where sample sizes, lengths of observation vectors, and weights in the loss function may all be unbalanced. 
First, we generalize Theorem 1 of Hamura and Kubokawa (2020b) to this unbalanced case and also obtain another general sufficient condition for a general shrinkage estimator to dominate the UMVU estimator. 
Then, using the method of maximum likelihood, a new empirical Bayes estimator is derived which has a simple form as well as improves on the UMVU estimator. 
Finally, we present still another dominance condition, which is applicable specifically to empirical Bayes estimators including those based on the method of moments. 

In Section \ref{sec:hb}, we consider the practically important problem of estimating the joint predictive density of several independent multinomial tables under the Kullback-Leibler divergence. 
The distribution of any one of them is specified by a set of negative multinomial probability vectors, with each cell probability given by %
the product of the corresponding elements of the vectors. 
The setting we consider is quite general in that two tables may be related through a set of common overlapping probability vectors. 
Two simple special cases are the prediction problems for independent multinomial vectors and for a single multinomial table. 
We show that the Bayesian predictive density with respect to the Jeffreys prior is dominated by that with respect to a generalization of the shrinkage prior considered by Hamura and Kubokawa (2020b) under suitable conditions. 
Whereas Komaki (2012) investigated asymptotic properties of Bayesian predictive densities for future multinomial observations based on current multinomial observations, the sample space is not a finite set in our setting and we investigate finite sample properties of Bayesian predictive densities. 
Although Hamura and Kubokawa (2019a) considered Bayesian predictive densities for a negative binomial model, where a future observation also is negative binomial and can take on an infinite number of values, they did not treat the problem of estimating the joint predictive density of multiple negative binomial observations. 

In Section \ref{sec:sim}, simple and illustrative simulation studies are performed. 
In Section \ref{subsec:sim_eb}, our proposed empirical Bayes estimator and the UMVU estimator given in Section \ref{sec:eb} are compared. 
In Section \ref{subsec:sim_hb}, the Bayesian predictive densities given in Section \ref{sec:hb} are compared. 

In Section \ref{sec:discussion}, predictive density estimation for the negative multinomial distribution is discussed. 
Although no dominance conditions are obtained, generalizing Theorem 2.1 of Hamura and Kubokawa (2019a), we derive two kinds of identities which relate prediction to estimation in the negative multinomial case. 
In particular, the risk function of an arbitrary Bayesian predictive density under the Kullback-Leibler divergence is expressed using the risk functions of an infinite number of corresponding Bayes estimators under a weighted version of Stein's loss.

\section{Empirical Bayes Point Estimation}
\label{sec:eb}
Let $N \in \mathbb{N} = \{ 1, 2, \dotsc \} $, $m_1 , \dots , m_N \in \mathbb{N}$, and $r_1 , \dots , r_N > 0$. 
For $\nu = 1, \dots , N$, let $\p _{\nu } = ( p_{i, \nu } )_{i = 1}^{m_{\nu }} \in D_{m_{\nu }} = \{ ( \mathring{p} _1 , \dots , \mathring{p} _{m_{\nu }} )^{\top } | \mathring{p} _1 , \dots , \mathring{p} _{m_{\nu }} > 0, \, \sum_{i = 1}^{m_{\nu }} \mathring{p} _i < 1 \} $ and let $p_{0, \nu } = 1 - p_{\cdot , \nu } = 1 - \sum_{i = 1}^{m_{\nu }} p_{i, \nu }$. 
Let $\X _1 , \dots , \X _N$ be independent negative multinomial variables such that for each $\nu = 1, \dots , N$, the probability mass function of $\X _{\nu }$ is given by 
\begin{align}
{\Ga \big( r_{\nu } + \sum_{i = 1}^{m_{\nu }} x_{i, \nu } \big) \over \Ga ( r_{\nu } ) \prod_{i = 1}^{m_{\nu }} x_{i, \nu } !} {p_{0, \nu }}^{r_{\nu }} \prod_{i = 1}^{m_{\nu }} {p_{i, \nu }}^{x_{i, \nu }} \non %
\end{align}
for $\x _{\nu } = ( x_{i, \nu } )_{i = 1}^{m_{\nu }} \in {\mathbb{N} _0}^{m_{\nu }}$, where $\mathbb{N} _0 = \{ 0, 1, 2, \dotsc \} $. 
As pointed out by Hamura and Kubokawa (2020b), $m_1 , \dots , m_N$ may be different for example when we consider marginal distributions of negative multinomial vectors of the same length. 
For some basic properties of the negative multinomial distribution, see Sibuya, Yoshimura, and Shimizu (1964) and Tsui (1986b). 

Now we assume that all the elements of $\p = ( \p _{\nu } )_{\nu = 1, \dots , N} \in D = D_{m_1} \times \dots \times D_{m_N}$ are unknown and consider the problem of estimating $\p $ on the basis of the minimal and complete sufficient statistic $\X = ( \X _{\nu } )_{\nu = 1, \dots , N} = (( X _{i, \nu } )_{i = 1}^{m_{\nu }} )_{\nu = 1, \dots , N}$ under the standardized squared loss function given by 
\begin{align}
L_{n, \c } ( \d , \p ) &= \sum_{\nu = 1}^{n} \sum_{i = 1}^{m_{\nu }} c_{i, \nu } {( d_{i, \nu } - p_{i, \nu } )^2 \over p_{i, \nu }} \label{eq:loss_SS} 
\end{align}
for $\d = (( d_{i, \nu } )_{i = 1}^{m_{\nu }} )_{\nu = 1, \dots , N} \in {\mathbb{R}}^{m_1} \times \dots \times {\mathbb{R}}^{m_N}$, where $n \in \{ 1, \dots , N \} $ and $\c = (( c_{i, \nu } )_{i = 1}^{m_{\nu }} )_{\nu = 1, \dots , N} \in [0, \infty )^{m_1} \times \dots \times [0, \infty )^{m_N}$. 

For $\nu = 1, \dots , N$, let $X_{\cdot , \nu } = \sum_{i = 1}^{m_{\nu }} X_{i, \nu }$. 
Then the UMVU estimator of $\p $ is $\pbh ^{\rm{U}} = (( \ph _{i, \nu }^{\rm{U}} )_{i = 1}^{m_{\nu }} )_{\nu = 1, \dots , N}$, where 
\begin{align}
\ph _{i, \nu }^{\rm{U}} &= {X_{i, \nu } \over r_{\nu } + X_{\cdot , \nu } - 1} %
\label{eq:UMVU} 
\end{align}
for $i = 1, \dots , m_{\nu }$ for $\nu = 1, \dots , N$. 
(We write $0 / 0 = 0$.) 
We first derive a general sufficient condition for the shrinkage estimator 
\begin{align}
\pbh ^{( \bde )} &= (( \ph _{i, \nu }^{( \bde )} )_{i = 1}^{m_{\nu }} )_{\nu = 1, \dots , N} = \Big( \Big( {X_{i, \nu } \over r_{\nu } + X_{\cdot , \nu } - 1 + \de _{\nu } ( X_{\cdot , \cdot } )} \Big) _{i = 1}^{m_{\nu }} \Big) _{\nu = 1, \dots , N} \label{eq:shrinkage_estimator} 
\end{align}
to dominate $\pbh ^{\rm{U}}$, where $\bde = ( \de _{\nu } )_{\nu = 1}^{N} \colon \mathbb{N} _0 \to (0, \infty )^N$ and $X_{\cdot , \cdot } = \sum_{\nu = 1}^{N} X_{\cdot , \nu } = \sum_{\nu = 1}^{N} \sum_{i = 1}^{m_{\nu }} X_{i, \nu }$. 
For notational simplicity, let $\underline{r} = \min_{1 \le \nu \le n} r_{\nu }$ and $\overline{r} = \max_{1 \le \nu \le n} r_{\nu }$. 
For $\nu = 1, \dots , N$, let $c_{\cdot , \nu } = \sum_{i = 1}^{m_{\nu }} c_{i, \nu }$. 
Let $\underline{c_{\cdot }} = \min_{1 \le \nu \le n} c_{\cdot , \nu }$ and $\overline{\overline{c}} = \max_{1 \le \nu \le n} \max_{1 \le i \le m_{\nu }} c_{i, \nu }$. 
Finally, let $\underline{\de } (x) = \min_{1 \le \nu \le n} \de _{\nu } (x)$ and $\overline{\de } (x) = \max_{1 \le \nu \le n} \de _{\nu } (x)$ for $x \in \mathbb{N} _0$ and let $\rho = \inf_{x \in \mathbb{N} \setminus \{ 1 \} } \underline{\de } (x) / \overline{\de } (x) \in [0, 1]$. 
\begin{thm}
\label{thm:EB} 
Assume that $r_{\nu } \ge 5 / 2$ for all $\nu = 1, \dots , n$ with $c_{\cdot , \nu } > 0$ and that $0 < 3 \overline{\overline{c}} \le \underline{c_{\cdot }}$. %
Suppose that for all $\nu = 1, \dots , n$ such that $c_{\cdot , \nu } > 0$ and for all $x \in \mathbb{N}$, we have 
\begin{align}
x \de _{\nu } (x) \le (x + 1) \de _{\nu } (x + 1) \text{.} \label{eq:assumption_EB_0} 
\end{align}
Suppose further that for all $x \in \mathbb{N}$, one of the following two conditions are satisfied: 
\begin{enumerate}
\item[{\rm{(i)}}]
\begin{itemize}
\item
$\overline{\overline{c}} \overline{\de } (x + 1) \le 2 ( \underline{r} / \overline{r} )^2 ( \underline{c_{\cdot }} - 3 \overline{\overline{c}} ) \rho $ implies 
\begin{align}
\Big\{ 2 \Big( {\underline{r} \over \overline{r}} \Big) ^2 ( \underline{c_{\cdot }} - 3 \overline{\overline{c}} ) \rho - \underline{c_{\cdot }} \Big\} \overline{\de } (x + 1) + 2 \underline{r} \Big( {\underline{r} \over \overline{r}} \Big) ^2 ( \underline{c_{\cdot }} - 3 \overline{\overline{c}} ) \rho \ge 0 \quad \text{and} \label{eq:assumption_EB_1_1} 
\end{align}
\item
$\overline{\overline{c}} \overline{\de } (x + 1) > 2 ( \underline{r} / \overline{r} )^2 ( \underline{c_{\cdot }} - 3 \overline{\overline{c}} ) \rho $ implies 
\begin{align}
n \Big[ \Big\{ 2 \Big( {\underline{r} \over \overline{r}} \Big) ^2 ( \underline{c_{\cdot }} - 3 \overline{\overline{c}} ) \rho - \underline{c_{\cdot }} \Big\} \overline{\de } (x + 1) + 2 \underline{r} \Big( {\underline{r} \over \overline{r}} \Big) ^2 ( \underline{c_{\cdot }} - 3 \overline{\overline{c}} ) \rho \Big] \ge x \Big\{ \overline{\overline{c}} \overline{\de } (x + 1) - 2 \Big( {\underline{r} \over \overline{r}} \Big) ^2 ( \underline{c_{\cdot }} - 3 \overline{\overline{c}} ) \rho \Big\} \text{.} \label{eq:assumption_EB_1_2} 
\end{align}
\end{itemize}
\item[{\rm{(ii)}}]
\begin{itemize}
\item
$\overline{\overline{c}} \overline{\de } (x + 1) \le 2 ( \underline{c_{\cdot }} - 3 \overline{\overline{c}} ) \rho $ implies 
\begin{align}
2 ( \underline{c_{\cdot }} - 3 \overline{\overline{c}} ) \rho - ( \underline{c_{\cdot }} - \underline{r} \overline{\overline{c}} ) \ge 0 \quad \text{and} \label{eq:assumption_EB_2_1} 
\end{align}
\item
$\overline{\overline{c}} \overline{\de } (x + 1) > 2 ( \underline{c_{\cdot }} - 3 \overline{\overline{c}} ) \rho $ implies 
\begin{align}
n \{ 2 ( \underline{c_{\cdot }} - 3 \overline{\overline{c}} ) \rho - ( \underline{c_{\cdot }} - \underline{r} \overline{\overline{c}} ) \} \overline{\de } (x + 1) \ge \Big( \sum_{\nu = 1}^{n} r_{\nu } + x \Big) \{ \overline{\overline{c}} \overline{\de } (x + 1) - 2 ( \underline{c_{\cdot }} - 3 \overline{\overline{c}} ) \rho \} \text{.} \label{eq:assumption_EB_2_2} 
\end{align}
\end{itemize}
\end{enumerate}
Then the shrinkage estimator $\pbh ^{( \bde )}$ given in (\ref{eq:shrinkage_estimator}) dominates the UMVU estimator $\pbh ^{\rm{U}}$ given by (\ref{eq:UMVU}) under the standardized squared loss (\ref{eq:loss_SS}). 
\end{thm}

Part (i) of Theorem \ref{thm:EB} is a generalization of Theorem 1 of Hamura and Kubokawa (2020b), who further obtained simpler conditions in specific cases. 
On the other hand, part (ii) is another result of this paper. 
It is worth noting that under the setting of Theorem \ref{thm:EB}, there may exist $\nu = 1, \dots , n$ such that $c_{i, \nu } = 0 < c_{i' , \nu }$ for some $i, i' = 1, \dots , m_{\nu }$. 

Next, we derive an empirical Bayes estimator based on the method of maximum likelihood. 
Consider the conjugate Dirichlet prior distribution 
\begin{align}
\prod_{\nu = 1}^{N} {\rm{Dir}}_{m_{\nu }} ( \p _{\nu } | \at _{\nu } v, \j ^{( m_{\nu } )} ) = \prod_{\nu = 1}^{N} \Big\{ {\Ga ( \at _{\nu } v + m_{\nu } ) \over \Ga ( \at _{\nu } v)} {p_{0, \nu }}^{\at _{\nu } v - 1} \Big\} \text{,} \non 
\end{align}
where $v \in (0, \infty )$ and where $\at _{\nu } \in (0, \infty )$ and $\j ^{( m_{\nu } )} = (1, \dots , 1)^{\top } \in \mathbb{R} ^{m_{\nu }}$ for $\nu = 1, \dots , N$. 
It corresponds to the Bayes estimator 
\begin{align}
\Big( \Big( {X_{i, \nu } \over r_{\nu } + X_{\cdot , \nu } - 1 + \at _{\nu } v + m_{\nu }} \Big) _{i = 1}^{m_{\nu }} \Big) _{\nu = 1, \dots , N} \non 
\end{align}
of $\p $. 
On the other hand, since the maximum likelihood estimator and the prior mean of $p_{0, \nu }$ is $r_{\nu } / ( r_{\nu } + X_{\cdot , \nu } )$ and $\at _{\nu } v / ( \at _{\nu } v + m_{\nu } )$ for $\nu = 1, \dots , N$, a reasonable estimator of $v$ would be 
\begin{align}
{1 \over X_{\cdot , \cdot }} \sum_{\nu = 1}^{N} {m_{\nu } r_{\nu } \over \at _{\nu }} \text{.} \non 
\end{align}
Thus, we obtain the empirical Bayes estimator 
\begin{align}
\pbh ^{( \abt )} &= \Big( \Big( {X_{i, \nu } \over r_{\nu } + X_{\cdot , \nu } - 1 + \de _{\nu }^{( \abt )} ( X_{\cdot , \cdot } )} \Big) _{i = 1}^{m_{\nu }} \Big) _{\nu = 1, \dots , N} \text{,} \label{eq:EB_ML} 
\end{align}
where $\abt = ( \at _{\nu } )_{\nu = 1}^{N}$ and where 
\begin{align}
\de _{\nu }^{( \abt )} ( X_{\cdot , \cdot } ) &= m_{\nu } + {\at _{\nu } \over X_{\cdot , \cdot }} \sum_{{\nu }' = 1}^{N} {m_{{\nu }'} r_{{\nu }'} \over \at _{{\nu }'}} \non 
\end{align}
if $X_{\cdot , \cdot } \ge 1$ while  $\de _{\nu }^{( \abt )} (0) \in (0, \infty )$ for $\nu = 1, \dots , N$. 
This estimator was not considered by Hamura and Kubokawa (2020b). 
It is of the form (\ref{eq:shrinkage_estimator}) and clearly satisfies condition (\ref{eq:assumption_EB_0}). 
Whether the other conditions hold or not depends on the choice of the hyperparameter $\abt $. 
For example, 
\begin{align}
\rho &= \begin{cases} \displaystyle \inf_{x \in \mathbb{N} \setminus \{ 1 \} } {( \min_{1 \le \nu \le n} m_{\nu } ) \big( 1 + \sum_{{\nu }' = 1}^{N} r_{{\nu }'} / x \big) \over ( \max_{1 \le \nu \le n} m_{\nu } ) \big( 1 + \sum_{{\nu }' = 1}^{N} r_{{\nu }'} / x \big) } = {\min_{1 \le \nu \le n} m_{\nu } \over \max_{1 \le \nu \le n} m_{\nu }} \text{,} & \text{if $\abt = ( m_{\nu } )_{\nu = 1}^{N}$} \text{,} \\ \displaystyle \inf_{x \in \mathbb{N} \setminus \{ 1 \} } {\min_{1 \le \nu \le n} m_{\nu } + \sum_{{\nu }' = 1}^{N} m_{{\nu }'} r_{{\nu }'} / x \over \max_{1 \le \nu \le n} m_{\nu } + \sum_{{\nu }' = 1}^{N} m_{{\nu }'} r_{{\nu }'} / x} = {\min_{1 \le \nu \le n} m_{\nu } \over \max_{1 \le \nu \le n} m_{\nu }} \text{,} & \text{if $\abt = \j ^{(N)}$} \text{,} \\ \displaystyle \inf_{x \in \mathbb{N} \setminus \{ 1 \} } {\min_{1 \le \nu \le n} \big( m_{\nu } + r_{\nu } \sum_{{\nu }' = 1}^{N} m_{{\nu }'} / x \big) \over \max_{1 \le \nu \le n} \big( m_{\nu } + r_{\nu } \sum_{{\nu }' = 1}^{N} m_{{\nu }'} / x \big) } \text{,} & \text{if $\abt = ( r_{\nu } )_{\nu = 1}^{N}$} \text{,} \end{cases} \non 
\end{align}
where $\j ^{(N)} = (1, \dots , 1)^{\top } \in \mathbb{R} ^N$. 

There are other empirical Bayes estimators. 
For example, since the prior mean of $E[ \X_{\cdot , \cdot } ] = \sum_{\nu = 1}^{N} \sum_{i = 1}^{m_{\nu }} r_{\nu } p_{i, \nu } / p_{0, \nu }$ is $\sum_{\nu = 1}^{N} \sum_{i = 1}^{m_{\nu }} r_{\nu } / (v - 1) = \sum_{\nu = 1}^{N} m_{\nu } r_{\nu } / (v - 1)$ when $\at _{\nu } = 1$ and $v > 1$ for all $\nu = 1, \dots , N$, one estimator of $v$ based on the method of moments would be 
\begin{align}
1 + {1 \over X_{\cdot , \cdot }} \sum_{\nu = 1}^{N} m_{\nu } r_{\nu } \text{.} \non 
\end{align}
We could also use $1 + \big( \sum_{\nu = 1}^{N} \sum_{i = 1}^{m_{\nu }} r_{\nu } \ct _{i, \nu } \big) / \sum_{\nu = 1}^{N} \sum_{i = 1}^{m_{\nu }} \ct _{i, \nu } X_{i, \nu }$ 
for $(( \ct _{i, \nu } )_{i = 1}^{m_{\nu }} )_{\nu = 1, \dots , N} \in (0, \infty )^{m_1} \times \dots \times (0, \infty )^{m_N}$. 
More generally, we consider the shrinkage estimator 
\begin{align}
\pbh ^{( \bbt , \cbt )} &= (( \ph _{i, \nu }^{( \bbt , \cbt )} )_{i = 1}^{m_{\nu }} )_{\nu = 1, \dots , N} = \Big( \Big( {X_{i, \nu } \over r_{\nu } + X_{\cdot , \nu } - 1 + \bt _{\nu } + 1 / \Xt ^{( \cbt ^{( \nu )} )}} \Big) _{i = 1}^{m_{\nu }} \Big) _{\nu = 1, \dots , N} \text{,} \label{eq:EB_affine} 
\end{align}
where $\bbt = ( \bt _{\nu } )_{\nu = 1}^{N} \in (0, \infty )^N$ and $\cbt = ( \cbt ^{( \nu )} )_{\nu = 1}^{N} = ((( \ct _{i, {\nu }'}^{( \nu )} )_{i = 1}^{m_{{\nu }'}} )_{{\nu }' = 1, \dots , N} )_{\nu  = 1}^{N} \in ((0, \infty )^{m_1} \times \dots \times (0, \infty )^{m_N} )^N$ and where $\Xt ^{( \cbt ^{( \nu )} )} = \sum_{{\nu }' = 1}^{N} \sum_{i = 1}^{m_{{\nu }'}} \ct _{i, {\nu }'}^{( \nu )} X_{i, {\nu }'}$ for $\nu = 1, \dots , N$. 

\begin{thm}
\label{thm:EB_affine} 
Under Assumption \ref{assumption:tEB_affine} given in the Appendix, the shrinkage estimator $\pbh ^{( \bbt , \cbt )}$ given in (\ref{eq:EB_affine}) dominates the UMVU estimator $\pbh ^{\rm{U}}$ given by (\ref{eq:UMVU}) under the standardized squared loss (\ref{eq:loss_SS}). 
\end{thm}

When $\Xt ^{( \cbt ^{(1)} )} = \dots = \Xt ^{( \cbt ^{(N)} )} = \ct X_{\cdot , \cdot }$, where $\ct \in (0, \infty )$, we have the following result. 

\begin{cor}
\label{cor:EB_affine} 
Assume that $\cbt ^{(1)} = \dots = \cbt ^{(N)} = ( \ct \j ^{( m_1 )} , \dots , \ct \j ^{( m_N )} )$. 
Then, under Assumption \ref{assumption:cEB_affine} given in the Appendix, $\pbh ^{( \bbt , \cbt )}$ dominates $\pbh ^{\rm{U}}$ under the loss (\ref{eq:loss_SS}). 
\end{cor}

In Corollary \ref{cor:EB_affine}, it is not necessarily assumed as in Theorem \ref{thm:EB} that $r_{\nu } \ge 5 / 2$ for all $\nu = 1, \dots , n$ with $c_{\cdot , \nu } > 0$. 
Moreover, for the balanced case with $r_1 \ge 1$, another dominance condition can be obtained by modifying the proof of Theorem \ref{thm:EB_affine} given in the Appendix. 
See Remark \ref{rem:1r} for details. 

Finally, in order to estimate $\p $, we could also use the hierarchical shrinkage prior introduced by Hamura and Kubokawa (2020b) or its generalization. 
However, since they considered essentially the same hierarchical Bayes estimator and gave important methods of evaluating the risk function, we do not discuss the approach further. 
The usefulness of hierarchical Bayes procedures will be shown in the next section. 

\section{Hierarchical Bayes Predictive Density Estimation}
\label{sec:hb}
In this section, we consider predictive density estimation for the multinomial distribution. 
Let $L \in \mathbb{N}$ and $d^{(1)} , \dots , d^{(L)} \in \{ 1, \dots , N \} $. 
For $\la = 1, \dots , L$, let $\nu _{1}^{( \la )} , \dots , \nu _{d^{( \la )}}^{( \la )} \in \mathbb{N}$ be such that $1 \le \nu _{1}^{( \la )} < \dots < \nu _{d^{( \la )}}^{( \la )} \le N$ and let $I_{0}^{( \la )} = \{ 0, 1, \dots , m_{\nu _{1}^{( \la )}} \} \times \dots \times \{ 0, 1, \dots , m_{\nu _{d^{( \la )}}^{( \la )}} \} $ and $\mathcal{W} ^{( \la )} = \big\{ ( \mathring{w} _{\i } )_{\i \in I_{0}^{( \la )}} \big| \mathring{w} _{\i } \in \mathbb{N} _0 \quad \text{for all $\i \in I_{0}^{( \la )}$} \quad \text{and} \quad  \sum_{\i \in I_{0}^{( \la )}} \mathring{w} _{\i } = l^{( \la )} \big\} $. 
Now let $l^{(1)} , \dots , l^{(L)} \in \mathbb{N}$ and let $\W ^{(1)} , \dots , \W ^{(L)}$ be independent multinomial variables such that for $\la = 1, \dots , L$, the probability mass function of $\W ^{( \la )}$ is given by 
\begin{align}
f_{\la } ( \w ^{( \la )} | \p ) &= {l^{( \la )} ! \over \prod_{\i \in I_{0}^{( \la )}} w_{\i }^{( \la )} !} \prod_{\i = ( i_h )_{h = 1}^{d^{( \la )}} \in I_{0}^{( \la )}} \Big\{ \prod_{h = 1}^{d^{( \la )}} p_{i_h , \nu _{h}^{( \la )}} \Big\} ^{w_{\i }^{( \la )}} \non %
\end{align}
for $\w ^{( \la )} = ( w _{\i }^{( \la )} )_{\i \in I_{0}^{( \la )}} \in \mathcal{W} ^{( \la )}$. %
We consider the problem of estimating the joint probability mass of $\W ^{(1)} , \dots , \W ^{(L)}$, namely $f( \w | \p ) = \prod_{\la = 1}^{L} f_{\la } ( \w ^{( \la )} | \p )$, $\w = ( \w ^{( \la )} )_{\la = 1, \dots , L} \in \mathcal{W} = \mathcal{W} ^{(1)} \times \dots \times \mathcal{W} ^{(L)}$, on the basis of $\X $ given in the previous section under the Kullback-Leibler divergence. 
The risk function of a predictive mass $\fh ( \cdot ; \X )$ is given by 
\begin{align}
E \Big[ \log {f( \W | \p ) \over \fh ( \W ; \X )} \Big] \text{,} \non 
\end{align}
where $\W = ( \W ^{( \la )} )_{\la = 1, \dots , L} = (( W _{\i }^{( \la )} )_{\i \in I_{0}^{( \la )}} )_{\la = 1, \dots , L}$. 

As noted in Remark 2.2 of Hamura and Kubokawa (2019a), defining a natural plug-in predictive mass is not necessarily easy. 
Therefore, in this section, we seek a good %
Bayesian predictive mass. 
As shown by Aitchison (1975), the Bayesian predictive mass $\fh ^{( \pi )} ( \cdot ; \X )$ associated with a prior $\p \sim \pi ( \p )$ is given by 
\begin{align}
\fh ^{( \pi )} ( \w ; \X ) &= E_{\pi } [ f( \w | \p ) | \X ] \text{.} \label{eq:Bayesian_predictive_mass} 
\end{align}

We first consider the natural conjugate Dirichlet distribution with density 
\begin{align}
\pi _{\a _0 , \a } ( \p ) %
&\propto \prod_{\nu = 1}^{N} \Big( %
{p_{0, \nu }}^{a_{0, \nu } - 1} \prod_{i = 1}^{m_{\nu }} {p_{i, \nu }}^{a_{i, \nu } - 1} \Big) \text{,} \label{eq:prior_Dir} 
\end{align}
where %
$\a _0 = ( a_{0, \nu } )_{\nu = 1}^{N} \in \mathbb{R} ^N$, $\a = ( \a _{\nu } )_{\nu = 1, \dots , N} = (( a_{i, \nu } )_{i = 1}^{m_{\nu }} )_{\nu = 1, \dots , N} \in (0, \infty )^{m_1} \times \dots \times (0, \infty )^{m_N}$, and $a_{\cdot , \nu } = \sum_{i = 1}^{m_{\nu }} a_{i, \nu }$ for $\nu = 1, \dots , N$. 
The Jeffreys prior is a special case of the Dirichlet prior. 

\begin{lem}
\label{lem:Jeff} 
The Dirichlet prior (\ref{eq:prior_Dir}) with $\a _0 = ((1 - m_{\nu } ) / 2)_{\nu = 1}^{N}$ and $\a = ( \j ^{( m_{\nu } )} / 2)_{\nu = 1, \dots , N}$ is the Jeffreys prior. 
\end{lem}

Next we consider the following conjugate shrinkage prior. 
Let 
\begin{align}
\pi _{\al , \be , \bga , \a _0 , \a } ( \p ) &= \int_{0}^{\infty } u^{\al - 1} e^{- \be u} \Big\{ \prod_{\nu = 1}^{N} \Big( {p_{0, \nu }}^{\ga _{\nu } u + a_{0, \nu } - 1} \prod_{i = 1}^{m_{\nu }} {p_{i, \nu }}^{a_{i, \nu } - 1} \Big) \Big\} du \text{,} \label{eq:prior_shrinkage} 
\end{align}
where $\al > 0$, $\be > 0$, and $\bga = ( \ga _{\nu } )_{\nu = 1}^{N} \in (0, \infty )^N$. 
This shrinkage prior is based on that of Section 3 of Hamura and Kubokawa (2020b) and is a slightly simplified version of the one mentioned in the discussion of their papar. 

Under the prior (\ref{eq:prior_Dir}), the posterior distribution of $\p $ given $\X = \x $ is proper for all $\x \in {\mathbb{N} _0}^{m_1} \times \dots \times {\mathbb{N} _0}^{m_N}$ if and only if $r_{\nu } + a_{0, \nu } > 0$ for all $\nu = 1, \dots , N$. 
Also, this condition implies that the posterior under (\ref{eq:prior_shrinkage}) is proper, since we have assumed that $\be \neq 0$ for simplicity. 

In order to derive the Bayesian predictive mass with respect to (\ref{eq:prior_Dir}) and that with respect to (\ref{eq:prior_shrinkage}) in Proposition \ref{prp:mass_estimator}, we first rewrite $f( \w | \p )$. 
Let $S( \la ) = \{ \nu _{1}^{( \la )} , \dots , \nu _{d^{( \la )}}^{( \la )} \} $ for $\la = 1, \dots , L$. 
For $\nu = 1, \dots , N$, let $\La ( \nu ) = \{ \la \in \{ 1, \dots , L \} | \nu \in S( \la ) \} $ and, for $\la \in \La ( \nu )$, let $\{ h_{\nu }^{( \la )} \} = \{ h \in \{ 1, \dots , d^{( \la )} \} | \nu = \nu _{h}^{( \la )} \} $ and let, for $i = 0, 1, \dots , m_{\nu }$, $I_{0}^{( \la )} (i, \nu ) = \{ ( i_h )_{h = 1}^{d^{( \la )}} \in I_{0}^{( \la )} | i_{h_{\nu }^{( \la )}} = i \} $. 

\begin{lem}
\label{lem:likelihood} 
For any $\w = (( w_{\i }^{( \la )} )_{\i \in I_{0}^{( \la )}} )_{\la = 1, \dots , L} \in \mathcal{W}$, we have 
\begin{align}
f( \w | \p ) &= \Big\{ \prod_{\la = 1}^{L} {l^{( \la )} ! \over \prod_{\i \in I_{0}^{( \la )}} w_{\i }^{( \la )} !} \Big\} \prod_{\nu = 1}^{N} \prod_{i = 0}^{m_{\nu }} {p_{i, \nu }}^{\sum_{\la \in \La ( \nu )} \sum_{\i \in I_{0}^{( \la )} (i, \nu )} w_{\i }^{( \la )}} \text{.} \non 
\end{align}
\end{lem}

Let 
\begin{align}
C( \w ) &= \prod_{\la = 1}^{L} {l^{( \la )} ! \over \prod_{\i \in I_{0}^{( \la )}} w_{\i }^{( \la )} !} \non 
\end{align}
for $\w = (( w_{\i }^{( \la )} )_{\i \in I_{0}^{( \la )}} )_{\la = 1, \dots , L} \in \mathcal{W}$. 
For $(i, \nu ) \in \mathbb{N} _0 \times \{ 1, \dots , N \} $ with $i \le m_{\nu }$, let 
\begin{align}
s_{i, \nu } ( \w ) = \sum_{\la \in \La ( \nu )} \sum_{\i \in I_{0}^{( \la )} (i, \nu )} w_{\i }^{( \la )} \non 
\end{align}
for $\w = (( w_{\i }^{( \la )} )_{\i \in I_{0}^{( \la )}} )_{\la = 1, \dots , L} \in \mathcal{W}$. 
Using (\ref{eq:Bayesian_predictive_mass}) and Lemma \ref{lem:likelihood}, the following expressions for $\fh ^{( \pi _{\a _0 , \a } )} ( \cdot ; \X )$ and $\fh ^{( \pi _{\al , \be , \bga , \a _0 , \a } )} ( \cdot ; \X )$ are obtained. 

\begin{prp}
\label{prp:mass_estimator} 
Suppose that $r_{\nu } + a_{0, \nu } > 0$ for all $\nu = 1, \dots , N$. 
\begin{enumerate}
\item[{\rm{(i)}}]
The Bayesian predictive mass $\fh ^{( \pi _{\a _0 , \a } )} ( \cdot ; \X )$ is given by 
\begin{align}
\fh ^{( \pi _{\a _0 , \a } )} ( \w ; \X ) &= C( \w ) \frac{ \displaystyle \prod_{\nu = 1}^{N} {\Ga ( s_{0, \nu } ( \w ) + r_{\nu } + a_{0, \nu } ) \prod_{i = 1}^{m_{\nu }} \Ga ( s_{i, \nu } ( \w ) + X_{i, \nu } + a_i ) \over \Ga \big( \sum_{\la \in \La ( \nu )} l^{( \la )} + r_{\nu } + a_{0, \nu } + X_{\cdot , \nu } + a_{\cdot , \nu } \big) } }{ \displaystyle \prod_{\nu = 1}^{N} {\Ga ( r_{\nu } + a_{0, \nu } ) \prod_{i = 1}^{m_{\nu }} \Ga ( X_{i, \nu } + a_i ) \over \Ga ( r_{\nu } + a_{0, \nu } + X_{\cdot , \nu } + a_{\cdot , \nu } )} } \text{.} \non 
\end{align}
\item[{\rm{(ii)}}]
The Bayesian predictive mass $\fh ^{( \pi _{\al , \be , \bga , \a _0 , \a } )} ( \cdot ; \X )$ is given by 
\begin{align}
&\fh ^{( \pi _{\al , \be , \bga , \a _0 , \a } )} ( \w ; \X ) \non \\
&= C( \w ) \frac{ \displaystyle \int_{0}^{\infty } u^{\al - 1} e^{- \be u} \Big\{ \prod_{\nu = 1}^{N} {\Ga ( \ga _{\nu } u + s_{0, \nu } ( \w ) + r_{\nu } + a_{0, \nu } ) \prod_{i = 1}^{m_{\nu }} \Ga ( s_{i, \nu } ( \w ) + X_{i, \nu } + a_i ) \over \Ga \big( \ga _{\nu } u + \sum_{\la \in \La ( \nu )} l^{( \la )} + r_{\nu } + a_{0, \nu } + X_{\cdot , \nu } + a_{\cdot , \nu } \big) } \Big\} du }{ \displaystyle \int_{0}^{\infty } u^{\al - 1} e^{- \be u} \Big\{ \prod_{\nu = 1}^{N} {\Ga ( \ga _{\nu } u + r_{\nu } + a_{0, \nu } ) \prod_{i = 1}^{m_{\nu }} \Ga ( X_{i, \nu } + a_i ) \over \Ga ( \ga _{\nu } u + r_{\nu } + a_{0, \nu } + X_{\cdot , \nu } + a_{\cdot , \nu } )} \Big\} du } \text{.} \non 
\end{align}
\end{enumerate}
\end{prp}

We now compare the risk functions of $\fh ^{( \pi _{\a _0 , \a } )} ( \cdot ; \X )$ and $\fh ^{( \pi _{\al , \be , \bga , \a _0 , \a } )} ( \cdot ; \X )$. 

\begin{thm}
\label{thm:multin} 
Assume that $r_{\nu } + a_{0, \nu } > 0$ for all $\nu = 1, \dots , N$. 
Assume that $r_{\nu } \ge 1$ for all $\nu = 1, \dots , N$. 
Suppose that 
\begin{align}
\Big\{ {( \al + 1) \ga _{\nu } \over \be + \ga _{\nu }} - a_{\cdot , \nu } \Big\} ( r_{\nu } - 1) \le x_{\nu } \Big\{ - {( \al + 1) \ga _{\nu } \over \be + \ga _{\nu }} - \sum_{\la \in \La ( \nu )} l^{( \la )} - a_{0, \nu } \Big\} \label{eq:assumption_multin} 
\end{align}
for all $x_{\nu } \in \mathbb{N}$ for all $\nu = 1, \dots , N$. 
Then $\fh ^{( \pi _{\al , \be , \bga , \a _0 , \a } )} ( \cdot ; \X )$ dominates $\fh ^{( \pi _{\a _0 , \a } )} ( \cdot ; \X )$. 
\end{thm}

\begin{cor}
\label{cor:multin} 
If $1 \le r_{\nu } > ( m_{\nu } - 1) / 2 > \sum_{\la \in \La ( \nu )} l^{( \la )}$ for all $\nu = 1, \dots , N$, then the Bayesian predictive mass with respect to the Jeffreys prior, namely $\fh ^{( \pi _{\a _0 , \a } )} ( \cdot ; \X )$ with $\a _0 = ((1 - m_{\nu } ) / 2)_{\nu = 1}^{N}$ and $\a = ( \j ^{( m_{\nu } )} / 2)_{\nu = 1, \dots , N}$, is inadmissible and dominated by the Bayesian predictive mass $\fh ^{( \pi _{\al , \be , \bga , \a _0 , \a } )} ( \cdot ; \X )$ with $\a _0 = ((1 - m_{\nu } ) / 2)_{\nu = 1}^{N}$ and $\a = ( \j ^{( m_{\nu } )} / 2)_{\nu = 1, \dots , N}$ for some $\al > 0$, $\be > 0$, and $\bga \in (0, \infty )^N$. 
\end{cor}

\section{Simulation Studies}
\label{sec:sim}
\subsection{Simulation study for the model in Section \ref{sec:eb}}
\label{subsec:sim_eb} 
In this section, we investigate through simulation the numerical performance of the risk functions of point estimators of $\p $ under the standardized squared error loss given by (\ref{eq:loss_SS}). 
Although there are a number of conceivable unbalanced settings, for the sake of simplicity, we only consider some of the most uncomplicated cases. 
In particular, we set $n = N = 2$, $m_1 = m_2 = 7$, and $\c = ( \j ^{(7)} , \j ^{(7)} )$ and focus on the effect of $r_1$, $r_2$, and $\p $. 
As in the Poisson case (see, for example, Hamura and Kubokawa (2019b, 2020c)), although the dominance conditions given in Section \ref{sec:eb} tend to be restrictive and may not be satisfied especially when $r_1$ and $r_2$ are highly unbalanced, our proposed estimator turns out to perform well in such cases also. 

We compare the UMVU estimator $\hat{\p } ^{\rm{U}}$ given by (\ref{eq:UMVU}) and the empirical Bayes estimator $\pbh ^{( \abt )}$ given in (\ref{eq:EB_ML}) with $\abt = \j ^{(N)}$, namely 
\begin{align} 
\pbh ^{\rm{EB}} &= \Big( \Big( {X_{i, \nu } \over r_{\nu } + X_{\cdot , \nu } - 1 + 7 + 7 \sum_{{\nu }' = 1}^{2} r_{{\nu }'} / X_{\cdot , \cdot }} \Big) _{i = 1}^{7} \Big) _{\nu = 1, 2} \text{.} \non 
\end{align}
Let $\p _0 (0) = (1, 1, 1, 1, 1, 1, 1)^{\top } / 8$, $\p _0 (1) = (1, 1, 1, 1, 10, 10, 10)^{\top } / 44$, and $\p _0 (2) = (10, 10, 10, 10, 1, 1, 1)^{\top } / 44$. 
We consider the following cases: 
\begin{itemize}
\item[(i)]
Let $r_1 = r_2 = 12$ and let $\p _1 = \p _2 = (1 - \om ) \p _0 (0) + \om \p _0 (1)$ for $\om = 0, 1 / 5, \dots , 4 / 5, 1$. 
\item[(ii)]
Let $r_1 = r_2 = 12$ and let $\p _1 = (1 - \om ) \p _0 (0) + \om \p _0 (1)$ and $\p _2 = (1 - \om ) \p _0 (0) + \om \p _0 (2)$ for $\om = 0, 1 / 5, \dots , 4 / 5, 1$. 
\item[(iii)]
Let $r_1 = 8$ and $r_2 = 16$ and let $\p _1 = \p _2 = (1 - \om ) \p _0 (0) + \om \p _0 (1)$ for $\om = 0, 1 / 5, \dots , 4 / 5, 1$. 
\item[(iv)]
Let $r_1 = 8$ and $r_2 = 16$ and let $\p _1 = (1 - \om ) \p _0 (0) + \om \p _0 (1)$ and $\p _2 = (1 - \om ) \p _0 (0) + \om \p _0 (2)$ for $\om = 0, 1 / 5, \dots , 4 / 5, 1$. 
\end{itemize}
In Cases (i) and (ii), $r_1$ and $r_2$ are balanced. 
On the other hand, they are highly unbalanced in Cases (iii) and (iv). 
The parameter vectors $\p _1$ and $\p _2$ are identical for all $\om = 0, 1 / 5, \dots , 4 / 5, 1$ in Cases (i) and (iii) and distinct for $\om = 1 / 5, \dots , 4 / 5, 1$ in Cases (ii) and (iv). 
We obtain approximated values of the risk functions of $\hat{\p } ^{\rm{U}}$ and $\pbh ^{\rm{EB}}$ by simulation with $100,000$ replications.

\begin{figure}
\centering
\includegraphics[width = 16cm]{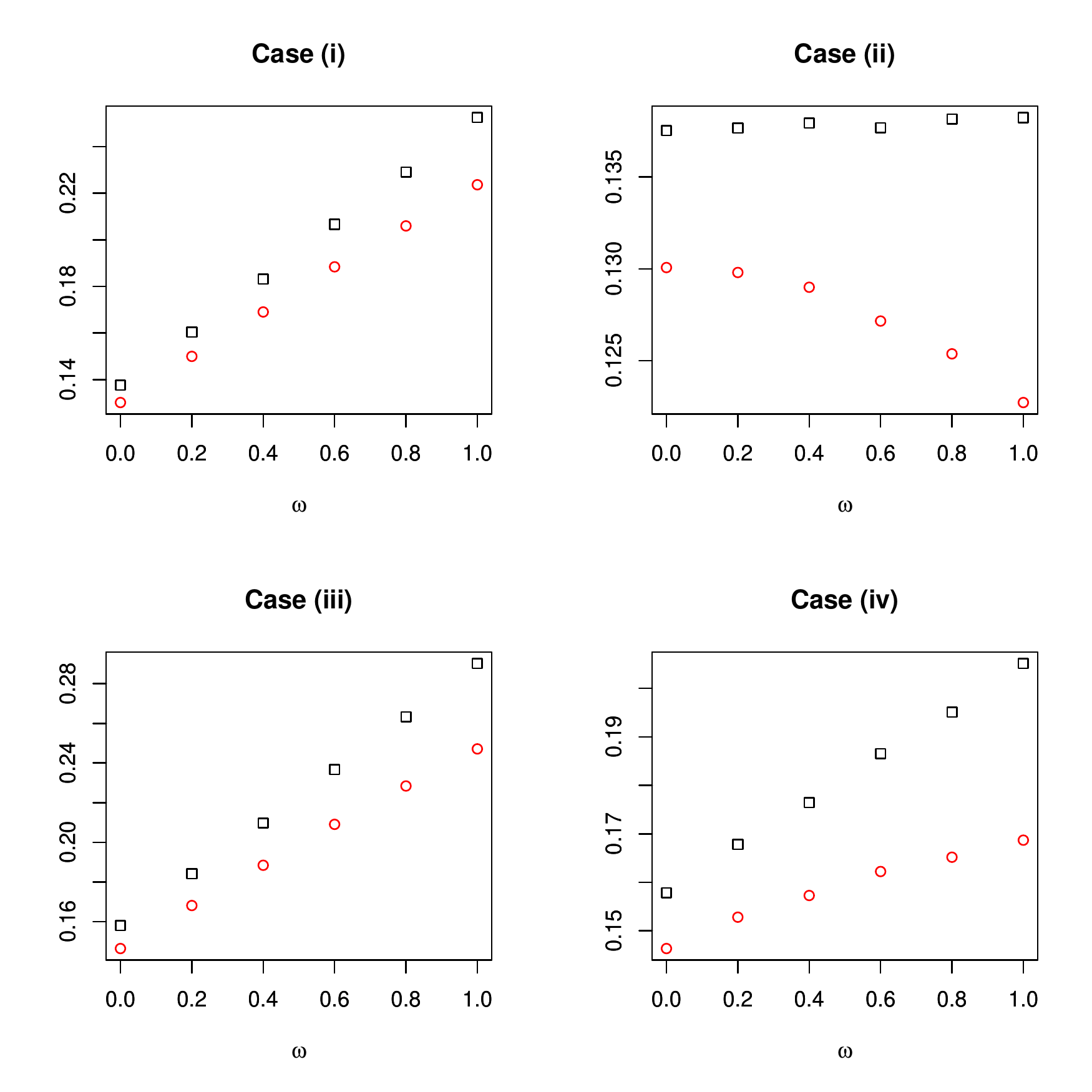}
\caption{Risks of the estimators $\hat{\p } ^{\rm{U}}$ and $\pbh ^{\rm{EB}}$ for $\om = 0, 1 / 5, \dots , 4 / 5, 1$ in Cases (i), (ii), (iii), and (iv). 
The black squares and red circles correspond to $\hat{\p } ^{\rm{U}}$ and $\pbh ^{\rm{EB}}$, respectively. }
\label{fig:eb} 
\end{figure}

The results are illustrated in Figure \ref{fig:eb}. 
It seems that $\pbh ^{\rm{EB}}$ dominates $\hat{\p } ^{\rm{U}}$ in every case. 
In Cases (i) and (iii), both $\hat{\p } ^{\rm{U}}$ and $\pbh ^{\rm{EB}}$ have large values of risks for large $\om $. 
In Case (ii), the risk values of $\hat{\p } ^{\rm{U}}$ are almost the same while those of $\pbh ^{\rm{EB}}$ are small for large $\om $. 
On the other hand, in Case (iv), where the amount of information from $\X _2$ is much larger than the amount of information from $\X _1$, the results are similar to those in Cases (i) and (iii). 
Overall, the risk values are smaller in Cases (i) and (ii) than in Cases (iii) and (iv) and larger in Cases (i) and (iii) than in Cases (ii) and (iv). 

\subsection{Simulation study for the model in Section \ref{sec:hb}}
\label{subsec:sim_hb} 
This section corresponds to Section \ref{sec:hb}. 
As in Section \ref{subsec:sim_eb}, we focus on simple cases and in particular consider low-dimensional settings for computational convenience. 
We set $N = 2$, $m_1 = m_2 = 3$, $L = 2$, $d^{(1)} = 1$, $d^{(2)} = 2$, $\nu _{1}^{(1)} = 1$, $\nu _{1}^{(2)} = 1$, $\nu _{2}^{(2)} = 2$, and $l^{(1)} = l^{(2)} = 1$. 
We note that $\p _1$ is related to both the vector $\W ^{(1)}$ and the matrix $\W ^{(2)}$. 
We investigate through simulation the numerical performance of the risk functions of $\fh ^{( \pi _{\a _0 , \a } )} ( \cdot ; \X )$ given in part (i) of Proposition \ref{prp:mass_estimator} and $\fh ^{( \pi _{\al , \be , \bga , \a _0 , \a } )} ( \cdot ; \X )$ given in part (ii) of Proposition \ref{prp:mass_estimator}; more specifically, we set $\a _0 = (- 1, -1)^{\top }$, $\a = ( \j ^{(3)} / 2, \j ^{(3)} / 2)$, $\al = 1$, $\be = 1$, and $\ga = (1, 1)^{\top }$ and compare the Bayesian predictive mass with respect to the Jeffreys prior, namely $\fh ^{\rm{J}} ( \cdot ; \X ) = \fh ^{( \pi _{(- 1, -1)^{\top } , ( \j ^{(3)} / 2, \j ^{(3)} / 2)} )} ( \cdot ; \X )$, and the Bayesian predictive mass $\fh ^{\rm{HB}} ( \cdot ; \X ) = \fh ^{( \pi _{1, 1, (1, 1)^{\top } , (- 1, -1)^{\top } , ( \j ^{(3)} / 2, \j ^{(3)} / 2)} )} ( \cdot ; \X )$. 
Let $\p (0) = ( (1, 1, 1)^{\top } / 4, (1, 1, 1)^{\top } / 4)$, $\p (1) = ( (1, 1, 2)^{\top } / 6, (1, 1, 2)^{\top } / 6)$, and $\p (2) = ( (1, 1, 2)^{\top } / 6, (2, 2, 1)^{\top } / 6)$. 
For each $\p = \p (0), \p (1), \p (2)$, we consider the following cases: (I) $r_1 = r_2 = 5$; (II) $r_1 = 4$ and $r_2 = 6$; (III) $r_1 = 6$ and $r_2 = 4$. 

We obtain approximated values of the risk functions of $\fh ^{\rm{J}} ( \cdot ; \X )$ and $\fh ^{\rm{HB}} ( \cdot ; \X )$ by simulation with $1,000$ replications. 
The Bayesian predictive mass $\fh ^{\rm{J}} ( \cdot ; \X )$ is computed by generating $2,000$ independent posterior samples while $\fh ^{\rm{HB}} ( \cdot ; \X )$ is computed based on a Gibbs sampler by generating $20,000$ approximate posterior samples after discarding the first $10,000$ samples. 
The percentage relative improvement in average loss (PRIAL) of $\fh ^{\rm{HB}} ( \cdot ; \X )$ over $\fh ^{\rm{J}} ( \cdot ; \X )$ is defined by 
\begin{align}
{\rm PRIAL} = 100 \Big\{ E \Big[ \log {f( \W | \p ) \over \fh ^{\rm{J}} ( \W ; \X )} \Big] - E \Big[ \log {f( \W | \p ) \over \fh ^{\rm{HB}} ( \W ; \X )} \Big] \Big\} / E \Big[ \log {f( \W | \p ) \over \fh ^{\rm{J}} ( \cdot ; \X )} \Big] \text{.} \non 
\end{align}

Table \ref{table:hb} reports values of the risks of $\fh ^{\rm{J}} ( \cdot ; \X )$ and $\fh ^{\rm{HB}} ( \cdot ; \X )$ with values of PRIAL given in parentheses. 
It can be seen from the values of PRIAL that $\fh ^{\rm{HB}} ( \cdot ; \X )$ has smaller values of risks than $\fh ^{\rm{J}} ( \cdot ; \X )$ in every case. 
When $\p = \p (0), \p (2)$, PRIAL is smallest in Case (II) and largest in Case (III). 
On the other hand, when $\p = \p (1)$, $\fh ^{\rm{HB}} ( \cdot ; \X )$ has the largest and smallest values of PRIAL in Cases (II) and (III), respectively.

\small
\begin{table}[!thb]
\caption{Risks of $\fh ^{\rm{J}} ( \cdot ; \X )$ (J) and $\fh ^{\rm{HB}} ( \cdot ; \X )$ (HB). 
Values of PRIAL of HB are given in parentheses. }
\begin{center}
$
{\renewcommand\arraystretch{1.1}\small
\begin{array}{c@{\hspace{2mm}}
              r@{\hspace{5mm}}
              r@{\hspace{2mm}}
              r@{\hspace{2mm}}
              r@{\hspace{2mm}}
              r@{\hspace{2mm}}
              r@{\hspace{2mm}}
              r@{\hspace{2mm}}
              r
             }
\hline
\text{Case} & \text{$\p $} & \text{J} & \text{HB} \\

\hline

\text{(I)} & \text{$\p (0)$}
&0.22 	

&0.22 	\, ( 1.13   )

\\
\text{(I)} & \text{$\p (1)$}
&0.23 	

&0.23 	  \, (1.08  )

\\
\text{(I)} & \text{$\p (2)$}
&0.27 	

&0.27 	  \, (  1.40)

\\
\hline

\text{(II)} & \text{$\p (0)$}
&0.28 	

&0.27 	 \, ( 1.00  )

\\

\text{(II)} & \text{$\p (1)$}
&0.32 

&0.31 	  \, ( 2.78 )

\\
\text{(II)} & \text{$\p (2)$}
&0.30 	

&0.30 	  \, (1.35  )

\\

\hline

\text{(III)} & \text{$\p (0)$}
&0.23 

&0.23  \, (	1.34   )

\\
\text{(III)} & \text{$\p (1)$}
&0.30 	

&0.29   \, ( 	0.52 )

\\
\text{(III)} & \text{$\p (2)$}
&0.25 	

&0.24 	 \, (  2.02 )

\\

\hline
\end{array}
}
$
\end{center}
\label{table:hb} 
\end{table}
\normalsize

\section{Discussion}
\label{sec:discussion}
In this paper, we considered the problems of estimating negative multinomial parameter vectors and the joint predictive density of multinomial tables on the basis of observations of negative multinomial variables in unbalanced settings. 
A related problem of mathematical interest is that of estimating the joint predictive density of future negative multinomial variables on the basis of the current negative multinomial observations. 
Although no dominance result has been obtained, we here derive identities which relate prediction to estimation in the negative multinomial case. 

Let $s_1 , \dots , s_n > 0$ and let $\Y _{\nu } = ( Y_{i, \nu } )_{i = 1}^{m_{\nu }}$, $\nu = 1, \dots , n$, be independent negative multinomial variables with mass functions 
\begin{align}
g_{\nu } ( \y _{\nu } | \p _{\nu } ) = {\Ga \big( s_{\nu } + \sum_{i = 1}^{m_{\nu }} y_{i, \nu } \big) \over \Ga ( s_{\nu } ) \prod_{i = 1}^{m_{\nu }} y_{i, \nu } !} {p_{0, \nu }}^{s_{\nu }} \prod_{i = 1}^{m_{\nu }} {p_{i, \nu }}^{y_{i, \nu }} \text{,} \label{eq:mass_Y} 
\end{align}
$\y _{\nu } = ( y_{i, \nu } )_{i = 1}^{m_{\nu }} \in {\mathbb{N} _0}^{m_{\nu }}$, $\nu = 1, \dots , n$, respectively. 
Consider the problem of estimating the predictive density $g( \y | \p ) = \prod_{\nu = 1}^{n} g_{\nu } ( \y _{\nu } | \p _{\nu } )$, $\y = ( \y _{\nu } )_{\nu = 1, \dots , n} \in {\mathbb{N} _0}^{m_1} \times \dots \times {\mathbb{N} _0}^{m_n}$, on the basis of $\X $ given in Section \ref{sec:eb} under the Kullback-Leibler divergence. 
As shown by Aitchison (1975), the Bayesian predictive mass $\gh ^{( \pi )} ( \cdot ; \X )$ with respect to a prior $\p \sim \pi ( \p )$ is given by 
\begin{align}
\gh ^{( \pi )} ( \y ; \X ) &= E_{\pi } [ g( \y | \p ) | \X ] \non \\
&= \Big\{ \prod_{\nu = 1}^{n} {\Ga \big( s_{\nu } + \sum_{i = 1}^{m_{\nu }} y_{i, \nu } \big) \over \Ga ( s_{\nu } ) \prod_{i = 1}^{m_{\nu }} y_{i, \nu } !} \Big\} \frac{ \int_D \pi ( \p ) \big\{ \prod_{\nu = 1}^{N} \big( {p_{0, \nu }}^{s_{\nu } + r_{\nu }} \prod_{i = 1}^{m_{\nu }} {p_{i, \nu }}^{y_{i, \nu } + X_{i, \nu }} \big) \big\} d\p }{ \int_D \pi ( \p ) \big\{ \prod_{\nu = 1}^{N} \big( {p_{0, \nu }}^{r_{\nu }} \prod_{i = 1}^{m_{\nu }} {p_{i, \nu }}^{X_{i, \nu }} \big) \big\} d\p } \text{,} \label{eq:bpm_Y} 
\end{align}
where $s_{\nu } = y_{1, \nu } = \dots = y_{m_{\nu } , \nu } = 0$ if $\nu \in \{ 1, \dots , N \} \cap [n + 1, \infty )$, and has risk given by 
\begin{align}
R( \p , \gh ^{( \pi )} ) &= E \Big[ \log {g( \Y | \p ) \over \gh ^{( \pi )} ( \Y ; \X )} \Big] \text{.} \label{eq:risk_Y} 
\end{align}

Let $t_1 , \dots , t_N \colon [0, 1] \to (0, \infty )$ be smooth, nondecreasing functions such that for all $\nu = 1, \dots , N$, 
\begin{align}
t_{\nu } (0) = r_{\nu } \quad \text{and} \quad t_{\nu } (1) = \begin{cases} r_{\nu } + s_{\nu } \text{,} & \text{if $\nu \le n$} \text{,} \\ r_{\nu } \text{,} & \text{if $\nu \ge n + 1$} \text{.} \end{cases} 
\end{align}
For each $\ta \in [0, 1]$, let $\Z _{\nu } ( \ta ) = ( Z_{i, \nu } ( \ta ))_{i = 1}^{m_{\nu }}$, $\nu = 1, \dots , N$, be independent negative multinomial variables with mass functions 
\begin{align}
{\Ga \big( t_{\nu } ( \ta ) + \sum_{i = 1}^{m_{\nu }} z_{i, \nu } \big) \over \Ga ( t_{\nu } ( \ta )) \prod_{i = 1}^{m_{\nu }} z_{i, \nu } !} {p_{0, \nu }}^{t_{\nu } ( \ta )} \prod_{i = 1}^{m_{\nu }} {p_{i, \nu }}^{z_{i, \nu }} \text{,} \non 
\end{align}
$( z_{i, \nu } )_{i = 1}^{m_{\nu }} \in {\mathbb{N} _0}^{m_{\nu }}$, $\nu = 1, \dots , N$, respectively, and let $\Z ( \ta ) = ( \Z _{\nu } ( \ta ))_{\nu = 1, \dots , N}$. 
Let $\mathcal{W} _{\nu , k} = \big\{ ( \mathring{w} _i )_{i = 1}^{m_{\nu }} \in {\mathbb{N} _0}^{m_{\nu }} \big| \sum_{i = 1}^{m_{\nu }} \mathring{w} _i = k \big\} $ for $\nu = 1, \dots , N$ and $k \in \mathbb{N} _0$. 
Let 
\begin{align}
L^{\rm{KL}} ( \dt , \th ) &= \dt - \th - \th \log ( \dt / \th ) \label{eq:loss_KL} 
\end{align}
for $\dt , \th \in (0, \infty )$. 
The following theorem shows that the risk function of an arbitrary Bayesian predictive mass can be expressed using the risk functions of the corresponding Bayes estimators of an infinite number of monomials of the unknown probabilities. 

\begin{thm}
\label{thm:NM} 
Let $\p \sim \pi ( \p )$ be a prior density. 
Then the risk of $\gh ^{( \pi )} ( \cdot ; \X )$ is expressed as 
\begin{align}
&R( \p , \gh ^{( \pi )} ) \non \\
&= \int_{0}^{1} \Big\{ \sum_{\nu = 1}^{n} {t_{\nu }}' ( \ta ) \sum_{k = 1}^{\infty } {1 \over k} \sum_{( w_i )_{i = 1}^{m_{\nu }} \in \mathcal{W} _{\nu , k}} {k ! \over \prod_{i = 1}^{m_{\nu }} w_i !} E \Big[ L^{\rm{KL}} \Big( E_{\pi } \Big[ \prod_{i = 1}^{m_{\nu }} {p_{i, \nu }}^{w_i} \Big| \Z ( \ta ) \Big] , \prod_{i = 1}^{m_{\nu }} {p_{i, \nu }}^{w_i} \Big) \Big] \Big\} d\ta \text{.} \non 
\end{align}
\end{thm}

Theorem 3 of Hamura and Kubokawa (2020b) is related to the monomials of degree $1$ in the above expression. 
In the negative binomial case, the ``intrinsic loss'' derived by Robert (1996) is not given by (\ref{eq:loss_KL}); see Remark 2.2 of Hamura and Kubokawa (2019a) for details. 

We also have the following somewhat simpler result. 
Let 
\begin{align}
\pi _{M, \bgat , \a _0 , \a } ( \p ) &= \int_{0}^{\infty } \Big[ \prod_{\nu = 1}^{N} \Big\{ {p_{0, \nu }}^{\gat _{\nu } (u) + a_{0, \nu } - 1} \prod_{i = 1}^{m_{\nu }} {p_{i, \nu }}^{a_{i, \nu } - 1} \Big\} \Big] dM(u) \text{,} \label{eq:general_shrinkage_prior} 
\end{align}
where $M$ is a measure  on $(0, \infty )$ while $\bgat = ( \gat _{\nu } )_{\nu = 1}^{N} \colon (0, \infty ) \to (0, \infty )^N$. 
Then Corollary \ref{cor:NM} gives an expression for the risk difference between the Bayesian predictive mass with respect to the prior (\ref{eq:general_shrinkage_prior}) and that with respect to the prior (\ref{eq:prior_Dir}). 

\begin{cor}
\label{cor:NM} 
The risk difference between $\gh ^{( \pi _{M, \bgat , \a _0 , \a } )} ( \cdot ; \X )$ and $\gh ^{( \pi _{\a _0 , \a } )} ( \cdot ; \X )$ is expressed as 
\begin{align}
&R( \p , \gh ^{( \pi _{M, \bgat , \a _0 , \a } )} ) - R( \p , \gh ^{( \pi _{\a _0 , \a } )} ) \non \\
&= \int_{0}^{1} \Big\{ \sum_{\nu = 1}^{n} {t_{\nu }}' ( \ta ) \sum_{k = 1}^{\infty } {1 \over k} E[ L^{\rm{KL}} ( E_{\pi _{M, \bgat , \a _0 , \a }} [ {p_{\cdot , \nu }}^k | \Z ( \ta ) ] , {p_{\cdot , \nu }}^k ) - L^{\rm{KL}} ( E_{\pi _{\a _0 , \a }} [ {p_{\cdot , \nu }}^k | \Z ( \ta ) ] , {p_{\cdot , \nu }}^k ) ] \Big\} d\ta \text{.} \non 
\end{align}
\end{cor}

Despite these identities, dominance conditions have not been obtained. 
It may be worth noting that $\log \{ \gh ^{( \pi _{\a _0 , \a } )} ( \Y ; \X ) / \gh ^{( \pi _{M, \bgat , \a _0 , \a } )} ( \Y ; \X ) \} $, whose expectation is the risk difference, is a function only of $X_{\cdot , \nu }$, $\nu = 1, \dots , N$, and $Y_{\cdot , \nu } = \sum_{i = 1}^{m_{\nu }} Y_{i, \nu }$, $\nu = 1, \dots , n$. 
Inadmissibility of $\gh ^{( \pi _{\a _0 , \a } )} ( \cdot ; \X )$ could be studied in a future paper.

\section{Appendix}
\label{sec:appendix} 
\subsection{Assumptions}
\label{subsec:assumptions} 
Let $\overline{c} _{\nu } = \max_{1 \le i \le m_{\nu }} c_{i, \nu }$ for $\nu = 1, \dots , N$. 
Let $\underline{\ct } _{\nu }^{( \nu )} = \min_{1 \le i \le m_{\nu }} \ct _{i, \nu }^{( \nu )}$, $\overline{\ct } _{\nu }^{( \nu )} = \max_{1 \le i \le m_{\nu }} \ct _{i, \nu }^{( \nu )}$, and $\Ct _{\nu } = ( \overline{\ct } _{\nu }^{( \nu )} / \underline{\ct } _{\nu }^{( \nu )} ) / \{ 1 + \bt _{\nu } ( \underline{\ct } _{\nu }^{( \nu )} + \overline{\ct } _{\nu }^{( \nu )} ) \} $ for $\nu = 1, \dots , N$. 
Let $A = \max_{1 \le \nu \le n} \overline{c} _{\nu } ( \Ct _{\nu } + 2)$, $\underline{\bt } = \min_{1 \le \nu \le n} \bt _{\nu }$, $\overline{\bt } = \max_{1 \le \nu \le n} \bt _{\nu }$, $\underline{\overline{\ct }} = \min_{1 \le \nu \le n} \overline{\ct } _{\nu }^{( \nu )}$, and $\overline{\overline{\ct }} = \max_{1 \le \nu \le n} \overline{\ct } _{\nu }^{( \nu )}$ and let $\ct _{*} = \min_{1 \le \nu \le N} \min_{1 \le {\nu }' \le N} \min_{1 \le i \le m_{{\nu }'}} \ct _{i, {\nu }'}^{( \nu )}$ and $\ct ^{*} = \max_{1 \le \nu \le N} \max_{1 \le {\nu }' \le N} \max_{1 \le i \le m_{{\nu }'}} \ct _{i, {\nu }'}^{( \nu )}$. 
Let $A_1 = \max_{1 \le \nu \le n} \{ \overline{c} _{\nu } (3 + 4 \bt _{\nu } \ct ) / (1 + 2 \bt _{\nu } \ct ) \} $. 

Assumption \ref{assumption:tEB_affine} and Assumption \ref{assumption:cEB_affine} correspond to Theorem \ref{thm:EB_affine} and Corollary \ref{cor:EB_affine}, respectively. 

\begin{assumption}
\label{assumption:tEB_affine} 
\hfill
\begin{enumerate}
\item[{\rm{(a)}}]
$\overline{\overline{c}} > 0$. 
\item[{\rm{(b)}}]
$r_{\nu } \ge \Ct _{\nu } + 1$ and $r_{\nu } + \bt _{\nu } \ge \Ct _{\nu } + 2$ for all $\nu = 1, \dots , n$ with $c_{\cdot , \nu } > 0$. 
\item[{\rm{(c)}}]
$\underline{c_{\cdot }} - A \ge 0$. 
\item[{\rm{(d)}}]
For all $x \in \mathbb{N}$, either 
\begin{itemize}
\item
$\overline{\overline{c}} \{ \overline{\bt } + 1 / ( \ct _{*} x + \underline{\overline{\ct }} ) \} - 2 ( {\underline{r} / \overline{r}} )^2 ( \underline{c_{\cdot }} - A) \{ \underline{\bt } \ct _{*} \underline{\overline{\ct }} / ( \overline{\bt } \ct ^{*} \overline{\overline{\ct }} ) \} \le 0$ implies 
\begin{align}
\underline{c_{\cdot }} \{ \overline{\bt } + 1 / ( \ct _{*} x + \underline{\overline{\ct }} ) \} - 2 \Big( {\underline{r} \over \overline{r}} \Big) ^2 ( \underline{c_{\cdot }} - A) {\underline{\bt } \ct _{*} \underline{\overline{\ct }} \over \overline{\bt } \ct ^{*} \overline{\overline{\ct }}} \{ \underline{r} + \overline{\bt } + 1 / ( \ct _{*} x + \underline{\overline{\ct }} ) \} \le 0 \quad \text{and} \non 
\end{align}
\item
$\overline{\overline{c}} \{ \overline{\bt } + 1 / ( \ct _{*} x + \underline{\overline{\ct }} ) \} - 2 ( {\underline{r} / \overline{r}} )^2 ( \underline{c_{\cdot }} - A) \{ \underline{\bt } \ct _{*} \underline{\overline{\ct }} / ( \overline{\bt } \ct ^{*} \overline{\overline{\ct }} ) \} > 0$ implies 
\begin{align}
&x \Big[ \overline{\overline{c}} \{ \overline{\bt } + 1 / ( \ct _{*} x + \underline{\overline{\ct }} ) \} - 2 \Big( {\underline{r} \over \overline{r}} \Big) ^2 ( \underline{c_{\cdot }} - A) {\underline{\bt } \ct _{*} \underline{\overline{\ct }} \over \overline{\bt } \ct ^{*} \overline{\overline{\ct }}} \Big] \non \\
&+ n \underline{c_{\cdot }} \{ \overline{\bt } + 1 / ( \ct _{*} x + \underline{\overline{\ct }} ) \} - 2 n \Big( {\underline{r} \over \overline{r}} \Big) ^2 ( \underline{c_{\cdot }} - A) {\underline{\bt } \ct _{*} \underline{\overline{\ct }} \over \overline{\bt } \ct ^{*} \overline{\overline{\ct }}} \{ \underline{r} + \overline{\bt } + 1 / ( \ct _{*} x + \underline{\overline{\ct }} ) \} \le 0 \non 
\end{align}
\end{itemize}
or 
\begin{itemize}
\item
$\overline{\overline{c}} \{ \overline{\bt } + 1 / ( \ct _{*} x + \underline{\overline{\ct }} ) \} - 2 ( \underline{c_{\cdot }} - A) \{ \underline{\bt } \ct _{*} \underline{\overline{\ct }} / ( \overline{\bt } \ct ^{*} \overline{\overline{\ct }} ) \} \le 0$ implies 
\begin{align}
( \underline{c_{\cdot }} - \overline{\overline{c}} \underline{r} ) - 2 ( \underline{c_{\cdot }} - A) {\underline{\bt } \ct _{*} \underline{\overline{\ct }} \over \overline{\bt } \ct ^{*} \overline{\overline{\ct }}} \le 0 \quad \text{and} \non 
\end{align}
\item
$\overline{\overline{c}} \{ \overline{\bt } + 1 / ( \ct _{*} x + \underline{\overline{\ct }} ) \} - 2 ( \underline{c_{\cdot }} - A) \{ \underline{\bt } \ct _{*} \underline{\overline{\ct }} / ( \overline{\bt } \ct ^{*} \overline{\overline{\ct }} ) \} > 0$ implies 
\begin{align}
&\Big( \sum_{\nu = 1}^{n} r_{\nu } + x \Big) \Big[ \overline{\overline{c}} \{ \overline{\bt } + 1 / ( \ct _{*} x + \underline{\overline{\ct }} ) \} - 2 ( \underline{c_{\cdot }} - A) {\underline{\bt } \ct _{*} \underline{\overline{\ct }} \over \overline{\bt } \ct ^{*} \overline{\overline{\ct }}} \Big] \non \\
&+ n ( \underline{c_{\cdot }} - \overline{\overline{c}} \underline{r} ) \{ \overline{\bt } + 1 / ( \ct _{*} x + \underline{\overline{\ct }} ) \} - 2 n ( \underline{c_{\cdot }} - A) {\underline{\bt } \ct _{*} \underline{\overline{\ct }} \over \overline{\bt } \ct ^{*} \overline{\overline{\ct }}} \{ \overline{\bt } + 1 / ( \ct _{*} x + \underline{\overline{\ct }} ) \} \le 0 \text{.} \non 
\end{align}
\end{itemize}
\end{enumerate}
\end{assumption}

\begin{assumption}
\label{assumption:cEB_affine} 
\hfill
\begin{enumerate}
\item[{\rm{(a)}}]
$\overline{\overline{c}} > 0$. 
\item[{\rm{(b)}}]
$r_{\nu } \ge 1 / (1 + 2 \bt _{\nu } \ct ) + 1$ and $r_{\nu } + \bt _{\nu } \ge 1 / (1 + 2 \bt _{\nu } \ct ) + 2$ for all $\nu = 1, \dots , n$ with $c_{\cdot , \nu } > 0$. 
\item[{\rm{(c)}}]
$\underline{c_{\cdot }} - A_1 \ge 0$. 
\item[{\rm{(d)}}]
For all $x \in \mathbb{N}$, either 
\begin{itemize}
\item
$\overline{\overline{c}} [ \overline{\bt } + 1 / \{ \ct (x + 1) \} ] - 2 ( {\underline{r} / \overline{r}} )^2 ( \underline{c_{\cdot }} - A_1 ) \underline{\bt } / \overline{\bt } \le 0$ implies 
\begin{align}
\underline{c_{\cdot }} [ \overline{\bt } + 1 / \{ \ct (x + 1) \} ] - 2 \Big( {\underline{r} \over \overline{r}} \Big) ^2 ( \underline{c_{\cdot }} - A_1 ) {\underline{\bt } \over \overline{\bt }} [ \underline{r} + \overline{\bt } + 1 / \{ \ct (x + 1) \} ] \le 0 \quad \text{and} \non 
\end{align}
\item
$\overline{\overline{c}} [ \overline{\bt } + 1 / \{ \ct (x + 1) \} ] - 2 ( {\underline{r} / \overline{r}} )^2 ( \underline{c_{\cdot }} - A_1 ) \underline{\bt } / \overline{\bt } > 0$ implies 
\begin{align}
&x \Big( \overline{\overline{c}} [ \overline{\bt } + 1 / \{ \ct (x + 1) \} ] - 2 \Big( {\underline{r} \over \overline{r}} \Big) ^2 ( \underline{c_{\cdot }} - A_1 ) {\underline{\bt } \over \overline{\bt }} \Big) \non \\
&+ n \underline{c_{\cdot }} [ \overline{\bt } + 1 / \{ \ct (x + 1) \} ] - 2 n \Big( {\underline{r} \over \overline{r}} \Big) ^2 ( \underline{c_{\cdot }} - A_1 ) {\underline{\bt } \over \overline{\bt }} [ \underline{r} + \overline{\bt } + 1 / \{ \ct (x + 1) \} ] \le 0 \non 
\end{align}
\end{itemize}
or 
\begin{itemize}
\item
$\overline{\overline{c}} [ \overline{\bt } + 1 / \{ \ct (x + 1) \} ] - 2 ( \underline{c_{\cdot }} - A_1 ) \underline{\bt } / \overline{\bt } \le 0$ implies 
\begin{align}
( \underline{c_{\cdot }} - \overline{\overline{c}} \underline{r} ) - 2 ( \underline{c_{\cdot }} - A_1 ) {\underline{\bt } \over \overline{\bt }} \le 0 \quad \text{and} \non 
\end{align}
\item
$\overline{\overline{c}} [ \overline{\bt } + 1 / \{ \ct (x + 1) \} ] - 2 ( \underline{c_{\cdot }} - A_1 ) \underline{\bt } / \overline{\bt } > 0$ implies 
\begin{align}
&\Big( \sum_{\nu = 1}^{n} r_{\nu } + x \Big) \Big( \overline{\overline{c}} [ \overline{\bt } + 1 / \{ \ct (x + 1) \} ] - 2 ( \underline{c_{\cdot }} - A_1 ) {\underline{\bt } \over \overline{\bt }} \Big) \non \\
&+ n ( \underline{c_{\cdot }} - \overline{\overline{c}} \underline{r} ) [ \overline{\bt } + 1 / \{ \ct (x + 1) \} ] - 2 n ( \underline{c_{\cdot }} - A_1 ) {\underline{\bt } \over \overline{\bt }} [ \overline{\bt } + 1 / \{ \ct (x + 1) \} ] \le 0 \text{.} \non 
\end{align}
\end{itemize}
\end{enumerate}
\end{assumption}

\subsection{Proofs}
\label{subsec:proofs} 
Here we prove Theorems \ref{thm:EB}, \ref{thm:EB_affine}, \ref{thm:multin}, and \ref{thm:NM}, Lemma \ref{lem:likelihood}, and Corollary \ref{cor:NM}. 
We use Lemma \ref{lem:hudson}, which is due to Hudson (1978). 

For $(i, \nu ), ( i' , {\nu }' ) \in \mathbb{N} \times \{ 1, \dots , N \} $ with $i \le m_{\nu }$ and $i' \le m_{{\nu }'}$, let $\de _{i, i' , \nu , {\nu }'}$ $= 1$ if $i = i'$ and $\nu = {\nu }'$ and $= 0$ otherwise. 
Let $\X _{\cdot } = ( X_{\cdot , \nu } )_{\nu = 1}^{N}$. 
For $\nu = 1, \dots , N$, let $\e _{\nu }^{(N)}$ be the $\nu $th unit vector in $\mathbb{R} ^N$, namely the $\nu $th column of the $N \times N$ identity matrix. 
For $\nu = 1, \dots , N$, let $\bm{0} ^{( m_{\nu } )} = (0, \dots , 0)^{\top } \in \mathbb{R} ^{m_{\nu }}$. 
For $\nu , {\nu }' = 1, \dots , N$, let $\de _{\nu , {\nu }'}^{(N)} = {\e _{\nu }^{(N)}}^{\top } \e _{{\nu }'}^{(N)}$. 

\begin{lem}
\label{lem:hudson} 
Let $\varphi \colon {\mathbb{N} _0}^{m_1} \times \dots \times {\mathbb{N} _0}^{m_N} \to \mathbb{R}$ and suppose that either $\varphi ( \x ) \ge 0$ for all $\x \in {\mathbb{N} _0}^{m_1} \times \dots \times {\mathbb{N} _0}^{m_N}$ or $E[ | \varphi( \X )| ] < \infty $. 
Then for all $(i, \nu ) \in \mathbb{N} \times \{ 1, \dots , N \} $ with $i \le m_{\nu }$, if $\varphi ( \x ) = 0$ for all $\x = (( x_{i' , {\nu }'} )_{i' = 1}^{m_{{\nu }'}} )_{{\nu }' = 1, \dots , N} \in {\mathbb{N} _0}^{m_1} \times \dots \times {\mathbb{N} _0}^{m_N}$ such that $x_{i, \nu } = 0$, we have 
\begin{align}
E \Big[ {\varphi ( \X ) \over p_{i, \nu }} \Big] = E \Big[ {r_{\nu } + X_{\cdot , \nu } \over X_{i, \nu } + 1} \varphi ( \X + \e _{i, \nu } ) \Big] \text{,} \non 
\end{align}
where $\X + \e _{i, \nu } = (( X _{i' , {\nu }'} + \de _{i, i' , \nu , {\nu }'} )_{i' = 1}^{m_{{\nu }'}} )_{{\nu }' = 1, \dots , N}$. 
\end{lem}

\noindent
{\bf Proof of Theorem \ref{thm:EB}.} \ \ Let $\De _{\c }^{( \bde )} = E[ L_{\c } ( \pbh ^{( \bde )} , \p ) ] - E[ L_{\c } ( \pbh ^{\rm{U}} , \p ) ]$. 
For $\nu = 1, \dots , N$, let 
\begin{align}
\phi _{\nu }^{( \bde )} ( \X _{\cdot } ) &= \begin{cases} \displaystyle {\de _{\nu } ( X_{\cdot , \cdot } ) \over r_{\nu } + X_{\cdot , \nu } - 1 + \de _{\nu } ( X_{\cdot , \cdot } )} \text{,} & \text{if $X_{\cdot , \nu } \ge 1$} \text{,} \\ \displaystyle 0 \text{,} & \text{if $X_{\cdot , \nu } = 0$} \text{,} \end{cases} \non 
\end{align}
so that $\ph _{i, \nu }^{( \bde )} = \ph _{i, \nu }^{\rm{U}} - \ph _{i, \nu }^{\rm{U}} \phi _{\nu }^{( \bde )}  ( \X _{\cdot } )$ for all $i = 1, \dots , m_{\nu }$. 
Then, by Lemma \ref{lem:hudson}, 
\begin{align}
\De _{\c }^{( \bde )} &= E \Big[ \sum_{\nu = 1}^{n} \sum_{i = 1}^{m_{\nu }} \Big[ c_{i, \nu } { ( \ph _{i, \nu }^{\rm{U}} )^2 \{ \phi _{\nu }^{( \bde )} ( \X _{\cdot } ) \} ^2 - 2 ( \ph _{i, \nu }^{\rm{U}} )^2 \phi _{\nu }^{( \bde )} ( \X _{\cdot } ) \over p_{i, \nu }} + 2 c_{i, \nu } \ph _{i, \nu }^{\rm{U}} \phi _{\nu }^{( \bde )} ( \X _{\cdot } ) \Big] \Big] \non \\
&= E \Big[ \sum_{\nu = 1}^{n} \sum_{i = 1}^{m_{\nu }} \Big( c_{i, \nu } {X_{i, \nu } + 1 \over r_{\nu } + X_{\cdot , \nu }} [ \{ \phi _{\nu }^{( \bde )} ( \X _{\cdot } + \e _{\nu }^{(N)} ) \} ^2 - 2 \phi _{\nu }^{( \bde )} ( \X _{\cdot } + \e _{\nu }^{(N)} )] \non \\
&\quad + 2 c_{i, \nu } {X_{i, \nu } \over r_{\nu } + X_{\cdot , \nu } - 1} \phi _{\nu }^{( \bde )} ( \X _{\cdot } ) \Big) \Big] \non \\
&= E \Big[ \sum_{\nu = 1}^{n} \{ I_{1, \nu }^{( \bde )} ( \X ) - 2 I_{2, \nu }^{( \bde )} ( \X ) + 2 I_{3, \nu }^{( \bde )} ( \X ) \} \Big] \text{,} \non 
\end{align}
where 
\begin{align}
I_{1, \nu }^{( \bde )} ( \x ) &= {\sum_{i = 1}^{m_{\nu }} c_{i, \nu } x_{i, \nu } + c_{\cdot , \nu } \over r_{\nu } + \sum_{i = 1}^{m_{\nu }} x_{i, \nu }} \Big\{ {\de _{\nu } \big( \sum_{\nu = 1}^{N} \sum_{i = 1}^{m_{\nu }} x_{i, \nu } + 1 \big) \over r_{\nu } + \sum_{i = 1}^{m_{\nu }} x_{i, \nu } + \de _{\nu } \big( \sum_{\nu = 1}^{N} \sum_{i = 1}^{m_{\nu }} x_{i, \nu } + 1 \big) } \Big\} ^2 \text{,} \non \\
I_{2, \nu }^{( \bde )} ( \x ) &= {\sum_{i = 1}^{m_{\nu }} c_{i, \nu } x_{i, \nu } + c_{\cdot , \nu } \over r_{\nu } + \sum_{i = 1}^{m_{\nu }} x_{i, \nu }} {\de _{\nu } \big( \sum_{\nu = 1}^{N} \sum_{i = 1}^{m_{\nu }} x_{i, \nu } + 1 \big) \over r_{\nu } + \sum_{i = 1}^{m_{\nu }} x_{i, \nu } + \de _{\nu } \big( \sum_{\nu = 1}^{N} \sum_{i = 1}^{m_{\nu }} x_{i, \nu } + 1 \big) } \text{,} \non \\
I_{3, \nu }^{( \bde )} ( \x ) &= {\big( \sum_{i = 1}^{m_{\nu }} c_{i, \nu } x_{i, \nu } \big) \de _{\nu } \big( \sum_{\nu = 1}^{N} \sum_{i = 1}^{m_{\nu }} x_{i, \nu } \big) \over \big( r_{\nu } + \sum_{i = 1}^{m_{\nu }} x_{i, \nu } - 1 \big) \big\{ r_{\nu } + \sum_{i = 1}^{m_{\nu }} x_{i, \nu } - 1 + \de _{\nu } \big( \sum_{\nu = 1}^{N} \sum_{i = 1}^{m_{\nu }} x_{i, \nu } \big) \big\} } \text{,} \non 
\end{align}
for $\x = (( x_{i, {\nu }'} )_{i = 1}^{m_{{\nu }'}} )_{{\nu }' = 1, \dots , N} \in {\mathbb{N} _0}^{m_1} \times \dots \times {\mathbb{N} _0}^{m_N}$ for each $\nu = 1, \dots , N$. 
Since $\overline{\overline{c}} > 0$, it follows that $\sum_{\nu = 1}^{n} \{ I_{1, \nu }^{( \bde )} (( \bm{0} ^{( m_{\nu } )} )_{\nu = 1, \dots , N} ) - 2 I_{2, \nu }^{( \bde )} (( \bm{0} ^{( m_{\nu } )} )_{\nu = 1, \dots , N} ) + 2 I_{3, \nu }^{( \bde )} (( \bm{0} ^{( m_{\nu } )} )_{\nu = 1, \dots , N} ) \} < 0$. 

Fix $\x = (( x_{i, \nu } )_{i = 1}^{m_{\nu }} )_{\nu = 1, \dots , N} \in ( {\mathbb{N} _0}^{m_1} \times \dots \times {\mathbb{N} _0}^{m_N} ) \setminus \{ ( \bm{0} ^{( m_{\nu } )} )_{\nu = 1, \dots , N} \} $. 
It is sufficient to show that $\sum_{\nu = 1}^{n} \{ I_{1, \nu }^{( \bde )} ( \x ) - 2 I_{2, \nu }^{( \bde )} ( \x ) + 2 I_{3, \nu }^{( \bde )} ( \x ) \} \le 0$. 
Let $x_{\cdot , \nu } = \sum_{i = 1}^{m_{\nu }} x_{i, \nu }$ for $\nu = 1, \dots , N$ and let $x_{\cdot , \cdot } = \sum_{\nu = 1}^{N} x_{\cdot , \nu }$. 
Let $\overline{c} _{\nu } = \max_{1 \le i \le m_{\nu }} c_{i, \nu }$ for $\nu = 1, \dots , N$. 
Then for all $\nu = 1, \dots , n$ such that $\sum_{i = 1}^{m_{\nu }} c_{i, \nu } x_{i, \nu } > 0$, since, by (\ref{eq:assumption_EB_0}), $\de _{\nu } ( x_{\cdot , \cdot } ) \le \{ ( x_{\cdot , \cdot } + 1) / x_{\cdot , \cdot } \} \de _{\nu } ( x_{\cdot , \cdot } + 1) \le \{ ( x_{\cdot , \nu } + 1) / x_{\cdot , \nu } \} \de _{\nu } ( x_{\cdot , \cdot } + 1)$, we have that 
\begin{align}
I_{3, \nu }^{( \bde )} ( \x ) &\le {\sum_{i = 1}^{m_{\nu }} c_{i, \nu } x_{i, \nu } \over r_{\nu } + x_{\cdot , \nu } - 1} {\de_{\nu } ( x_{\cdot , \cdot } + 1) \over \{ x_{\cdot , \nu } / ( x_{\cdot , \nu } + 1) \} ( r_{\nu } + x_{\cdot , \nu } - 1) + \de_{\nu } ( x_{\cdot , \cdot } + 1)} \non 
\end{align}
and hence that 
\begin{align}
- I_{2, \nu }^{( \bde )} ( \x ) + I_{3, \nu }^{( \bde )} ( \x ) &\le - {c_{\cdot , \nu } \over r_{\nu } + x_{\cdot , \nu }} {\de _{\nu } ( x_{\cdot , \cdot } + 1) \over r_{\nu } + x_{\cdot , \nu } + \de _{\nu } ( x_{\cdot , \cdot } + 1) } \non \\
&\quad + \Big( \sum_{i = 1}^{m_{\nu }} c_{i, \nu } x_{i, \nu } \Big) \de_{\nu } ( x_{\cdot , \cdot } + 1) \Big[ - {1 \over r_{\nu } + x_{\cdot , \nu }} {1 \over r_{\nu } + x_{\cdot , \nu } + \de_{\nu } ( x_{\cdot , \cdot } + 1)} \non \\
&\quad + {1 \over r_{\nu } + x_{\cdot , \nu } - 1} {1 \over \{ x_{\cdot , \nu } / ( x_{\cdot , \nu } + 1) \} ( r_{\nu } + x_{\cdot , \nu } - 1) + \de _{\nu } ( x_{\cdot , \cdot } + 1)} \Big] \non \\
&\le - {c_{\cdot , \nu } \over r_{\nu } + x_{\cdot , \nu }} {\de _{\nu } ( x_{\cdot , \cdot } + 1) \over r_{\nu } + x_{\cdot , \nu } + \de _{\nu } ( x_{\cdot , \cdot } + 1) } \non \\
&\quad + \overline{c} _{\nu } x_{\cdot , \nu } \de_{\nu } ( x_{\cdot , \cdot } + 1) \Big[ - {1 \over r_{\nu } + x_{\cdot , \nu }} {1 \over r_{\nu } + x_{\cdot , \nu } + \de_{\nu } ( x_{\cdot , \cdot } + 1)} \non \\
&\quad + {1 \over r_{\nu } + x_{\cdot , \nu } - 1} {1 \over \{ x_{\cdot , \nu } / ( x_{\cdot , \nu } + 1) \} ( r_{\nu } + x_{\cdot , \nu } - 1) + \de _{\nu } ( x_{\cdot , \cdot } + 1)} \Big] \text{,} \non 
\end{align}
where 
\begin{align}
{1 \over r_{\nu } + x_{\cdot , \nu } - 1} {1 \over \{ x_{\cdot , \nu } / ( x_{\cdot , \nu } + 1) \} ( r_{\nu } + x_{\cdot , \nu } - 1) + \de _{\nu } ( x_{\cdot , \cdot } + 1)} \le {x_{\cdot , \nu } + 3 \over r_{\nu } + x_{\cdot , \nu }} {1 / x_{\cdot , \nu } \over r_{\nu } + x_{\cdot , \nu } + \de _{\nu } ( x_{\cdot , \cdot } + 1)} \non 
\end{align}
by the assumption that $r_{\nu } \ge 5 / 2$ for all $\nu = 1, \dots , n$ with $c_{\cdot , \nu } > 0$. 
Thus, for any $\nu = 1, \dots , n$, 
\begin{align}
&I_{1, \nu }^{( \bde )} ( \x ) - 2 I_{2, \nu }^{( \bde )} ( \x ) + 2 I_{3, \nu }^{( \bde )} ( \x ) \non \\
&\le {\overline{c} _{\nu } x_{\cdot , \nu } + c_{\cdot , \nu } \over r_{\nu } + x_{\cdot , \nu }} \Big\{ {\de _{\nu } ( x_{\cdot , \cdot } + 1) \over r_{\nu } + x_{\cdot , \nu } + \de _{\nu } ( x_{\cdot , \cdot } + 1)} \Big\} ^2 + 2 {3 \overline{c} _{\nu } - c_{\cdot , \nu } \over r_{\nu } + x_{\cdot , \nu }} {\de _{\nu } ( x_{\cdot , \cdot } + 1) \over r_{\nu } + x_{\cdot , \nu } + \de _{\nu } ( x_{\cdot , \cdot } + 1)} \non \\
&= {\de _{\nu } ( x_{\cdot , \cdot } + 1) [( \overline{c} _{\nu } x_{\cdot , \nu } + c_{\cdot , \nu } ) \de _{\nu } ( x_{\cdot , \cdot } + 1) - 2 ( c_{\cdot , \nu } - 3 \overline{c} _{\nu } ) \{ r_{\nu } + x_{\cdot , \nu } + \de _{\nu } ( x_{\cdot , \cdot } + 1) \} ] \over ( r_{\nu } + x_{\cdot , \nu } ) \{ r_{\nu } + x_{\cdot , \nu } + \de _{\nu } ( x_{\cdot , \cdot } + 1) \} ^2} \non \\
&\le {\de _{\nu } ( x_{\cdot , \cdot } + 1) [( \overline{\overline{c}} x_{\cdot , \nu } + \underline{c_{\cdot }} ) \de _{\nu } ( x_{\cdot , \cdot } + 1) - 2 ( \underline{c_{\cdot }} - 3 \overline{\overline{c}} ) \{ r_{\nu } + x_{\cdot , \nu } + \de _{\nu } ( x_{\cdot , \cdot } + 1) \} ] \over ( r_{\nu } + x_{\cdot , \nu } ) \{ r_{\nu } + x_{\cdot , \nu } + \de _{\nu } ( x_{\cdot , \cdot } + 1) \} ^2} \non \\%
&\le {\overline{\overline{c}} x_{\cdot , \nu } + \underline{c_{\cdot }} \over r_{\nu } + x_{\cdot , \nu }} \Big\{ {\overline{\de } ( x_{\cdot , \cdot } + 1) \over r_{\nu } + x_{\cdot , \nu } + \overline{\de } ( x_{\cdot , \cdot } + 1)} \Big\} ^2 - 2 {\underline{c_{\cdot }} - 3 \overline{\overline{c}} \over r_{\nu } + x_{\cdot , \nu }} {\underline{\de } ( x_{\cdot , \cdot } + 1) \over r_{\nu } + x_{\cdot , \nu } + \underline{\de } ( x_{\cdot , \cdot } + 1)} \label{eq:tEBp1} %
\end{align}
by the assumption that $3 \overline{\overline{c}} \le \underline{c_{\cdot }}$. 

For part (i), we have by (\ref{eq:tEBp1}) that for any $\nu = 1, \dots , n$, 
\begin{align}
&I_{1, \nu }^{( \bde )} ( \x ) - 2 I_{2, \nu }^{( \bde )} ( \x ) + 2 I_{3, \nu }^{( \bde )} ( \x ) \non \\
&\le {\overline{\overline{c}} x_{\cdot , \nu } + \underline{c_{\cdot }} \over \underline{r} + x_{\cdot , \nu }} \Big\{ {\overline{\de } ( x_{\cdot , \cdot } + 1) \over \underline{r} + x_{\cdot , \nu } + \overline{\de } ( x_{\cdot , \cdot } + 1)} \Big\} ^2 - 2 {\underline{c_{\cdot }} - 3 \overline{\overline{c}} \over \overline{r} + x_{\cdot , \nu }} {\underline{\de } ( x_{\cdot , \cdot } + 1) \over \overline{r} + x_{\cdot , \nu } + \underline{\de } ( x_{\cdot , \cdot } + 1)} \non \\
&\le {1 \over \underline{r} + x_{\cdot , \nu }} {\overline{\de } ( x_{\cdot , \cdot } + 1) \over \{ \underline{r} + x_{\cdot , \nu } + \overline{\de } ( x_{\cdot , \cdot } + 1) \} ^2} \non \\
&\quad \times [ x_{\cdot , \nu } \{ \overline{\overline{c}} \overline{\de } ( x_{\cdot , \cdot } + 1) - 2 ( \underline{r} / \overline{r} )^2 ( \underline{c_{\cdot }} - 3 \overline{\overline{c}} ) \rho \} + \underline{c_{\cdot }} \overline{\de } ( x_{\cdot , \cdot } + 1) - 2 ( \underline{r} / \overline{r} )^2 ( \underline{c_{\cdot }} - 3 \overline{\overline{c}} ) \rho \{ \underline{r} + \overline{\de } ( x_{\cdot , \cdot } + 1) \} ] \text{,} \non 
\end{align}
which is nonpositive by (\ref{eq:assumption_EB_1_1}) if $\overline{\overline{c}} \overline{\de } ( x_{\cdot , \cdot } + 1) - 2 ( \underline{r} / \overline{r} )^2 ( \underline{c_{\cdot }} - 3 \overline{\overline{c}} ) \rho \le 0$. 
On the other hand, if $\overline{\overline{c}} \overline{\de } ( x_{\cdot , \cdot } + 1) - 2 ( \underline{r} / \overline{r} )^2 ( \underline{c_{\cdot }} - 3 \overline{\overline{c}} ) \rho > 0$, then, by the covariance inequality, 
\begin{align}
&\sum_{\nu = 1}^{n} \{ I_{1, \nu }^{( \bde )} ( \x ) - 2 I_{2, \nu }^{( \bde )} ( \x ) + 2 I_{3, \nu }^{( \bde )} ( \x ) \} \non \\
&\le {1 \over n} \Big[ \sum_{\nu = 1}^{n} {1 \over \underline{r} + x_{\cdot , \nu }} {\overline{\de } ( x_{\cdot , \cdot } + 1) \over \{ \underline{r} + x_{\cdot , \nu } + \overline{\de } ( x_{\cdot , \cdot } + 1) \} ^2} \Big] \non \\
&\quad \times [ x_{\cdot , \cdot } \{ \overline{\overline{c}} \overline{\de } ( x_{\cdot , \cdot } + 1) - 2 ( \underline{r} / \overline{r} )^2 ( \underline{c_{\cdot }} - 3 \overline{\overline{c}} ) \rho \} + n \underline{c_{\cdot }} \overline{\de } ( x_{\cdot , \cdot } + 1) - 2 n ( \underline{r} / \overline{r} )^2 ( \underline{c_{\cdot }} - 3 \overline{\overline{c}} ) \rho \{ \underline{r} + \overline{\de } ( x_{\cdot , \cdot } + 1) \} ] \text{,} \non 
\end{align}
which is nonpositive by (\ref{eq:assumption_EB_1_2}). 
This proves part (i). 

For part (ii), it follows from (\ref{eq:tEBp1}) that for all $\nu = 1, \dots , n$, 
\begin{align}
&I_{1, \nu }^{( \bde )} ( \x ) - 2 I_{2, \nu }^{( \bde )} ( \x ) + 2 I_{3, \nu }^{( \bde )} ( \x ) \non \\
&\le {1 \over r_{\nu } + x_{\cdot , \nu }} {\overline{\de } ( x_{\cdot , \cdot } + 1) \over \{ r_{\nu } + x_{\cdot , \nu } + \overline{\de } ( x_{\cdot , \cdot } + 1) \} ^2} [( \overline{\overline{c}} x_{\cdot , \nu } + \underline{c_{\cdot }} ) \overline{\de } ( x_{\cdot , \cdot } + 1) - 2 ( \underline{c_{\cdot }} - 3 \overline{\overline{c}} ) \rho \{ r_{\nu } + x_{\cdot , \nu } + \overline{\de } ( x_{\cdot , \cdot } + 1) \} ] \non \\
&\le {1 \over r_{\nu } + x_{\cdot , \nu }} {\overline{\de } ( x_{\cdot , \cdot } + 1) \over \{ r_{\nu } + x_{\cdot , \nu } + \overline{\de } ( x_{\cdot , \cdot } + 1) \} ^2} \non \\
&\quad \times [( r_{\nu } + x_{\cdot , \nu } ) \{ \overline{\overline{c}} \overline{\de } ( x_{\cdot , \cdot } + 1) - 2 ( \underline{c_{\cdot }} - 3 \overline{\overline{c}} ) \rho \} + \{ \underline{c_{\cdot }} - \underline{r} \overline{\overline{c}} - 2 ( \underline{c_{\cdot }} - 3 \overline{\overline{c}} ) \rho \} \overline{\de } ( x_{\cdot , \cdot } + 1) ] \text{,} \non 
\end{align}
which is nonpositive by (\ref{eq:assumption_EB_2_1}) if $\overline{\overline{c}} \overline{\de } ( x_{\cdot , \cdot } + 1) - 2 ( \underline{c_{\cdot }} - 3 \overline{\overline{c}} ) \rho \le 0$. 
If $\overline{\overline{c}} \overline{\de } ( x_{\cdot , \cdot } + 1) - 2 ( \underline{c_{\cdot }} - 3 \overline{\overline{c}} ) \rho > 0$, then, by the covariance inequality, 
\begin{align}
&\sum_{\nu = 1}^{n} \{ I_{1, \nu }^{( \bde )} ( \x ) - 2 I_{2, \nu }^{( \bde )} ( \x ) + 2 I_{3, \nu }^{( \bde )} ( \x ) \} \non \\
&\le {1 \over n} \Big[ \sum_{\nu = 1}^{n} {1 \over r_{\nu } + x_{\cdot , \nu }} {\overline{\de } ( x_{\cdot , \cdot } + 1) \over \{ r_{\nu } + x_{\cdot , \nu } + \overline{\de } ( x_{\cdot , \cdot } + 1) \} ^2} \Big] \non \\
&\quad \times \Big[ \Big( \sum_{\nu = 1}^{n} r_{\nu } + x_{\cdot , \cdot } \Big) \{ \overline{\overline{c}} \overline{\de } ( x_{\cdot , \cdot } + 1) - 2 ( \underline{c_{\cdot }} - 3 \overline{\overline{c}} ) \rho \} + n \{ \underline{c_{\cdot }} - \underline{r} \overline{\overline{c}} - 2 ( \underline{c_{\cdot }} - 3 \overline{\overline{c}} ) \rho \} \overline{\de } ( x_{\cdot , \cdot } + 1) \Big] \text{,} \non 
\end{align}
which is nonpositive by (\ref{eq:assumption_EB_2_2}). 
This proves part (ii). 
\hfill$\Box$

\bigskip

\noindent
{\bf Proof of Theorem \ref{thm:EB_affine}.} \ \ Let $\De _{\c }^{( \bbt , \cbt )} = E[ L_{\c } ( \pbh ^{( \bbt , \cbt )} , \p ) ] - E[ L_{\c } ( \pbh ^{\rm{U}} , \p ) ]$. 
For $\nu = 1, \dots , N$, let 
\begin{align}
\tilde{\de } _{\nu }^{( \bbt , \cbt )} ( \Xt ^{( \cbt ^{( \nu )} )} ) &= \begin{cases} \bt _{\nu } + 1 / \Xt ^{( \cbt ^{( \nu )} )} \text{,} & \text{if $\Xt ^{( \cbt ^{( \nu )} )} > 0$} \text{,} \\ 0 \text{,} & \text{if $\Xt ^{( \cbt ^{( \nu )} )} = 0$} \text{,} \end{cases} \non 
\end{align}
so that 
\begin{align}
\ph _{i, \nu }^{( \bbt , \cbt )} &= \ph _{i, \nu }^{\rm{U}} - {\ph _{i, \nu }^{\rm{U}} \tilde{\de } _{\nu }^{( \bbt , \cbt )} ( \Xt ^{( \cbt ^{( \nu )} )} ) \over r_{\nu } + X_{\cdot , \nu } - 1 + \tilde{\de } _{\nu }^{( \bbt , \cbt )} ( \Xt ^{( \cbt ^{( \nu )} )} )} \non 
\end{align}
for all $i = 1, \dots , m_{\nu }$. 
By Lemma \ref{lem:hudson}, we have 
\begin{align}
\De _{\c }^{( \bbt , \cbt )} &= E \Big[ \sum_{\nu = 1}^{n} \sum_{i = 1}^{m_{\nu }} \Big( {c_{i, \nu } \over p_{i, \nu }} \Big[ {( \ph _{i, \nu }^{\rm{U}} )^2 \{ \tilde{\de } _{\nu }^{( \bbt , \cbt )} ( \Xt ^{( \cbt ^{( \nu )} )} ) \} ^2 \over \{ r_{\nu } + X_{\cdot , \nu } - 1 + \tilde{\de } _{\nu }^{( \bbt , \cbt )} ( \Xt ^{( \cbt ^{( \nu )} )} ) \} ^2} \non \\
&\quad - 2 {( \ph _{i, \nu }^{\rm{U}} )^2 \tilde{\de } _{\nu }^{( \bbt , \cbt )} ( \Xt ^{( \cbt ^{( \nu )} )} ) \over r_{\nu } + X_{\cdot , \nu } - 1 + \tilde{\de } _{\nu }^{( \bbt , \cbt )} ( \Xt ^{( \cbt ^{( \nu )} )} )} \Big] + 2 c_{i, \nu } {\ph _{i, \nu }^{\rm{U}} \tilde{\de } _{\nu }^{( \bbt , \cbt )} ( \Xt ^{( \cbt ^{( \nu )} )} ) \over r_{\nu } + X_{\cdot , \nu } - 1 + \tilde{\de } _{\nu }^{( \bbt , \cbt )} ( \Xt ^{( \cbt ^{( \nu )} )} )} \Big) \Big] \non \\
&= E \Big[ \sum_{\nu = 1}^{n} \sum_{i = 1}^{m_{\nu }} \Big( c_{i, \nu } {X_{i, \nu } + 1 \over r_{\nu } + X_{\cdot , \nu }} \Big[ {\{ \bt _{\nu } + 1 / ( \Xt ^{( \cbt ^{( \nu )} )} + \ct _{i, \nu }^{( \nu )} ) \} ^2 \over \{ r_{\nu } + X_{\cdot , \nu } + \bt _{\nu } + 1 / ( \Xt ^{( \cbt ^{( \nu )} )} + \ct _{i, \nu }^{( \nu )} ) \} ^2} \non \\
&\quad - 2 {\bt _{\nu } + 1 / ( \Xt ^{( \cbt ^{( \nu )} )} + \ct _{i, \nu }^{( \nu )} ) \over r_{\nu } + X_{\cdot , \nu } + \bt _{\nu } + 1 / ( \Xt ^{( \cbt ^{( \nu )} )} + \ct _{i, \nu }^{( \nu )} )} \Big] + 2 c_{i, \nu } {\ph _{i, \nu }^{\rm{U}} \tilde{\de } _{\nu }^{( \bbt , \cbt )} ( \Xt ^{( \cbt ^{( \nu )} )} ) \over r_{\nu } + X_{\cdot , \nu } - 1 + \tilde{\de } _{\nu }^{( \bbt , \cbt )} ( \Xt ^{( \cbt ^{( \nu )} )} )} \Big) \Big] \non \\
&\le E \Big[ \sum_{\nu = 1}^{n} \sum_{i = 1}^{m_{\nu }} \Big( c_{i, \nu } {X_{i, \nu } + 1 \over r_{\nu } + X_{\cdot , \nu }} \Big[ {\{ \bt _{\nu } + 1 / ( \Xt ^{( \cbt ^{( \nu )} )} + \overline{\ct } _{\nu }^{( \nu )} ) \} ^2 \over \{ r_{\nu } + X_{\cdot , \nu } + \bt _{\nu } + 1 / ( \Xt ^{( \cbt ^{( \nu )} )} + \overline{\ct } _{\nu }^{( \nu )} ) \} ^2} \non \\
&\quad - 2 {\bt _{\nu } + 1 / ( \Xt ^{( \cbt ^{( \nu )} )} + \overline{\ct } _{\nu }^{( \nu )} ) \over r_{\nu } + X_{\cdot , \nu } + \bt _{\nu } + 1 / ( \Xt ^{( \cbt ^{( \nu )} )} + \overline{\ct } _{\nu }^{( \nu )} )} \Big] + 2 c_{i, \nu } {\ph _{i, \nu }^{\rm{U}} \tilde{\de } _{\nu }^{( \bbt , \cbt )} ( \Xt ^{( \cbt ^{( \nu )} )} ) \over r_{\nu } + X_{\cdot , \nu } - 1 + \tilde{\de } _{\nu }^{( \bbt , \cbt )} ( \Xt ^{( \cbt ^{( \nu )} )} )} \Big) \Big] \text{.} \non 
\end{align}
Fix $(( x_{i, \nu } )_{i = 1}^{m_{\nu }} )_{\nu = 1, \dots , N} \in ( {\mathbb{N} _0}^{m_1} \times \dots \times {\mathbb{N} _0}^{m_N} ) \setminus \{ ( \bm{0} ^{( m_{\nu } )} )_{\nu = 1, \dots , N} \} $ and let $x_{\cdot , \nu } = \sum_{i = 1}^{m_{\nu }} x_{i, \nu }$ and $\xt ^{( \cbt ^{( \nu )} )} = \sum_{{\nu }' = 1}^{N} \sum_{i = 1}^{m_{{\nu }'}} \ct _{i, {\nu }'}^{( \nu )} x_{i, {\nu }'}$ for $\nu = 1, \dots , N$. 
As in the proof of Theorem \ref{thm:EB}, it is sufficient to show that $\sum_{\nu = 1}^{n} I_{\nu }^{( \bbt , \cbt )} \le 0$, where 
\begin{align}
I_{\nu }^{( \bbt , \cbt )} &= \sum_{i = 1}^{m_{\nu }} \Big( c_{i, \nu } {x_{i, \nu } + 1 \over r_{\nu } + x_{\cdot , \nu }} \Big[ {\{ \bt _{\nu } + 1 / ( \xt ^{( \cbt ^{( \nu )} )} + \overline{\ct } _{\nu }^{( \nu )} ) \} ^2 \over \{ r_{\nu } + x_{\cdot , \nu } + \bt _{\nu } + 1 / ( \xt ^{( \cbt ^{( \nu )} )} + \overline{\ct } _{\nu }^{( \nu )} ) \} ^2} \non \\
&\quad - 2 {\bt _{\nu } + 1 / ( \xt ^{( \cbt ^{( \nu )} )} + \overline{\ct } _{\nu }^{( \nu )} ) \over r_{\nu } + x_{\cdot , \nu } + \bt _{\nu } + 1 / ( \xt ^{( \cbt ^{( \nu )} )} + \overline{\ct } _{\nu }^{( \nu )} )} \Big] + {2 c_{i, \nu } x_{i, \nu } ( \bt _{\nu } + 1 / \xt ^{( \cbt ^{( \nu )} )} ) \over ( r_{\nu } + x_{\cdot , \nu } - 1) ( r_{\nu } + x_{\cdot , \nu } - 1 + \bt _{\nu } + 1 / \xt ^{( \cbt ^{( \nu )} )} )} \Big) \non 
\end{align}
for $\nu = 1, \dots , n$. 
It can be verified that for all $\nu = 1, \dots , n$, 
\begin{align}
&I_{\nu }^{( \bbt , \cbt )} \non \\
&\le {\overline{c} _{\nu } x_{\cdot , \nu } + c_{\cdot , \nu } \over r_{\nu } + x_{\cdot , \nu }} {\{ \bt _{\nu } + 1 / ( \xt ^{( \cbt ^{( \nu )} )} + \overline{\ct } _{\nu }^{( \nu )} ) \} ^2 \over \{ r_{\nu } + x_{\cdot , \nu } + \bt _{\nu } + 1 / ( \xt ^{( \cbt ^{( \nu )} )} + \overline{\ct } _{\nu }^{( \nu )} ) \} ^2} - 2 {c_{\cdot , \nu } \over r_{\nu } + x_{\cdot , \nu }} {\bt _{\nu } + 1 / ( \xt ^{( \cbt ^{( \nu )} )} + \overline{\ct } _{\nu }^{( \nu )} ) \over r_{\nu } + x_{\cdot , \nu } + \bt _{\nu } + 1 / ( \xt ^{( \cbt ^{( \nu )} )} + \overline{\ct } _{\nu }^{( \nu )} )} \non \\
&\quad - 2 {\overline{c} _{\nu } x_{\cdot , \nu } \over r_{\nu } + x_{\cdot , \nu }} {\bt _{\nu } + 1 / ( \xt ^{( \cbt ^{( \nu )} )} + \overline{\ct } _{\nu }^{( \nu )} ) \over r_{\nu } + x_{\cdot , \nu } + \bt _{\nu } + 1 / ( \xt ^{( \cbt ^{( \nu )} )} + \overline{\ct } _{\nu }^{( \nu )} )} + {2 \overline{c} _{\nu } x_{\cdot , \nu } ( \bt _{\nu } + 1 / \xt ^{( \cbt ^{( \nu )} )} ) \over ( r_{\nu } + x_{\cdot , \nu } - 1) ( r_{\nu } + x_{\cdot , \nu } - 1 + \bt _{\nu } + 1 / \xt ^{( \cbt ^{( \nu )} )} )} \text{.} \non 
\end{align}
Now for all $\nu = 1, \dots , n$ such that $\overline{c} _{\nu } x_{\cdot , \nu } > 0$, since 
\begin{align}
&x_{\cdot , \nu } ( \bt _{\nu } + 1 / \xt ^{( \cbt ^{( \nu )} )} ) - ( x_{\cdot , \nu } + \Ct _{\nu } ) \{ \bt _{\nu } + 1 / ( \xt ^{( \cbt ^{( \nu )} )} + \overline{\ct } _{\nu }^{( \nu )} ) \} \non \\
&= \overline{\ct } _{\nu }^{( \nu )} x_{\cdot , \nu } / \{ \xt ^{( \cbt ^{( \nu )} )} ( \xt ^{( \cbt ^{( \nu )} )} + \overline{\ct } _{\nu }^{( \nu )} ) \} - \Ct _{\nu } \{ \bt _{\nu } ( \xt ^{( \cbt ^{( \nu )} )} + \overline{\ct } _{\nu }^{( \nu )} ) + 1 \} / ( \xt ^{( \cbt ^{( \nu )} )} + \overline{\ct } _{\nu }^{( \nu )} ) \non \\
&\le \overline{\ct } _{\nu }^{( \nu )} x_{\cdot , \nu } / \{ \underline{\ct } _{\nu }^{( \nu )} x_{\cdot , \nu } ( \xt ^{( \cbt ^{( \nu )} )} + \overline{\ct } _{\nu }^{( \nu )} ) \} - \Ct _{\nu } \{ \bt _{\nu } ( \underline{\ct } _{\nu }^{( \nu )} + \overline{\ct } _{\nu }^{( \nu )} ) + 1 \} / ( \xt ^{( \cbt ^{( \nu )} )} + \overline{\ct } _{\nu }^{( \nu )} ) = 0 \text{,} \non 
\end{align}
it follows that 
\begin{align}
&{2 \overline{c} _{\nu } x_{\cdot , \nu } ( \bt _{\nu } + 1 / \xt ^{( \cbt ^{( \nu )} )} ) \over ( r_{\nu } + x_{\cdot , \nu } - 1) ( r_{\nu } + x_{\cdot , \nu } - 1 + \bt _{\nu } + 1 / \xt ^{( \cbt ^{( \nu )} )} )} \non \\
&\le {2 \overline{c} _{\nu } ( x_{\cdot , \nu } + \Ct _{\nu } ) \over r_{\nu } + x_{\cdot , \nu } - 1} {\bt _{\nu } + 1 / ( \xt ^{( \cbt ^{( \nu )} )} + \overline{\ct } _{\nu }^{( \nu )} ) \over r_{\nu } + x_{\cdot , \nu } - 1 + \bt _{\nu } + 1 / ( \xt ^{( \cbt ^{( \nu )} )} + \overline{\ct } _{\nu }^{( \nu )} )} \non \\
&\le {2 \overline{c} _{\nu } ( x_{\cdot , \nu } + \Ct _{\nu } + 1) \over r_{\nu } + x_{\cdot , \nu }} {\bt _{\nu } + 1 / ( \xt ^{( \cbt ^{( \nu )} )} + \overline{\ct } _{\nu }^{( \nu )} ) \over r_{\nu } + x_{\cdot , \nu } - 1 + \bt _{\nu } + 1 / ( \xt ^{( \cbt ^{( \nu )} )} + \overline{\ct } _{\nu }^{( \nu )} )} \non \\
&\le {2 \overline{c} _{\nu } ( x_{\cdot , \nu } + \Ct _{\nu } + 2) \over r_{\nu } + x_{\cdot , \nu }} {\bt _{\nu } + 1 / ( \xt ^{( \cbt ^{( \nu )} )} + \overline{\ct } _{\nu }^{( \nu )} ) \over r_{\nu } + x_{\cdot , \nu } + \bt _{\nu } + 1 / ( \xt ^{( \cbt ^{( \nu )} )} + \overline{\ct } _{\nu }^{( \nu )} )} \label{tEB_affinep0} 
\end{align}
by assumption. 
Therefore, letting $x_{\cdot , \cdot } = \sum_{\nu = 1}^{N} x_{\cdot , \nu }$ and noting that $\underline{c_{\cdot }} - A \ge 0$, we have for all $\nu = 1, \dots , n$, 
\begin{align}
I_{\nu }^{( \bbt , \cbt )} &\le {\overline{c} _{\nu } x_{\cdot , \nu } + c_{\cdot , \nu } \over r_{\nu } + x_{\cdot , \nu }} {\{ \bt _{\nu } + 1 / ( \xt ^{( \cbt ^{( \nu )} )} + \overline{\ct } _{\nu }^{( \nu )} ) \} ^2 \over \{ r_{\nu } + x_{\cdot , \nu } + \bt _{\nu } + 1 / ( \xt ^{( \cbt ^{( \nu )} )} + \overline{\ct } _{\nu }^{( \nu )} ) \} ^2} \non \\
&\quad + 2 {\overline{c} _{\nu } ( \Ct _{\nu } + 2) - c_{\cdot , \nu } \over r_{\nu } + x_{\cdot , \nu }} {\bt _{\nu } + 1 / ( \xt ^{( \cbt ^{( \nu )} )} + \overline{\ct } _{\nu }^{( \nu )} ) \over r_{\nu } + x_{\cdot , \nu } + \bt _{\nu } + 1 / ( \xt ^{( \cbt ^{( \nu )} )} + \overline{\ct } _{\nu }^{( \nu )} )} \non \\
&= {1 \over r_{\nu } + x_{\cdot , \nu }} {\bt _{\nu } + 1 / ( \xt ^{( \cbt ^{( \nu )} )} + \overline{\ct } _{\nu }^{( \nu )} ) \over \{ r_{\nu } + x_{\cdot , \nu } + \bt _{\nu } + 1 / ( \xt ^{( \cbt ^{( \nu )} )} + \overline{\ct } _{\nu }^{( \nu )} ) \} ^2} \non \\
&\quad \times [( \overline{c} _{\nu } x_{\cdot , \nu } + c_{\cdot , \nu } ) \{ \bt _{\nu } + 1 / ( \xt ^{( \cbt ^{( \nu )} )} + \overline{\ct } _{\nu }^{( \nu )} ) \} \non \\
&\quad - 2 \{ c_{\cdot , \nu } - \overline{c} _{\nu } ( \Ct _{\nu } + 2) \} \{ r_{\nu } + x_{\cdot , \nu } + \bt _{\nu } + 1 / ( \xt ^{( \cbt ^{( \nu )} )} + \overline{\ct } _{\nu }^{( \nu )} ) \} ] \non \\
&\le {1 \over r_{\nu } + x_{\cdot , \nu }} {\bt _{\nu } + 1 / ( \xt ^{( \cbt ^{( \nu )} )} + \overline{\ct } _{\nu }^{( \nu )} ) \over \{ r_{\nu } + x_{\cdot , \nu } + \bt _{\nu } + 1 / ( \xt ^{( \cbt ^{( \nu )} )} + \overline{\ct } _{\nu }^{( \nu )} ) \} ^2} \non \\
&\quad \times [( \overline{\overline{c}} x_{\cdot , \nu } + \underline{c_{\cdot }} ) \{ \bt _{\nu } + 1 / ( \xt ^{( \cbt ^{( \nu )} )} + \overline{\ct } _{\nu }^{( \nu )} ) \} - 2 ( \underline{c_{\cdot }} - A) \{ r_{\nu } + x_{\cdot , \nu } + \bt _{\nu } + 1 / ( \xt ^{( \cbt ^{( \nu )} )} + \overline{\ct } _{\nu }^{( \nu )} ) \} ] \non \\
&\le {\overline{\overline{c}} x_{\cdot , \nu } + \underline{c_{\cdot }} \over r_{\nu } + x_{\cdot , \nu }} {\{ \overline{\bt } + 1 / ( \ct _{*} x_{\cdot , \cdot } + \underline{\overline{\ct }} ) \} ^2 \over \{ r_{\nu } + x_{\cdot , \nu } + \overline{\bt } + 1 / ( \ct _{*} x_{\cdot , \cdot } + \underline{\overline{\ct }} ) \} ^2} - 2 {\underline{c_{\cdot }} - A \over r_{\nu } + x_{\cdot , \nu }} {\underline{\bt } + 1 / ( \ct ^{*} x_{\cdot , \cdot } + \overline{\overline{\ct }} ) \over r_{\nu } + x_{\cdot , \nu } + \underline{\bt } + 1 / ( \ct ^{*} x_{\cdot , \cdot } + \overline{\overline{\ct }} )} \text{,} \non 
\end{align}
which implies that 
\begin{align}
I_{\nu }^{( \bbt , \cbt )} &\le {\overline{\overline{c}} x_{\cdot , \nu } + \underline{c_{\cdot }} \over \underline{r} + x_{\cdot , \nu }} {\{ \overline{\bt } + 1 / ( \ct _{*} x_{\cdot , \cdot } + \underline{\overline{\ct }} ) \} ^2 \over \{ \underline{r} + x_{\cdot , \nu } + \overline{\bt } + 1 / ( \ct _{*} x_{\cdot , \cdot } + \underline{\overline{\ct }} ) \} ^2} - 2 {\underline{c_{\cdot }} - A \over \overline{r} + x_{\cdot , \nu }} {\underline{\bt } + 1 / ( \ct ^{*} x_{\cdot , \cdot } + \overline{\overline{\ct }} ) \over \overline{r} + x_{\cdot , \nu } + \underline{\bt } + 1 / ( \ct ^{*} x_{\cdot , \cdot } + \overline{\overline{\ct }} )} \non \\
&\le {1 \over \underline{r} + x_{\cdot , \nu }} {\overline{\bt } + 1 / ( \ct _{*} x_{\cdot , \cdot } + \underline{\overline{\ct }} ) \over \{ \underline{r} + x_{\cdot , \nu } + \overline{\bt } + 1 / ( \ct _{*} x_{\cdot , \cdot } + \underline{\overline{\ct }} ) \} ^2} \non \\
&\quad \times \Big[ ( \overline{\overline{c}} x_{\cdot , \nu } + \underline{c_{\cdot }} ) \{ \overline{\bt } + 1 / ( \ct _{*} x_{\cdot , \cdot } + \underline{\overline{\ct }} ) \} - 2 \Big( {\underline{r} \over \overline{r}} \Big) ^2 ( \underline{c_{\cdot }} - A) {\underline{\bt } \ct _{*} \underline{\overline{\ct }} \over \overline{\bt } \ct ^{*} \overline{\overline{\ct }}} \{ \underline{r} + x_{\cdot , \nu } + \overline{\bt } + 1 / ( \ct _{*} x_{\cdot , \cdot } + \underline{\overline{\ct }} ) \} \Big] \non \\
&= {1 \over \underline{r} + x_{\cdot , \nu }} {\overline{\bt } + 1 / ( \ct _{*} x_{\cdot , \cdot } + \underline{\overline{\ct }} ) \over \{ \underline{r} + x_{\cdot , \nu } + \overline{\bt } + 1 / ( \ct _{*} x_{\cdot , \cdot } + \underline{\overline{\ct }} ) \} ^2} \non \\
&\quad \times \Big( x_{\cdot , \nu } \Big[ \overline{\overline{c}} \{ \overline{\bt } + 1 / ( \ct _{*} x_{\cdot , \cdot } + \underline{\overline{\ct }} ) \} - 2 \Big( {\underline{r} \over \overline{r}} \Big) ^2 ( \underline{c_{\cdot }} - A) {\underline{\bt } \ct _{*} \underline{\overline{\ct }} \over \overline{\bt } \ct ^{*} \overline{\overline{\ct }}} \Big] \non \\
&\quad + \underline{c_{\cdot }} \{ \overline{\bt } + 1 / ( \ct _{*} x_{\cdot , \cdot } + \underline{\overline{\ct }} ) \} - 2 \Big( {\underline{r} \over \overline{r}} \Big) ^2 ( \underline{c_{\cdot }} - A) {\underline{\bt } \ct _{*} \underline{\overline{\ct }} \over \overline{\bt } \ct ^{*} \overline{\overline{\ct }}} \{ \underline{r} + \overline{\bt } + 1 / ( \ct _{*} x_{\cdot , \cdot } + \underline{\overline{\ct }} ) \} \Big) \label{tEB_affinep1} 
\end{align}
and that 
\begin{align}
I_{\nu }^{( \bbt , \cbt )} &\le {\overline{\overline{c}} ( r_{\nu } + x_{\cdot , \nu } ) + \underline{c_{\cdot }} - \overline{\overline{c}} \underline{r} \over r_{\nu } + x_{\cdot , \nu }} {\{ \overline{\bt } + 1 / ( \ct _{*} x_{\cdot , \cdot } + \underline{\overline{\ct }} ) \} ^2 \over \{ r_{\nu } + x_{\cdot , \nu } + \overline{\bt } + 1 / ( \ct _{*} x_{\cdot , \cdot } + \underline{\overline{\ct }} ) \} ^2} - 2 {\underline{c_{\cdot }} - A \over r_{\nu } + x_{\cdot , \nu }} {\underline{\bt } + 1 / ( \ct ^{*} x_{\cdot , \cdot } + \overline{\overline{\ct }} ) \over r_{\nu } + x_{\cdot , \nu } + \underline{\bt } + 1 / ( \ct ^{*} x_{\cdot , \cdot } + \overline{\overline{\ct }} )} \non \\
&\le {1 \over r_{\nu } + x_{\cdot , \nu }} {\overline{\bt } + 1 / ( \ct _{*} x_{\cdot , \cdot } + \underline{\overline{\ct }} ) \over \{ r_{\nu } + x_{\cdot , \nu } + \overline{\bt } + 1 / ( \ct _{*} x_{\cdot , \cdot } + \underline{\overline{\ct }} ) \} ^2} \non \\
&\quad \times \Big[ \{ \overline{\overline{c}} ( r_{\nu } + x_{\cdot , \nu } ) + \underline{c_{\cdot }} - \overline{\overline{c}} \underline{r} \} \{ \overline{\bt } + 1 / ( \ct _{*} x_{\cdot , \cdot } + \underline{\overline{\ct }} ) \} - 2 ( \underline{c_{\cdot }} - A) {\underline{\bt } \ct _{*} \underline{\overline{\ct }} \over \overline{\bt } \ct ^{*} \overline{\overline{\ct }}} \{ r_{\nu } + x_{\cdot , \nu } + \overline{\bt } + 1 / ( \ct _{*} x_{\cdot , \cdot } + \underline{\overline{\ct }} ) \} \Big] \non \\
&= {1 \over r_{\nu } + x_{\cdot , \nu }} {\overline{\bt } + 1 / ( \ct _{*} x_{\cdot , \cdot } + \underline{\overline{\ct }} ) \over \{ r_{\nu } + x_{\cdot , \nu } + \overline{\bt } + 1 / ( \ct _{*} x_{\cdot , \cdot } + \underline{\overline{\ct }} ) \} ^2} \non \\
&\quad \times \Big( ( r_{\nu } + x_{\cdot , \nu } ) \Big[ \overline{\overline{c}} \{ \overline{\bt } + 1 / ( \ct _{*} x_{\cdot , \cdot } + \underline{\overline{\ct }} ) \} - 2 ( \underline{c_{\cdot }} - A) {\underline{\bt } \ct _{*} \underline{\overline{\ct }} \over \overline{\bt } \ct ^{*} \overline{\overline{\ct }}} \Big] \non \\
&\quad + ( \underline{c_{\cdot }} - \overline{\overline{c}} \underline{r} ) \{ \overline{\bt } + 1 / ( \ct _{*} x_{\cdot , \cdot } + \underline{\overline{\ct }} ) \} - 2 ( \underline{c_{\cdot }} - A) {\underline{\bt } \ct _{*} \underline{\overline{\ct }} \over \overline{\bt } \ct ^{*} \overline{\overline{\ct }}} \{ \overline{\bt } + 1 / ( \ct _{*} x_{\cdot , \cdot } + \underline{\overline{\ct }} ) \} \Big) \text{.} \label{tEB_affinep2} 
\end{align}
By (\ref{tEB_affinep1}) and (\ref{tEB_affinep2}) and by the covariance inequality, we conclude as in the proof of Theorem \ref{thm:EB} that $\sum_{\nu = 1}^{n} I_{\nu }^{( \bbt , \cbt )} \le 0$. 
\hfill$\Box$

\begin{remark}
\label{rem:1r} 
Suppose that $m_1 = \dots = m_N$, that $r_1 = \dots = r_N$, and that $\c = ( \j ^{( m_{\nu } )} )_{\nu = 1, \dots , N}$. 
Then, by modifying the above proof, we can show that if $r_1 \ge 1$, the UMVU estimator is dominated by an empirical Bayes estimator for sufficiently large $m_1$, which is related to the problem of Section 5.1 of Hamura and Kubokawa (2020b). 
For example, the empirical Bayes estimator (\ref{eq:EB_ML}) with $\mathring{\a } = \j ^{(N)}$ corresponds to $\bbt = m_1 \j ^{(N)}$ and $\cbt = (((1 / (N m_1 r_1 ))_{i = 1}^{m_{{\nu }'}} )_{{\nu }' = 1, \dots , N} )_{\nu = 1}^{N}$. 
In this case, 
\begin{align}
I_{\nu }^{( \bbt , \cbt )} &= {x_{\cdot , \nu } + m_1 \over r_1 + x_{\cdot , \nu }} \Big[ \Big\{ {m_1 + N m_1 r_1 / ( x_{\cdot , \cdot } + 1) \over r_1 + x_{\cdot , \nu } + m_1 + N m_1 r_1 / ( x_{\cdot , \cdot } + 1)} \Big\} ^2 \non \\
&\quad - 2 {m_1 + N m_1 r_1 / ( x_{\cdot , \cdot } + 1) \over r_1 + x_{\cdot , \nu } + m_1 + N m_1 r_1 / ( x_{\cdot , \cdot } + 1)} \Big] + {2 x_{\cdot , \nu } ( m_1 + N m_1 r_1 / x_{\cdot , \cdot } ) \over ( r_1 + x_{\cdot , \nu } - 1) ( r_1 + x_{\cdot , \nu } - 1 + m_1 + N m_1 r_1 / x_{\cdot , \cdot } )} \non 
\end{align}
for $\nu = 1, \dots , n$. 
Now suppose that $r_1 \ge 1$ and that $r_1 + m_1 \ge 4$. 
Then for all $\nu = 1, \dots , n$ such that $x_{\cdot , \nu } \ge 1$, (\ref{tEB_affinep0}) can be replaced by 
\begin{align}
&{2 x_{\cdot , \nu } ( m_1 + N m_1 r_1 / x_{\cdot , \cdot } ) \over ( r_1 + x_{\cdot , \nu } - 1) ( r_1 + x_{\cdot , \nu } - 1 + m_1 + N m_1 r_1 / x_{\cdot , \cdot } )} \non \\
&\le {2 ( x_{\cdot , \nu } + 1) \over r_1 + x_{\cdot , \nu }} {m_1 + N m_1 r_1 / x_{\cdot , \cdot } \over r_1 + x_{\cdot , \nu } - 1 + m_1 + N m_1 r_1 / x_{\cdot , \cdot }} \non \\
&\le {2 ( x_{\cdot , \nu } + 3) \over r_1 + x_{\cdot , \nu }} {m_1 + N m_1 r_1 / ( x_{\cdot , \cdot } + 1) \over r_1 + x_{\cdot , \nu } - 1 + m_1 + N m_1 r_1 / x_{\cdot , \cdot }} \non \\
&\le {2 ( x_{\cdot , \nu } + 4) \over r_1 + x_{\cdot , \nu }} {m_1 + N m_1 r_1 / ( x_{\cdot , \cdot } + 1) \over r_1 + x_{\cdot , \nu } + m_1 + N m_1 r_1 / x_{\cdot , \cdot }} \text{,} \non 
\end{align}
where the second inequality holds even if $x_{\cdot , \nu } = x_{\cdot , \cdot }$ since $x_{\cdot , \cdot } \ge 1$. 
This leads to a dominance condition which is satisfied when $m_1$ is sufficiently large. 
\end{remark}

\noindent
{\bf Proof of Lemma \ref{lem:likelihood}.} \ \ We have 
\begin{align}
{f( \w | \p ) \over C( \w )} &= \prod_{\la = 1}^{L} \prod_{\i = ( i_h )_{h = 1}^{d^{( \la )}} \in I_{0}^{( \la )}} \Big\{ \prod_{h = 1}^{d^{( \la )}} p_{i_h , \nu _{h}^{( \la )}} \Big\} ^{w_{\i }^{( \la )}} = \prod_{\la = 1}^{L} \prod_{h = 1}^{d^{( \la )}} \prod_{\i = ( i_h )_{h = 1}^{d^{( \la )}} \in I_{0}^{( \la )}} {p_{i_h , \nu _{h}^{( \la )}}}^{w_{\i }^{( \la )}} \non \\
&= \prod_{\nu = 1}^{N} \prod_{i = 0}^{m_{\nu }} \prod_{\la \in \La ( \nu )} \prod_{\i \in I_{0}^{( \la )} (i, \nu )} {p_{i, \nu }}^{w_{\i }^{( \la )}} = \prod_{\nu = 1}^{N} \prod_{i = 0}^{m_{\nu }} {p_{i, \nu }}^{\sum_{\la \in \La ( \nu )} \sum_{\i \in I_{0}^{( \la )} (i, \nu )} w_{\i }^{( \la )}} \text{,} \non 
\end{align}
which is the desired result. 
\hfill$\Box$

\bigskip

\noindent
{\bf Proof of Theorem \ref{thm:multin}.} \ \ In this proof, if $\varphi $ is a continuous function from $(0, \infty )$ to $[0, \infty )$, we write 
\begin{align}
&\int_{0}^{\infty } d\mu (u) = \int_{0}^{\infty } u^{\al - 1} e^{- \be u} \Big\{ \prod_{\nu = 1}^{N} {\Ga ( \ga _{\nu } u + r_{\nu } + a_{0, \nu } ) \Ga ( r_{\nu } + a_{0, \nu } + X_{\cdot , \nu } + a_{\cdot , \nu } ) \over \Ga ( \ga _{\nu } u + r_{\nu } + a_{0, \nu } + X_{\cdot , \nu } + a_{\cdot , \nu } ) \Ga ( r_{\nu } + a_{0, \nu } )} \Big\} du \text{,} \non \\
&\int_{0}^{\infty } \varphi (u) d\mu (u) = \int_{0}^{\infty } \varphi (u) u^{\al - 1} e^{- \be u} \Big\{ \prod_{\nu = 1}^{N} {\Ga ( \ga _{\nu } u + r_{\nu } + a_{0, \nu } ) \Ga ( r_{\nu } + a_{0, \nu } + X_{\cdot , \nu } + a_{\cdot , \nu } ) \over \Ga ( \ga _{\nu } u + r_{\nu } + a_{0, \nu } + X_{\cdot , \nu } + a_{\cdot , \nu } ) \Ga ( r_{\nu } + a_{0, \nu } )} \Big\} du \text{,} \quad \text{and} \non \\
&E^U [ \varphi (U) ] = \int_{0}^{\infty } \varphi (u) d\mu (u) / \int_{0}^{\infty } d\mu (u) \text{.} \non 
\end{align}
Let $\De ^{( \al , \be , \bga , \a _0 , \a )} %
= E[ \log \{ f( \W | \p ) / \fh ^{( \pi _{\al , \be , \bga , \a _0 , \a } )} ( \W ; \X ) \} ] - E[ \log \{ f( \W | \p ) / \fh ^{( \pi _{\a _0 , \a } )} ( \W ; \X ) \} ]$. 
Then, by Proposition \ref{prp:mass_estimator}, 
\begin{align}
&\De ^{( \al , \be , \bga , \a _0 , \a )} = E \Big[ - \log {\fh ^{( \pi _{\al , \be , \bga , \a _0 , \a } )} ( \W ; \X ) \over \fh ^{( \pi _{\a _0 , \a } )} ( \W ; \X )} \Big] \non \\
&= E \Big[ - \log E^U \Big[ \prod_{\nu = 1}^{N} \Big\{ {\Ga ( \ga _{\nu } U + s_{0, \nu } ( \W ) + r_{\nu } + a_{0, \nu } ) \Ga \big( \sum_{\la \in \La ( \nu )} l^{( \la )} + r_{\nu } + a_{0, \nu } + X_{\cdot , \nu } + a_{\cdot , \nu } \big) \over \Ga \big( \ga _{\nu } U + \sum_{\la \in \La ( \nu )} l^{( \la )} + r_{\nu } + a_{0, \nu } + X_{\cdot , \nu } + a_{\cdot , \nu } \big) \Ga ( s_{0, \nu } ( \W ) + r_{\nu } + a_{0, \nu } )} \non \\
&\quad \times {\Ga ( \ga _{\nu } U + r_{\nu } + a_{0, \nu } + X_{\cdot , \nu } + a_{\cdot , \nu } ) \Ga ( r_{\nu } + a_{0, \nu } ) \over \Ga ( \ga _{\nu } U + r_{\nu } + a_{0, \nu } ) \Ga ( r_{\nu } + a_{0, \nu } + X_{\cdot , \nu } + a_{\cdot , \nu } )} \Big\} \Big] \Big] \text{.} \label{tmultinp1} 
\end{align}

For $\nu = 1, \dots , N$, let $\pt _{0, \nu } = p_{0, \nu }$ and $\pt _{1, \nu } = p_{\cdot , \nu } = \sum_{i = 1}^{m_{\nu }} p_{i, \nu }$ for notational convenience. 
For $\la = 1, \dots , L$, let $\widetilde{\mathcal{W}} ^{( \la )} = \big\{ ( \mathring{w} _{\tilde{\i }} )_{\tilde{\i } \in \{ 0, 1 \} ^{d^{( \la )}}} \big| \mathring{w} _{\tilde{\i }} \in \mathbb{N} _0 \quad \text{for all $\tilde{\i } \in \{ 0, 1 \} ^{d^{( \la )}}$} \quad \text{and} \quad  \sum_{\tilde{\i } \in \{ 0, 1 \} ^{d^{( \la )}}} \mathring{w} _{\tilde{\i }} = 1 \big\} $. 
Let $\Wbt ^{( \la )} (j) = ( \Wt _{\tilde{\i } }^{( \la )} (j) )_{\tilde{\i } \in \{ 0, 1 \} ^{d^{( \la )}}}$, $j = 1, \dots , l^{( \la )}$, $\la = 1, \dots , L$, be independent multinomial random variables with mass functions 
\begin{align}
\prod_{\tilde{\i } = ( \tilde{i} _h )_{h = 1}^{d^{( \la )}} \in \{ 0, 1 \} ^{d^{( \la )}}} \Big\{ \prod_{h = 1}^{d^{( \la )}} \pt _{\tilde{i} _h , \nu _{h}^{( \la )}} \Big\} ^{\wt _{\tilde{\i } }^{( \la )} (j)} \text{,} \non 
\end{align}
$( \wt _{\tilde{\i } }^{( \la )} (j) )_{\tilde{\i } \in \{ 0, 1 \} ^{d^{( \la )}}} \in \widetilde{\mathcal{W}} ^{( \la )}$, $j = 1, \dots , l^{( \la )}$, $\la = 1, \dots , L$, respectively. 
For  $\nu = 1, \dots , N$, let $\It _{0}^{( \la )} ( \nu ) = \It _{0}^{( \la )} (0, \nu ) = \{ ( \tilde{i} _h )_{h = 1}^{d^{( \la )}} \in \{ 0, 1 \} ^{d^{( \la )}} | \tilde{i} _{h_{\nu }^{( \la )}} = 0 \} $ for $\la \in \La ( \nu )$. 
Notice that 
\begin{align}
\Big( \Big( \sum_{\i \in I_{0}^{( \la )} (0, \nu )} W_{\i }^{( \la )} \Big) _{\la \in \La ( \nu )} \Big) _{\nu = 1, \dots , N} \stackrel{\rm{d}}{=} \Big( \Big( \sum_{\tilde{\i } \in \It _{0}^{( \la )} ( \nu )} \sum_{j = 1}^{l^{( \la )}} \Wt _{\tilde{\i }}^{( \la )} (j) \Big) _{\la \in \La ( \nu )} \Big) _{\nu = 1, \dots , N} \text{.} \label{tmultinp2} 
\end{align}
Then it follows from (\ref{tmultinp1}) and (\ref{tmultinp2}) that 
\begin{align}
\De ^{( \al , \be , \bga , \a _0 , \a )} &= E \Big[ - \log E^U \Big[ \prod_{\nu = 1}^{N} \Big\{ {\Ga \big( \ga _{\nu } U + \sum_{\la \in \La ( \nu )} \sum_{\i \in I_{0}^{( \la )} (0, \nu )} W_{\i }^{( \la )} + r_{\nu } + a_{0, \nu } \big) \over \Ga \big( \ga _{\nu } U + \sum_{\la \in \La ( \nu )} l^{( \la )} + r_{\nu } + a_{0, \nu } + X_{\cdot , \nu } + a_{\cdot , \nu } \big) } \non \\
&\quad \times {\Ga \big( \sum_{\la \in \La ( \nu )} l^{( \la )} + r_{\nu } + a_{0, \nu } + X_{\cdot , \nu } + a_{\cdot , \nu } \big) \over \Ga \big( \sum_{\la \in \La ( \nu )} \sum_{\i \in I_{0}^{( \la )} (0, \nu )} W_{\i }^{( \la )} + r_{\nu } + a_{0, \nu } \big) } \non \\
&\quad \times {\Ga ( \ga _{\nu } U + r_{\nu } + a_{0, \nu } + X_{\cdot , \nu } + a_{\cdot , \nu } ) \Ga ( r_{\nu } + a_{0, \nu } ) \over \Ga ( \ga _{\nu } U + r_{\nu } + a_{0, \nu } ) \Ga ( r_{\nu } + a_{0, \nu } + X_{\cdot , \nu } + a_{\cdot , \nu } )} \Big\} \Big] \Big] \non \\
&= E \Big[ - \log E^U \Big[ \prod_{\nu = 1}^{N} \Big\{ {\Ga \big( \ga _{\nu } U + \sum_{\la \in \La ( \nu )} \sum_{\tilde{\i } \in \It _{0}^{( \la )} ( \nu )} \sum_{j = 1}^{l^{( \la )}} \Wt _{\tilde{\i }}^{( \la )} (j) + r_{\nu } + a_{0, \nu } \big) \over \Ga \big( \ga _{\nu } U + \sum_{\la \in \La ( \nu )} l^{( \la )} + r_{\nu } + a_{0, \nu } + X_{\cdot , \nu } + a_{\cdot , \nu } \big) } \non \\
&\quad \times {\Ga \big( \sum_{\la \in \La ( \nu )} l^{( \la )} + r_{\nu } + a_{0, \nu } + X_{\cdot , \nu } + a_{\cdot , \nu } \big) \over \Ga \big( \sum_{\la \in \La ( \nu )} \sum_{\tilde{\i } \in \It _{0}^{( \la )} ( \nu )} \sum_{j = 1}^{l^{( \la )}} \Wt _{\tilde{\i }}^{( \la )} (j) + r_{\nu } + a_{0, \nu } \big) } \non \\
&\quad \times {\Ga ( \ga _{\nu } U + r_{\nu } + a_{0, \nu } + X_{\cdot , \nu } + a_{\cdot , \nu } ) \Ga ( r_{\nu } + a_{0, \nu } ) \over \Ga ( \ga _{\nu } U + r_{\nu } + a_{0, \nu } ) \Ga ( r_{\nu } + a_{0, \nu } + X_{\cdot , \nu } + a_{\cdot , \nu } )} \Big\} \Big] \Big] \text{.} \non 
\end{align}
Therefore, 
\begin{align}
\De ^{( \al , \be , \bga , \a _0 , \a )} &= \sum_{((( \wt _{\tilde{\i } }^{( \la )} (j) )_{\tilde{\i } \in \{ 0, 1 \} ^{d^{( \la )}}} )_{j = 1, \dots , l^{( \la )}} )_{\la = 1, \dots , L} \in ( \widetilde{\mathcal{W}} ^{(1)} \times \dots \times \widetilde{\mathcal{W}} ^{(1)} ) \times \dots \times ( \widetilde{\mathcal{W}} ^{(L)} \times \dots \times \widetilde{\mathcal{W}} ^{(L)} )} \Big( \non \\
&\quad \Big[ \prod_{\la = 1}^{L} \prod_{j = 1}^{l^{( \la )}} \prod_{\tilde{\i } = ( \tilde{i} _h )_{h = 1}^{d^{( \la )}} \in \{ 0, 1 \} ^{d^{( \la )}}} \Big\{ \prod_{h = 1}^{d^{( \la )}} \pt _{\tilde{i} _h , \nu _{h}^{( \la )}} \Big\} ^{\wt _{\tilde{\i } }^{( \la )} (j)} \Big] \non \\
&\quad \times E \Big[ - \log E^U \Big[ \prod_{\nu = 1}^{N} \Big\{ {\Ga \big( \ga _{\nu } U + \sum_{\la \in \La ( \nu )} \sum_{\tilde{\i } \in \It _{0}^{( \la )} ( \nu )} \sum_{j = 1}^{l^{( \la )}} \wt _{\tilde{\i }}^{( \la )} (j) + r_{\nu } + a_{0, \nu } \big) \over \Ga \big( \ga _{\nu } U + \sum_{\la \in \La ( \nu )} l^{( \la )} + r_{\nu } + a_{0, \nu } + X_{\cdot , \nu } + a_{\cdot , \nu } \big) } \non \\
&\quad \times {\Ga \big( \sum_{\la \in \La ( \nu )} l^{( \la )} + r_{\nu } + a_{0, \nu } + X_{\cdot , \nu } + a_{\cdot , \nu } \big) \over \Ga \big( \sum_{\la \in \La ( \nu )} \sum_{\tilde{\i } \in \It _{0}^{( \la )} ( \nu )} \sum_{j = 1}^{l^{( \la )}} \wt _{\tilde{\i }}^{( \la )} (j) + r_{\nu } + a_{0, \nu } \big) } \non \\
&\quad \times {\Ga ( \ga _{\nu } U + r_{\nu } + a_{0, \nu } + X_{\cdot , \nu } + a_{\cdot , \nu } ) \Ga ( r_{\nu } + a_{0, \nu } ) \over \Ga ( \ga _{\nu } U + r_{\nu } + a_{0, \nu } ) \Ga ( r_{\nu } + a_{0, \nu } + X_{\cdot , \nu } + a_{\cdot , \nu } )} \Big\} \Big] \Big] \Big) \non \\
&= \sum_{\tilde{i} _{1}^{(1)} (1) = 0}^{1} \pt _{\tilde{i} _{1}^{(1)} (1), \nu _{1}^{(1)}} \dotsm \sum_{\tilde{i} _{d^{(1)}}^{(1)} (1) = 0}^{1} \pt _{\tilde{i} _{d^{(1)}}^{(1)} (1), \nu _{d^{(1)}}^{(1)}} \non \\
&\quad \dotsm \sum_{\tilde{i} _{1}^{(1)} ( l^{(1)} ) = 0}^{1} \pt _{\tilde{i} _{1}^{(1)} ( l^{(1)} ), \nu _{1}^{(1)}} \dotsm \sum_{\tilde{i} _{d^{(1)}}^{(1)} ( l^{(1)} ) = 0}^{1} \pt _{\tilde{i} _{d^{(1)}}^{(1)} ( l^{(1)} ), \nu _{d^{(1)}}^{(1)}} \non \\
&\quad \dotsm \non \\
&\quad \sum_{\tilde{i} _{1}^{(L)} (1) = 0}^{1} \pt _{\tilde{i} _{1}^{(L)} (1), \nu _{1}^{(L)}} \dotsm \sum_{\tilde{i} _{d^{(L)}}^{(L)} (1) = 0}^{1} \pt _{\tilde{i} _{d^{(L)}}^{(L)} (1), \nu _{d^{(L)}}^{(L)}} \non \\
&\quad \dotsm \sum_{\tilde{i} _{1}^{(L)} ( l^{(L)} ) = 0}^{1} \pt _{\tilde{i} _{1}^{(L)} ( l^{(L)} ), \nu _{1}^{(L)}} \dotsm \sum_{\tilde{i} _{d^{(L)}}^{(L)} ( l^{(L)} ) = 0}^{1} \pt _{\tilde{i} _{d^{(L)}}^{(L)} ( l^{(L)} ), \nu _{d^{(L)}}^{(L)}} E \Big[ \non \\
&\quad - \log E^U \Big[ \prod_{\nu = 1}^{N} \Big\{ {\Ga \big( \ga _{\nu } U + \sum_{\la \in \La ( \nu )} \sum_{\tilde{\i } \in \It _{0}^{( \la )} ( \nu )} \sum_{j = 1}^{l^{( \la )}} \tilde{\de } ^{( \la )} ( \tilde{\i } , ( \tilde{i} _{h}^{( \la )} (j))_{h = 1}^{d^{( \la )}} ) + r_{\nu } + a_{0, \nu } \big) \over \Ga \big( \ga _{\nu } U + \sum_{\la \in \La ( \nu )} l^{( \la )} + r_{\nu } + a_{0, \nu } + X_{\cdot , \nu } + a_{\cdot , \nu } \big) } \non \\
&\quad \times {\Ga \big( \sum_{\la \in \La ( \nu )} l^{( \la )} + r_{\nu } + a_{0, \nu } + X_{\cdot , \nu } + a_{\cdot , \nu } \big) \over \Ga \big( \sum_{\la \in \La ( \nu )} \sum_{\tilde{\i } \in \It _{0}^{( \la )} ( \nu )} \sum_{j = 1}^{l^{( \la )}} \tilde{\de } ^{( \la )} ( \tilde{\i } , ( \tilde{i} _{h}^{( \la )} (j))_{h = 1}^{d^{( \la )}} ) + r_{\nu } + a_{0, \nu } \big) } \non \\
&\quad \times {\Ga ( \ga _{\nu } U + r_{\nu } + a_{0, \nu } + X_{\cdot , \nu } + a_{\cdot , \nu } ) \Ga ( r_{\nu } + a_{0, \nu } ) \over \Ga ( \ga _{\nu } U + r_{\nu } + a_{0, \nu } ) \Ga ( r_{\nu } + a_{0, \nu } + X_{\cdot , \nu } + a_{\cdot , \nu } )} \Big\} \Big] \Big] \text{,} \non 
\end{align}
where $\tilde{\de } ^{( \la )} ( \tilde{\i } , {\tilde{\i }}' )$ $= 1$ if $\tilde{\i } = {\tilde{\i }}'$ and $= 0$ if $\tilde{\i } \neq {\tilde{\i }}'$ for $\tilde{\i } , {\tilde{\i }}' \in \{ 0, 1 \} ^{d^{( \la )}}$ for $\la = 1, \dots , L$. 
Furthermore, since 
\begin{align}
\sum_{\la \in \La ( \nu )} \sum_{\tilde{\i } \in \It _{0}^{( \la )} ( \nu )} \sum_{j = 1}^{l^{( \la )}} \tilde{\de } ^{( \la )} ( \tilde{\i } , ( \tilde{i} _{h}^{( \la )} (j))_{h = 1}^{d^{( \la )}} ) &= \sum_{\la \in \La ( \nu )} \sum_{j = 1}^{l^{( \la )}} \{ 1 - \tilde{i} _{h_{\nu }^{( \la )}}^{( \la )} (j) \} \non 
\end{align}
for all $((( \tilde{i} _{h}^{( \la )} (j))_{h = 1}^{d^{( \la )}} )_{j = 1, \dots , l^{( \la )}} )_{\la = 1, \dots , L} \in ( \{ 0, 1 \} ^{d^{(1)}} \times \dots \times \{ 0, 1 \} ^{d^{(1)}} ) \times \dots \times ( \{ 0, 1 \} ^{d^{(L)}} \times \dots \times \{ 0, 1 \} ^{d^{(L)}} )$ for all $\nu = 1, \dots , N$, we can rewrite the risk difference as 
\begin{align}
\De ^{( \al , \be , \bga , \a _0 , \a )} &= \sum_{\tilde{i} _{1}^{(1)} (1) = 0}^{1} \pt _{\tilde{i} _{1}^{(1)} (1), \nu _{1}^{(1)}} \dotsm \sum_{\tilde{i} _{d^{(1)}}^{(1)} (1) = 0}^{1} \pt _{\tilde{i} _{d^{(1)}}^{(1)} (1), \nu _{d^{(1)}}^{(1)}} \non \\
&\quad \dotsm \sum_{\tilde{i} _{1}^{(1)} ( l^{(1)} ) = 0}^{1} \pt _{\tilde{i} _{1}^{(1)} ( l^{(1)} ), \nu _{1}^{(1)}} \dotsm \sum_{\tilde{i} _{d^{(1)}}^{(1)} ( l^{(1)} ) = 0}^{1} \pt _{\tilde{i} _{d^{(1)}}^{(1)} ( l^{(1)} ), \nu _{d^{(1)}}^{(1)}} \non \\
&\quad \dotsm \non \\
&\quad \sum_{\tilde{i} _{1}^{(L)} (1) = 0}^{1} \pt _{\tilde{i} _{1}^{(L)} (1), \nu _{1}^{(L)}} \dotsm \sum_{\tilde{i} _{d^{(L)}}^{(L)} (1) = 0}^{1} \pt _{\tilde{i} _{d^{(L)}}^{(L)} (1), \nu _{d^{(L)}}^{(L)}} \non \\
&\quad \dotsm \sum_{\tilde{i} _{1}^{(L)} ( l^{(L)} ) = 0}^{1} \pt _{\tilde{i} _{1}^{(L)} ( l^{(L)} ), \nu _{1}^{(L)}} \dotsm \sum_{\tilde{i} _{d^{(L)}}^{(L)} ( l^{(L)} ) = 0}^{1} \pt _{\tilde{i} _{d^{(L)}}^{(L)} ( l^{(L)} ), \nu _{d^{(L)}}^{(L)}} E \Big[ \non \\
&\quad - \log E^U \Big[ F \Big(U, 
((( \tilde{i} _{h}^{( \la )} (j))_{h = 1}^{d^{( \la )}} )_{j = 1, \dots , l^{( \la )}} )_{\la = 1, \dots , L} , \Big( \sum_{\la \in \La ( \nu )} l^{( \la )} \Big) _{\nu = 1}^{N} \Big) \Big] \Big] \text{,} \label{tmultinp2.05} 
\end{align}
where 
\begin{align}
F(u, \tilde{\i } , \k ) &= \prod_{\nu = 1}^{N} \Big[ {\Ga \big( \ga _{\nu } u + \sum_{\la \in \La ( \nu )} \sum_{j = 1}^{l^{( \la )}} \{ 1 - \tilde{i} _{h_{\nu }^{( \la )}}^{( \la )} (j) \} + r_{\nu } + a_{0, \nu } \big) \over \Ga ( \ga _{\nu } u + k_{\nu } + r_{\nu } + a_{0, \nu } + X_{\cdot , \nu } + a_{\cdot , \nu } ) } \non \\
&\quad \times {\Ga ( k_{\nu } + r_{\nu } + a_{0, \nu } + X_{\cdot , \nu } + a_{\cdot , \nu } ) \over \Ga \big( \sum_{\la \in \La ( \nu )} \sum_{j = 1}^{l^{( \la )}} \{ 1 - \tilde{i} _{h_{\nu }^{( \la )}}^{( \la )} (j) \} + r_{\nu } + a_{0, \nu } \big) } \non \\
&\quad \times {\Ga ( \ga _{\nu } u + r_{\nu } + a_{0, \nu } + X_{\cdot , \nu } + a_{\cdot , \nu } ) \Ga ( r_{\nu } + a_{0, \nu } ) \over \Ga ( \ga _{\nu } u + r_{\nu } + a_{0, \nu } ) \Ga ( r_{\nu } + a_{0, \nu } + X_{\cdot , \nu } + a_{\cdot , \nu } )} \Big] \non 
\end{align}
for $u \in (0, \infty )$, $\tilde{\i } = ((( \tilde{i} _{h}^{( \la )} (j))_{h = 1}^{d^{( \la )}} )_{j = 1, \dots , l^{( \la )}} )_{\la = 1, \dots , L} \in ( \{ 0, 1 \} ^{d^{(1)}} \times \dots \times \{ 0, 1 \} ^{d^{(1)}} ) \times \dots \times ( \{ 0, 1 \} ^{d^{(L)}} \times \dots \times \{ 0, 1 \} ^{d^{(L)}} )$, and $\k = ( k_{\nu } )_{\nu = 1}^{N} \in {\mathbb{N} _0}^N$. 

Now fix $\la ^{*} = 1, \dots , L$, $h^{*} = 1, \dots , d^{( \la ^{*} )}$, and $j^{*} = 1, \dots , l^{( \la ^{*} )}$. 
For each $(j, h, \la ) \in \mathbb{N} \times \mathbb{N} \times \{ 1, \dots , L \} $ satisfying $j \le l^{( \la )}$, $h \le d^{( \la )}$, and $(j, h, \la ) \neq ( j^{*} , h^{*} , \la ^{*} )$, fix $\tilde{i} _{h}^{( \la )} (j) \in \{ 0, 1 \} $. 
Let ${\nu }^{*} = \nu _{h^{*}}^{( {\la }^{*} )}$. 
For $u \in (0, \infty )$, $\tilde{i} \in \{ 0, 1 \} $, and $\k \in {\mathbb{N} _0}^N$, let $F^{*} (u, \tilde{i} , \k )$ denote $F(u, (( \tilde{i} _{h}^{( \la )} (j))_{h = 1}^{d^{( \la )}} )_{j = 1, \dots , l^{( \la )}} )_{\la = 1, \dots , L} , \k )$ with $\tilde{i} _{h^{*}}^{( \la ^{*} )} ( j^{*} ) = \tilde{i}$. 
For each $\nu = 1, \dots , N$, let $\st _{\nu }^{*} ( \tilde{i} )$ denote $\sum_{\la \in \La ( \nu )} \sum_{j = 1}^{l^{( \la )}} \{ 1 - \tilde{i} _{h_{\nu }^{( \la )}}^{( \la )} (j) \} $ with $\tilde{i} _{h^{*}}^{( \la ^{*} )} ( j^{*} ) = \tilde{i}$ for $\tilde{i} \in \{ 0, 1 \} $. 
Finally, fix $\k = ( k_{\nu } )_{\nu = 1}^{N} \in {\mathbb{N} _0}^N$ such that $\st _{\nu }^{*} ( \tilde{i} ) \le k_{\nu } \le \sum_{\la \in \La ( \nu )} l^{( \la )}$ for all $\nu = 1, \dots , N$ for any $\tilde{i} \in \{ 0, 1 \} $. 
Then, by Lemma \ref{lem:hudson}, 
\begin{align}
&\sum_{\tilde{i} = 0}^{1} \pt _{\tilde{i} , {\nu }^{*}} E[ - \log E^U [ F^{*} (U, \tilde{i} , \k ) ] ] \non \\
&= E[ - \log E^U [ F^{*} (U, 0, \k ) ] ] + \pt _{1, {\nu }^{*}} E \Big[ \log {E^U [ F^{*} (U, 0, \k ) ] \over E^U [ F^{*} (U, 1, \k ) ]} \Big] \non \\
&= E \Big[ - \log {\int_{0}^{\infty } F^{*} (u, 0, \k ) d\mu (u) \over \int_{0}^{\infty } d\mu (u)} \Big] + \pt _{1, {\nu }^{*}} E \Big[ \log {\int_{0}^{\infty } F^{*} (u, 0, \k ) d\mu (u) \over \int_{0}^{\infty } F^{*} (u, 1, \k ) d\mu (u)} \Big] \non \\
&= E \Big[ - \log {\int_{0}^{\infty } F^{*} (u, 0, \k ) d\mu (u) \over \int_{0}^{\infty } d\mu (u)} \Big] + E \Big[ {X_{\cdot , {\nu }^{*}} \over r_{{\nu }^{*}} + X_{\cdot , {\nu }^{*}} - 1} \non \\
&\quad \times \log \Big\{ \int_{0}^{\infty } F^{*} (u, 0, \k ) {\ga _{{\nu }^{*}} u + k_{{\nu }^{*}} + r_{{\nu }^{*}} + a_{0, {\nu }^{*}} + X_{\cdot , {\nu }^{*}} + a_{\cdot , {\nu }^{*}} - 1 \over k_{{\nu }^{*}} + r_{{\nu }^{*}} + a_{0, {\nu }^{*}} + X_{\cdot , {\nu }^{*}} + a_{\cdot , {\nu }^{*}} - 1} d\mu (u) \non \\
&\quad / \int_{0}^{\infty } F^{*} (u, 1, \k ) {\ga _{{\nu }^{*}} u + k_{{\nu }^{*}} + r_{{\nu }^{*}} + a_{0, {\nu }^{*}} + X_{\cdot , {\nu }^{*}} + a_{\cdot , {\nu }^{*}} - 1 \over k_{{\nu }^{*}} + r_{{\nu }^{*}} + a_{0, {\nu }^{*}} + X_{\cdot , {\nu }^{*}} + a_{\cdot , {\nu }^{*}} - 1} d\mu (u) \Big\} \Big] \text{.} \non 
\end{align}
In the following, if $\varphi $ is a continuous function from $(0, \infty )$ to $[0, \infty )$, we write 
\begin{align}
&\int_{0}^{\infty } d\mut (u) = \int_{0}^{\infty } F^{*} (u, 1, \k ) {\ga _{{\nu }^{*}} u + k_{{\nu }^{*}} + r_{{\nu }^{*}} + a_{0, {\nu }^{*}} + X_{\cdot , {\nu }^{*}} + a_{\cdot , {\nu }^{*}} - 1 \over k_{{\nu }^{*}} + r_{{\nu }^{*}} + a_{0, {\nu }^{*}} + X_{\cdot , {\nu }^{*}} + a_{\cdot , {\nu }^{*}} - 1} d\mu (u) \text{,} \non \\
&\int_{0}^{\infty } \varphi (u) d\mut (u) = \int_{0}^{\infty } \varphi (u) F^{*} (u, 1, \k ) {\ga _{{\nu }^{*}} u + k_{{\nu }^{*}} + r_{{\nu }^{*}} + a_{0, {\nu }^{*}} + X_{\cdot , {\nu }^{*}} + a_{\cdot , {\nu }^{*}} - 1 \over k_{{\nu }^{*}} + r_{{\nu }^{*}} + a_{0, {\nu }^{*}} + X_{\cdot , {\nu }^{*}} + a_{\cdot , {\nu }^{*}} - 1} d\mu (u) \text{,} \quad \text{and} \non \\
&\Et ^U [ \varphi (U) ] = \int_{0}^{\infty } \varphi (u) d\mut (u) / \int_{0}^{\infty } d\mut (u) \text{.} \non 
\end{align}
Then we have 
\begin{align}
&\sum_{\tilde{i} = 0}^{1} \pt _{\tilde{i} , {\nu }^{*}} E[ - \log E^U [ F^{*} (U, \tilde{i} , \k ) ] ] \non \\
&= E \Big[ - \log {\int_{0}^{\infty } d\mut (u) \over \int_{0}^{\infty } d\mu (u)} - \log {\int_{0}^{\infty } F^{*} (u, 0, \k ) d\mu (u) \over \int_{0}^{\infty } d\mut (u)} + {X_{\cdot , {\nu }^{*}} \over r_{{\nu }^{*}} + X_{\cdot , {\nu }^{*}} - 1} \log \Et ^U \Big[ {F^{*} (U, 0, \k ) \over F^{*} (U, 1, \k )} \Big] \Big] \non \\
&= E \Big[ - \log {\int_{0}^{\infty } d\mut (u) \over \int_{0}^{\infty } d\mu (u)} - \log \Et ^U \Big[ {F^{*} (U, 0, \k ) \over F^{*} (U, 1, \k )} {k_{{\nu }^{*}} + r_{{\nu }^{*}} + a_{0, {\nu }^{*}} + X_{\cdot , {\nu }^{*}} + a_{\cdot , {\nu }^{*}} - 1 \over \ga _{{\nu }^{*}} U + k_{{\nu }^{*}} + r_{{\nu }^{*}} + a_{0, {\nu }^{*}} + X_{\cdot , {\nu }^{*}} + a_{\cdot , {\nu }^{*}} - 1} \Big] \non \\
&\quad + {X_{\cdot , {\nu }^{*}} \over r_{{\nu }^{*}} + X_{\cdot , {\nu }^{*}} - 1} \log \Et ^U \Big[ {F^{*} (U, 0, \k ) \over F^{*} (U, 1, \k )} \Big] \Big] \text{.} \label{tmultinp2.1} 
\end{align}
Notice that for all $u \in (0, \infty )$, 
\begin{align}
{F^{*} (u, 0, \k ) \over F^{*} (u, 1, \k )} &= \prod_{\nu = 1}^{N} {\Ga ( \ga _{\nu } u + \st _{\nu }^{*} (0) + r_{\nu } + a_{0, \nu } ) \Ga ( \st _{\nu }^{*} (1) + r_{\nu } + a_{0, \nu } ) \over \Ga ( \ga _{\nu } u + \st _{\nu }^{*} (1) + r_{\nu } + a_{0, \nu } ) \Ga ( \st _{\nu }^{*} (0) + r_{\nu } + a_{0, \nu } )} \non \\
&= {\Ga ( \ga _{{\nu }^{*}} u + \st _{{\nu }^{*}}^{*} (0) + r_{{\nu }^{*}} + a_{0, {\nu }^{*}} ) \Ga ( \st _{{\nu }^{*}}^{*} (1) + r_{{\nu }^{*}} + a_{0, {\nu }^{*}} ) \over \Ga ( \ga _{{\nu }^{*}} u + \st _{{\nu }^{*}}^{*} (1) + r_{{\nu }^{*}} + a_{0, {\nu }^{*}} ) \Ga ( \st _{{\nu }^{*}}^{*} (0) + r_{{\nu }^{*}} + a_{0, {\nu }^{*}} )} \non \\
&= {\ga _{{\nu }^{*}} u + \st _{{\nu }^{*}}^{*} (1) + r_{{\nu }^{*}} + a_{0, {\nu }^{*}} \over \st _{{\nu }^{*}}^{*} (1) + r_{{\nu }^{*}} + a_{0, {\nu }^{*}}} \non 
\end{align}
since $\st _{{\nu }^{*}}^{*} (0) = \st _{{\nu }^{*}}^{*} (1) + 1$. 
It follows that 
\begin{align}
&\log \Et ^U \Big[ {F^{*} (U, 0, \k ) \over F^{*} (U, 1, \k )} {k_{{\nu }^{*}} + r_{{\nu }^{*}} + a_{0, {\nu }^{*}} + X_{\cdot , {\nu }^{*}} + a_{\cdot , {\nu }^{*}} - 1 \over \ga _{{\nu }^{*}} U + k_{{\nu }^{*}} + r_{{\nu }^{*}} + a_{0, {\nu }^{*}} + X_{\cdot , {\nu }^{*}} + a_{\cdot , {\nu }^{*}} - 1} \Big] \non \\
&= \log \Et ^U \Big[ {k_{{\nu }^{*}} + r_{{\nu }^{*}} + a_{0, {\nu }^{*}} + X_{\cdot , {\nu }^{*}} + a_{\cdot , {\nu }^{*}} - 1 \over \st _{{\nu }^{*}}^{*} (1) + r_{{\nu }^{*}} + a_{0, {\nu }^{*}}} {\ga _{{\nu }^{*}} U + \st _{{\nu }^{*}}^{*} (1) + r_{{\nu }^{*}} + a_{0, {\nu }^{*}} \over \ga _{{\nu }^{*}} U + k_{{\nu }^{*}} + r_{{\nu }^{*}} + a_{0, {\nu }^{*}} + X_{\cdot , {\nu }^{*}} + a_{\cdot , {\nu }^{*}} - 1} \Big] \non \\
&= \log \Et ^U \Big[ \Big\{ 1 + {k_{{\nu }^{*}} - \st _{{\nu }^{*}}^{*} (1) + X_{\cdot , {\nu }^{*}} + a_{\cdot , {\nu }^{*}} - 1 \over \st _{{\nu }^{*}}^{*} (1) + r_{{\nu }^{*}} + a_{0, {\nu }^{*}}} \Big\} \Big\{ 1 - {k_{{\nu }^{*}} - \st _{{\nu }^{*}}^{*} (1) + X_{\cdot , {\nu }^{*}} + a_{\cdot , {\nu }^{*}} - 1 \over \ga _{{\nu }^{*}} U + k_{{\nu }^{*}} + r_{{\nu }^{*}} + a_{0, {\nu }^{*}} + X_{\cdot , {\nu }^{*}} + a_{\cdot , {\nu }^{*}} - 1} \Big\} \Big] \non \\
&= \log \Et ^U \Big[ 1 + {\{ k_{{\nu }^{*}} - \st _{{\nu }^{*}}^{*} (1) + X_{\cdot , {\nu }^{*}} + a_{\cdot , {\nu }^{*}} - 1 \} \ga _{{\nu }^{*}} U \over \{ \st _{{\nu }^{*}}^{*} (1) + r_{{\nu }^{*}} + a_{0, {\nu }^{*}} \} ( \ga _{{\nu }^{*}} U + k_{{\nu }^{*}} + r_{{\nu }^{*}} + a_{0, {\nu }^{*}} + X_{\cdot , {\nu }^{*}} + a_{\cdot , {\nu }^{*}} - 1)} \Big] \label{tmultinp2.2} 
\end{align}
and that 
\begin{align}
&{X_{\cdot , {\nu }^{*}} \over r_{{\nu }^{*}} + X_{\cdot , {\nu }^{*}} - 1} \log \Et ^U \Big[ {F^{*} (U, 0, \k ) \over F^{*} (U, 1, \k )} \Big] \non \\
&= {X_{\cdot , {\nu }^{*}} \over r_{{\nu }^{*}} + X_{\cdot , {\nu }^{*}} - 1} \log \Et ^U \Big[ 1 + {\ga _{{\nu }^{*}} U \over \st _{{\nu }^{*}}^{*} (1) + r_{{\nu }^{*}} + a_{0, {\nu }^{*}}} \Big] \non \\
&\le \log \Et ^U \Big[ 1 + {X_{\cdot , {\nu }^{*}} \over r_{{\nu }^{*}} + X_{\cdot , {\nu }^{*}} - 1} {\ga _{{\nu }^{*}} U \over \st _{{\nu }^{*}}^{*} (1) + r_{{\nu }^{*}} + a_{0, {\nu }^{*}}} \Big] \text{,} \label{tmultinp2.3} 
\end{align}
where the inequality follows since $0 \le X_{\cdot , {\nu }^{*}} / ( r_{{\nu }^{*}} + X_{\cdot , {\nu }^{*}} - 1) \le 1$ by assumption. 
By integration by parts, 
\begin{align}
&( \al + 1) \int_{0}^{\infty } u d\mut (u) = \int_{0}^{\infty } \Big[ ( \al + 1) u^{\al } e^{- \be u} \Big\{ \prod_{\nu = 1}^{N} {\Ga ( \ga _{\nu } u + r_{\nu } + a_{0, \nu } ) \Ga ( r_{\nu } + a_{0, \nu } + X_{\cdot , \nu } + a_{\cdot , \nu } ) \over \Ga ( \ga _{\nu } u + r_{\nu } + a_{0, \nu } + X_{\cdot , \nu } + a_{\cdot , \nu } ) \Ga ( r_{\nu } + a_{0, \nu } )} \Big\} \non \\
&\quad \times F^{*} (u, 1, \k ) {\ga _{{\nu }^{*}} u + k_{{\nu }^{*}} + r_{{\nu }^{*}} + a_{0, {\nu }^{*}} + X_{\cdot , {\nu }^{*}} + a_{\cdot , {\nu }^{*}} - 1 \over k_{{\nu }^{*}} + r_{{\nu }^{*}} + a_{0, {\nu }^{*}} + X_{\cdot , {\nu }^{*}} + a_{\cdot , {\nu }^{*}} - 1} \Big] du \non \\
&= \int_{0}^{\infty } \Big( u^{\al + 1} e^{- \be u} \Big\{ \prod_{\nu = 1}^{N} {\Ga ( \ga _{\nu } u + r_{\nu } + a_{0, \nu } ) \Ga ( r_{\nu } + a_{0, \nu } + X_{\cdot , \nu } + a_{\cdot , \nu } ) \over \Ga ( \ga _{\nu } u + r_{\nu } + a_{0, \nu } + X_{\cdot , \nu } + a_{\cdot , \nu } ) \Ga ( r_{\nu } + a_{0, \nu } )} \Big\} \non \\
&\quad \times F^{*} (u, 1, \k ) {\ga _{{\nu }^{*}} u + k_{{\nu }^{*}} + r_{{\nu }^{*}} + a_{0, {\nu }^{*}} + X_{\cdot , {\nu }^{*}} + a_{\cdot , {\nu }^{*}} - 1 \over k_{{\nu }^{*}} + r_{{\nu }^{*}} + a_{0, {\nu }^{*}} + X_{\cdot , {\nu }^{*}} + a_{\cdot , {\nu }^{*}} - 1} \non \\
&\quad \times \Big[ \be + \sum_{\nu = 1}^{N} \ga _{\nu } \{ \psi ( \ga _{\nu } u + r_{\nu } + a_{0, \nu } + X_{\cdot , \nu } + a_{\cdot , \nu } ) - \psi ( \ga _{\nu } u + r_{\nu } + a_{0, \nu } ) \} \non \\
&\quad - \sum_{\nu = 1}^{N} \ga _{\nu } \{ \psi ( \ga _{\nu } u + \st _{\nu }^{*} (1) + r_{\nu } + a_{0, \nu } ) - \psi ( \ga _{\nu } u + k_{\nu } + r_{\nu } + a_{0, \nu } + X_{\cdot , \nu } + a_{\cdot , \nu } ) \non \\
&\quad + \psi ( \ga _{\nu } u + r_{\nu } + a_{0, \nu } + X_{\cdot , \nu } + a_{\cdot , \nu } ) - \psi ( \ga _{\nu } u + r_{\nu } + a_{0, \nu } ) \} \non \\
&\quad - {\ga _{{\nu }^{*}} \over \ga _{{\nu }^{*}} u + k_{{\nu }^{*}} + r_{{\nu }^{*}} + a_{0, {\nu }^{*}} + X_{\cdot , {\nu }^{*}} + a_{\cdot , {\nu }^{*}} - 1} \Big] \Big) du \non \\
&= \int_{0}^{\infty } \Big( u^2 \Big[ \be + \sum_{\nu = 1}^{N} \ga _{\nu } \{ \psi ( \ga _{\nu } u + k_{\nu } + r_{\nu } + a_{0, \nu } + X_{\cdot , \nu } + a_{\cdot , \nu } ) - \psi ( \ga _{\nu } u + \st _{\nu }^{*} (1) + r_{\nu } + a_{0, \nu } ) \} \non \\
&\quad - {\ga _{{\nu }^{*}} \over \ga _{{\nu }^{*}} u + k_{{\nu }^{*}} + r_{{\nu }^{*}} + a_{0, {\nu }^{*}} + X_{\cdot , {\nu }^{*}} + a_{\cdot , {\nu }^{*}} - 1} \Big] \Big) d\mut (u) \text{.} \non 
\end{align}
Therefore, by Lemma 7 of Hamura and Kubokawa (2020b), 
\begin{align}
( \al + 1) \int_{0}^{\infty } u d\mut (u) &\ge \int_{0}^{\infty } u^2 [ \be + \ga _{{\nu }^{*}} \{ \psi ( \ga _{{\nu }^{*}} u + k_{{\nu }^{*}} + r_{{\nu }^{*}} + a_{0, {\nu }^{*}} + X_{\cdot , {\nu }^{*}} + a_{\cdot , {\nu }^{*}} - 1) \non \\
&\quad - \psi ( \ga _{{\nu }^{*}} u + \st _{{\nu }^{*}}^{*} (1) + r_{{\nu }^{*}} + a_{0, {\nu }^{*}} ) \} ] d\mut (u) \non \\
&\ge \int_{0}^{\infty } u^2 \Big\{ \be + \ga _{{\nu }^{*}} {k_{{\nu }^{*}} - \st _{{\nu }^{*}}^{*} (1) + X_{\cdot , {\nu }^{*}} + a_{\cdot , {\nu }^{*}} - 1 \over \ga _{{\nu }^{*}} u + k_{{\nu }^{*}} + r_{{\nu }^{*}} + a_{0, {\nu }^{*}} + X_{\cdot , {\nu }^{*}} + a_{\cdot , {\nu }^{*}} - 1} \Big\} d\mut (u) \non \\
&\ge ( \be + \ga _{{\nu }^{*}} ) \int_{0}^{\infty } u^2 {k_{{\nu }^{*}} - \st _{{\nu }^{*}}^{*} (1) + X_{\cdot , {\nu }^{*}} + a_{\cdot , {\nu }^{*}} - 1 \over \ga _{{\nu }^{*}} u + k_{{\nu }^{*}} + r_{{\nu }^{*}} + a_{0, {\nu }^{*}} + X_{\cdot , {\nu }^{*}} + a_{\cdot , {\nu }^{*}} - 1} d\mut (u) \text{,} \non 
\end{align}
where the third inequality follows since $k_{{\nu }^{*}} \ge \st _{{\nu }^{*}}^{*} (0) = \st _{{\nu }^{*}}^{*} (1) + 1$, and this implies that 
\begin{align}
\Et ^U \Big[ {U^2 \over \ga _{{\nu }^{*}} U + k_{{\nu }^{*}} + r_{{\nu }^{*}} + a_{0, {\nu }^{*}} + X_{\cdot , {\nu }^{*}} + a_{\cdot , {\nu }^{*}} - 1} \Big] &\le {( \al + 1) / ( \be + \ga _{{\nu }^{*}} ) \over k_{{\nu }^{*}} - \st _{{\nu }^{*}}^{*} (1) + X_{\cdot , {\nu }^{*}} + a_{\cdot , {\nu }^{*}} - 1} \Et ^U [ U ] \text{.} \label{tmultinp3} 
\end{align}
When $X_{{\nu }^{*}} \ge 1$, we have, by (\ref{eq:assumption_multin}), 
\begin{align}
\Big\{ {( \al + 1) \ga _{{\nu }^{*}} \over \be + \ga _{{\nu }^{*}}} - a_{\cdot , {\nu }^{*}} \Big\} ( r_{{\nu }^{*}} - 1) \le X_{{\nu }^{*}} \Big\{ - {( \al + 1) \ga _{{\nu }^{*}} \over \be + \ga _{{\nu }^{*}}} - k_{{\nu }^{*}} - a_{0, {\nu }^{*}} \Big\} \text{,} \non 
\end{align}
which implies that 
\begin{align}
\ga _{{\nu }^{*}} {( \al + 1) / ( \be + \ga _{{\nu }^{*}} ) \over k_{{\nu }^{*}} - \st _{{\nu }^{*}}^{*} (1) + X_{\cdot , {\nu }^{*}} + a_{\cdot , {\nu }^{*}} - 1} \le 1 - {X_{\cdot , {\nu }^{*}} \over r_{{\nu }^{*}} + X_{\cdot , {\nu }^{*}} - 1} {k_{{\nu }^{*}} + r_{{\nu }^{*}} + a_{0, {\nu }^{*}} + X_{\cdot , {\nu }^{*}} + a_{\cdot , {\nu }^{*}} - 1 \over k_{{\nu }^{*}} - \st _{{\nu }^{*}}^{*} (1) + X_{\cdot , {\nu }^{*}} + a_{\cdot , {\nu }^{*}} - 1} \label{tmultinp4} 
\end{align}
since $k_{{\nu }^{*}} \ge \st _{{\nu }^{*}}^{*} (1) + 1$. 
From (\ref{tmultinp3}) and (\ref{tmultinp4}), it follows that when $X_{{\nu }^{*}} \ge 1$, 
\begin{align}
&\Et ^U \Big[ {\ga _{{\nu }^{*}} U^2 \over \ga _{{\nu }^{*}} U + k_{{\nu }^{*}} + r_{{\nu }^{*}} + a_{0, {\nu }^{*}} + X_{\cdot , {\nu }^{*}} + a_{\cdot , {\nu }^{*}} - 1} \Big] \non \\
&\le \ga _{{\nu }^{*}} {( \al + 1) / ( \be + \ga _{{\nu }^{*}} ) \over k_{{\nu }^{*}} - \st _{{\nu }^{*}}^{*} (1) + X_{\cdot , {\nu }^{*}} + a_{\cdot , {\nu }^{*}} - 1} \Et ^U [ U ] \non \\
&\le \Big\{ 1 - {X_{\cdot , {\nu }^{*}} \over r_{{\nu }^{*}} + X_{\cdot , {\nu }^{*}} - 1} {k_{{\nu }^{*}} + r_{{\nu }^{*}} + a_{0, {\nu }^{*}} + X_{\cdot , {\nu }^{*}} + a_{\cdot , {\nu }^{*}} - 1 \over k_{{\nu }^{*}} - \st _{{\nu }^{*}}^{*} (1) + X_{\cdot , {\nu }^{*}} + a_{\cdot , {\nu }^{*}} - 1} \Big\} \Et ^U [ U ] \text{,} \non 
\end{align}
which can be rewritten as 
\begin{align}
&{X_{\cdot , {\nu }^{*}} \over r_{{\nu }^{*}} + X_{\cdot , {\nu }^{*}} - 1} \Et ^U \Big[ {U \over k_{{\nu }^{*}} - \st _{{\nu }^{*}}^{*} (1) + X_{\cdot , {\nu }^{*}} + a_{\cdot , {\nu }^{*}} - 1} \Big] \non \\
&\le \Et ^U \Big[ {U \over k_{{\nu }^{*}} + r_{{\nu }^{*}} + a_{0, {\nu }^{*}} + X_{\cdot , {\nu }^{*}} + a_{\cdot , {\nu }^{*}} - 1} \Big( 1 - {\ga _{{\nu }^{*}} U \over \ga _{{\nu }^{*}} U + k_{{\nu }^{*}} + r_{{\nu }^{*}} + a_{0, {\nu }^{*}} + X_{\cdot , {\nu }^{*}} + a_{\cdot , {\nu }^{*}} - 1} \Big) \Big] \non 
\end{align}
or 
\begin{align}
&\Et ^U \Big[ {X_{\cdot , {\nu }^{*}} \over r_{{\nu }^{*}} + X_{\cdot , {\nu }^{*}} - 1} {\ga _{{\nu }^{*}} U \over \st _{{\nu }^{*}}^{*} (1) + r_{{\nu }^{*}} + a_{0, {\nu }^{*}}} \Big] \non \\
&\le \Et ^U \Big[ {\{ k_{{\nu }^{*}} - \st _{{\nu }^{*}}^{*} (1) + X_{\cdot , {\nu }^{*}} + a_{\cdot , {\nu }^{*}} - 1 \} \ga _{{\nu }^{*}} U \over \{ \st _{{\nu }^{*}}^{*} (1) + r_{{\nu }^{*}} + a_{0, {\nu }^{*}} \} ( \ga _{{\nu }^{*}} U + k_{{\nu }^{*}} + r_{{\nu }^{*}} + a_{0, {\nu }^{*}} + X_{\cdot , {\nu }^{*}} + a_{\cdot , {\nu }^{*}} - 1)} \Big] \text{.} \label{tmultinp5} 
\end{align}
Thus, by (\ref{tmultinp2.1}), (\ref{tmultinp2.2}), (\ref{tmultinp2.3}), and (\ref{tmultinp5}), 
\begin{align}
\sum_{\tilde{i} = 0}^{1} \pt _{\tilde{i} , {\nu }^{*}} E[ - \log E^U [ F^{*} (U, \tilde{i} , \k ) ] ] &< E \Big[ - \log {\int_{0}^{\infty } d\mut (u) \over \int_{0}^{\infty } d\mu (u)} \Big] \non \\
&= E[ - \log E^U [ F^{*} (U, 1, \k - \e _{{\nu }^{*}}^{(N)} ) ] ] \text{.} \label{tmultinp6} 
\end{align}

Finally, applying (\ref{tmultinp6}) to (\ref{tmultinp2.05}) sequentially, we obtain 
\begin{align}
\De ^{( \al , \be , \bga , \a _0 , \a )} < \dots < 0 \text{.} \non 
\end{align}
This completes the proof. 
\hfill$\Box$

\bigskip

\noindent
{\bf Proof of Theorem \ref{thm:NM}.} \ \ By (\ref{eq:mass_Y}), (\ref{eq:bpm_Y}), and (\ref{eq:risk_Y}), 
\begin{align}
R( \p , \gh ^{( \pi )} ) &= E \Big[ \log \Big\{ \prod_{\nu = 1}^{n} \Big( {p_{0, \nu }}^{s_{\nu }} \prod_{i = 1}^{m_{\nu }} {p_{i, \nu }}^{Y_{i, \nu }} \Big) \Big\} \Big] \non \\
&\quad + E \Big[ - \log \frac{ \int_D \pi ( \p ) \big\{ \prod_{\nu = 1}^{N} \big( {p_{0, \nu }}^{s_{\nu } + r_{\nu }} \prod_{i = 1}^{m_{\nu }} {p_{i, \nu }}^{Y_{i, \nu } + X_{i, \nu }} \big) \big\} d\p }{ \int_D \pi ( \p ) \big\{ \prod_{\nu = 1}^{N} \big( {p_{0, \nu }}^{r_{\nu }} \prod_{i = 1}^{m_{\nu }} {p_{i, \nu }}^{X_{i, \nu }} \big) \big\} d\p } \Big] \text{,} \label{tNMp1} 
\end{align}
where $Y_{1, \nu } = \dots = Y_{m_{\nu } , \nu } = 0$ if $\nu \in \{ 1, \dots , N \} \cap [n + 1, \infty )$. 
The first term on the right of (\ref{tNMp1}) is 
\begin{align}
&E \Big[ \log \Big\{ \prod_{\nu = 1}^{n} \Big( {p_{0, \nu }}^{s_{\nu }} \prod_{i = 1}^{m_{\nu }} {p_{i, \nu }}^{Y_{i, \nu }} \Big) \Big\} \Big] \non \\
&= \sum_{\nu = 1}^{n} \Big( s_{\nu } \log p_{0, \nu } + \sum_{i = 1}^{m_{\nu }} s_{\nu } {p_{i, \nu } \over p_{0, \nu }} \log p_{i, \nu } \Big) \non \\
&= \sum_{\nu = 1}^{n} s_{\nu } \sum_{k = 1}^{\infty } {1 \over k} \Big( - {p_{\cdot , \nu }}^k + {p_{\cdot , \nu }}^k \sum_{i = 1}^{m_{\nu }} k {p_{i, \nu } \over p_{\cdot , \nu }} \log p_{i, \nu } \Big) \non \\
&= \sum_{\nu = 1}^{n} s_{\nu } \sum_{k = 1}^{\infty } {1 \over k} \sum_{( w_i  )_{i = 1}^{m_{\nu }} \in \mathcal{W} _{\nu , k}} {k ! \over \prod_{i = 1}^{m_{\nu }} w_i !} \Big\{ - \prod_{i = 1}^{m_{\nu }} {p_{i, \nu }}^{w_i} + \Big( \prod_{i = 1}^{m_{\nu }} {p_{i, \nu }}^{w_i} \Big) \sum_{i = 1}^{m_{\nu }} w_i \log p_{i, \nu } \Big\} \text{.} \label{tNMp2} 
\end{align}
On the other hand, since $t_{\nu }$ is a constant if $\nu \in \{ 1, \dots , N \} \cap [n + 1, \infty )$, 
\begin{align}
&E \Big[ - \log \frac{ \int_D \pi ( \p ) \big\{ \prod_{\nu = 1}^{N} \big( {p_{0, \nu }}^{s_{\nu } + r_{\nu }} \prod_{i = 1}^{m_{\nu }} {p_{i, \nu }}^{Y_{i, \nu } + X_{i, \nu }} \big) \big\} d\p }{ \int_D \pi ( \p ) \big\{ \prod_{\nu = 1}^{N} \big( {p_{0, \nu }}^{r_{\nu }} \prod_{i = 1}^{m_{\nu }} {p_{i, \nu }}^{X_{i, \nu }} \big) \big\} d\p } \Big] = \int_{0}^{1} \Big\{ {\partial \over \partial \ta } E[ - \log G( \ta , \Z ( \ta )) ] \Big\} d\ta \non \\
&= \int_{0}^{1} E \Big[ \sum_{\nu = 1}^{n} {t_{\nu }}' ( \ta ) \Big\{ \sum_{k = 1}^{Z_{\cdot , \nu } ( \ta )} {1 \over t_{\nu } ( \ta ) + k - 1} + \log p_{0, \nu } \Big\} \{ - \log G( \ta , \Z ( \ta )) \} - \frac{ \displaystyle {\partial G \over \partial \ta } ( \ta , \Z ( \ta )) }{ \displaystyle G( \ta , \Z ( \ta )) } \Big] d\ta \text{,} \label{tNMp3} 
\end{align}
where 
\begin{align}
G( \ta , (( z_{i, \nu } )_{i = 1}^{m_{\nu }} )_{\nu = 1, \dots , N} ) &= \int_D \pi ( \p ) \Big[ \prod_{\nu = 1}^{N} \Big\{ {p_{0, \nu }}^{t_{\nu } ( \ta )} \prod_{i = 1}^{m_{\nu }} {p_{i, \nu }}^{z_{i, \nu }} \Big\} \Big] d\p \non 
\end{align}
for $(( z_{i, \nu } )_{i = 1}^{m_{\nu }} )_{\nu = 1, \dots , N} \in {\mathbb{N} _0}^{m_1} \times \dots \times {\mathbb{N} _0}^{m_N}$ and where $Z_{\cdot , \nu } ( \ta ) = \sum_{i = 1}^{m_{\nu }} Z_{i, \nu } ( \ta )$ for $\nu = 1, \dots , N$ for each $\ta \in [0, 1]$. 

Fix $\ta \in [0, 1]$. 
Then 
\begin{align}
&E \Big[ \Big\{ {\partial G \over \partial \ta } ( \ta , \Z ( \ta )) \Big\} / G( \ta , \Z ( \ta )) \Big] \non \\
&= E \Big[ \int_D \pi ( \p ) \Big[ \Big\{ \sum_{\nu = 1}^{N} {t_{\nu }}' ( \ta ) \log p_{0, \nu } \Big\} \prod_{\nu = 1}^{N} \Big\{ {p_{0, \nu }}^{t_{\nu } ( \ta )} \prod_{i = 1}^{m_{\nu }} {p_{i, \nu }}^{Z_{i, \nu } ( \ta )} \Big\} \Big] d\p / G( \ta , \Z ( \ta )) \Big] \non \\
&= - \sum_{\nu = 1}^{n} {t_{\nu }}' ( \ta ) \sum_{k = 1}^{\infty } {1 \over k} E \Big[ \int_D \pi ( \p ) \Big[ {p_{\cdot , \nu }}^k \prod_{{\nu }' = 1}^{N} \Big\{ {p_{0, {\nu }'}}^{t_{{\nu }'} ( \ta )} \prod_{i = 1}^{m_{{\nu }'}} {p_{i, {\nu }'}}^{Z_{i, {\nu }'} ( \ta )} \Big\} \Big] d\p / G( \ta , \Z ( \ta )) \Big] \non \\
&= - \sum_{\nu = 1}^{n} {t_{\nu }}' ( \ta ) \sum_{k = 1}^{\infty } {1 \over k} \sum_{( w_i  )_{i = 1}^{m_{\nu }} \in \mathcal{W} _{\nu , k}} {k ! \over \prod_{i = 1}^{m_{\nu }} w_i !} E \Big[ \int_D \pi ( \p ) \Big[ \Big( \prod_{i = 1}^{m_{\nu }} {p_{i, \nu }}^{w_i} \Big) \prod_{{\nu }' = 1}^{N} \Big\{ {p_{0, {\nu }'}}^{t_{{\nu }'} ( \ta )} \prod_{i = 1}^{m_{{\nu }'}} {p_{i, {\nu }'}}^{Z_{i, {\nu }'} ( \ta )} \Big\} \Big] d\p \non \\
&\quad / G( \ta , \Z ( \ta )) \Big] \text{.} \label{tNMp4} 
\end{align}
On the other hand, by Lemmas 2.1 and 2.2 of Hamura and Kubokawa (2019a), we have for any $\nu = 1, \dots , n$, 
\begin{align}
&E \Big[ \Big\{ \sum_{k = 1}^{Z_{\cdot , \nu } ( \ta )} {1 \over t_{\nu } ( \ta ) + k - 1} + \log p_{0, \nu } \Big\} \{ - \log G( \ta , \Z ( \ta )) \} \Big] \non \\
&= E \Big[ \Big\{ \sum_{k = 1}^{Z_{\cdot , \nu } ( \ta )} {1 \over k} {Z_{\cdot , \nu } ( \ta ) \dotsm \{ Z_{\cdot , \nu } ( \ta ) - k + 1 \} \over \{ t_{\nu } ( \ta ) + Z_{\cdot , \nu } ( \ta ) - 1 \} \dotsm \{ t_{\nu } ( \ta ) + Z_{\cdot , \nu } ( \ta ) - k \} } + \log p_{0, \nu } \Big\} \{ - \log G( \ta , \Z ( \ta )) \} \Big] \non \\
&= \sum_{k = 1}^{\infty } {1 \over k} {p_{\cdot , \nu }}^k E[ E[ - \log G( \ta , \Z ( \ta )) | \Z _{\cdot } ( \ta ) + k \e _{\nu }^{(N)} ] - \{ - \log G( \ta , \Z ( \ta )) \} ] \text{,} \non 
\end{align}
where $\Z _{\cdot } ( \ta ) = ( Z_{\cdot , \nu } ( \ta ))_{\nu = 1}^{N}$. 
Now, fix $k \in \mathbb{N}$. 
Let $\W _{\nu }$, %
$\nu = 1, \dots , N$, be mutually independent multinomial variables such that for each $\nu = 1, \dots , N$, the probability mass function of $\W _{\nu } | Z_{\cdot , \nu } ( \ta )$ is given by 
\begin{align}
{Z_{\cdot , \nu } ( \ta ) ! \over \prod_{i = 1}^{m_{\nu }} w_{i, \nu } !} \prod_{i = 1}^{m_{\nu }} \Big( {p_{i, \nu } \over p_{\cdot , \nu }} \Big) ^{w_{i, \nu }} \non 
\end{align}
for $( w_{i, \nu } )_{i = 1}^{m_{\nu }} \in \mathcal{W} _{\nu , Z_{\cdot , \nu } ( \ta )}$. 
Let $\W _{\nu }^{*}$, $\nu = 1, \dots , N$, be independent multinomial variable with mass functions 
\begin{align}
{k ! \over \prod_{i = 1}^{m_{\nu }} w_{i, \nu }^{*} !} \prod_{i = 1}^{m_{\nu }} \Big( {p_{i, \nu } \over p_{\cdot , \nu }} \Big) ^{w_{i, \nu }^{*}} \text{,} \non 
\end{align}
$( w_{i, \nu }^{*} )_{i = 1}^{m_{\nu }} \in \mathcal{W} _{\nu , k}$, $\nu = 1, \dots , N$, respectively. 
Then, for any $\nu = 1, \dots , N$, 
\begin{align}
&E[ - \log G( \ta , \Z ( \ta )) | \Z _{\cdot } ( \ta ) + k \e _{\nu }^{(N)} ] \non \\
&= E[ - \log G( \ta , ( \W _{{\nu }'} + \de _{\nu , {\nu }'}^{(N)} \W _{{\nu }'}^{*} )_{{\nu }' = 1, \dots , N} ) | \Z _{\cdot } ( \ta ) ] \non \\
&= \sum_{( w_{i, \nu }^{*} )_{i = 1}^{m_{\nu }} \in \mathcal{W} _{\nu , k}} {k ! \over \prod_{i = 1}^{m_{\nu }} w_{i, \nu }^{*} !} \Big\{ \prod_{i = 1}^{m_{\nu }} \Big( {p_{i, \nu } \over p_{\cdot , \nu }} \Big) ^{w_{i, \nu }^{*}} \Big\} E[ - \log G( \ta , ( \W _{{\nu }'} + \de _{\nu , {\nu }'}^{(N)} ( w_{i, \nu }^{*} )_{i = 1}^{m_{\nu }} )_{{\nu }' = 1, \dots , N} ) | \Z _{\cdot } ( \ta ) ] \non \\
&= \sum_{( w_{i, \nu }^{*} )_{i = 1}^{m_{\nu }} \in \mathcal{W} _{\nu , k}} {k ! \over \prod_{i = 1}^{m_{\nu }} w_{i, \nu }^{*} !} \Big\{ \prod_{i = 1}^{m_{\nu }} \Big( {p_{i, \nu } \over p_{\cdot , \nu }} \Big) ^{w_{i, \nu }^{*}} \Big\} E[ - \log G( \ta , ( \Z _{{\nu }'} ( \ta ) + \de _{\nu , {\nu }'}^{(N)} ( w_{i, \nu }^{*} )_{i = 1}^{m_{\nu }} )_{{\nu }' = 1, \dots , N} ) | \Z _{\cdot } ( \ta ) ] \non 
\end{align}
and therefore 
\begin{align}
&E[ E[ - \log G( \ta , \Z ( \ta )) | \Z _{\cdot } ( \ta ) + k \e _{\nu }^{(N)} ] ] \non \\
&= {1 \over {p_{\cdot , \nu }}^k} \sum_{( w_i )_{i = 1}^{m_{\nu }} \in \mathcal{W} _{\nu , k}} {k ! \over \prod_{i = 1}^{m_{\nu }} w_i !} \Big( \prod_{i = 1}^{m_{\nu }} {p_{i, \nu }}^{w_i} \Big) E[ - \log G( \ta , ( \Z _{{\nu }'} ( \ta ) + \de _{\nu , {\nu }'}^{(N)} ( w_i )_{i = 1}^{m_{\nu }} )_{{\nu }' = 1, \dots , N} ) ] \text{.} \non 
\end{align}
Since $k$ is arbitrarily chosen, it follows that 
\begin{align}
&E \Big[ \Big\{ \sum_{k = 1}^{Z_{\cdot , \nu } ( \ta )} {1 \over t_{\nu } ( \ta ) + k - 1} + \log p_{0, \nu } \Big\} \{ - \log G( \ta , \Z ( \ta )) \} \Big] \non \\
&= \sum_{k = 1}^{\infty } {1 \over k} \Big\{ \sum_{( w_i )_{i = 1}^{m_{\nu }} \in \mathcal{W} _{\nu , k}} {k ! \over \prod_{i = 1}^{m_{\nu }} w_i !} \Big( \prod_{i = 1}^{m_{\nu }} {p_{i, \nu }}^{w_i} \Big) E[ - \log G( \ta , ( \Z _{{\nu }'} ( \ta ) + \de _{\nu , {\nu }'}^{(N)} ( w_i )_{i = 1}^{m_{\nu }} )_{{\nu }' = 1, \dots , N} ) ] \non \\
&\quad - \sum_{( w_i )_{i = 1}^{m_{\nu }} \in \mathcal{W} _{\nu , k}} {k ! \over \prod_{i = 1}^{m_{\nu }} w_i !} \Big( \prod_{i = 1}^{m_{\nu }} {p_{i, \nu }}^{w_i} \Big) E[ - \log G( \ta , \Z ( \ta )) ] \Big\} \non \\
&= \sum_{k = 1}^{\infty } {1 \over k} \sum_{( w_i )_{i = 1}^{m_{\nu }} \in \mathcal{W} _{\nu , k}} {k ! \over \prod_{i = 1}^{m_{\nu }} w_i !} \Big( \prod_{i = 1}^{m_{\nu }} {p_{i, \nu }}^{w_i} \Big) E \Big[ - \log {G( \ta , ( \Z _{{\nu }'} ( \ta ) + \de _{\nu , {\nu }'}^{(N)} ( w_i )_{i = 1}^{m_{\nu }} )_{{\nu }' = 1, \dots , N} ) \over G( \ta , \Z ( \ta ))} \Big] \text{.} \label{tNMp5} 
\end{align}

Finally, combining (\ref{tNMp1}), (\ref{tNMp2}), (\ref{tNMp3}), (\ref{tNMp4}), and (\ref{tNMp5}), we obtain 
\begin{align}
R( \p , \gh ^{( \pi )} ) %
&= \int_{0}^{1} \Big[ \sum_{\nu = 1}^{n} {t_{\nu }}' ( \ta ) \sum_{k = 1}^{\infty } {1 \over k} \sum_{( w_i )_{i = 1}^{m_{\nu }} \in \mathcal{W} _{\nu , k}} \Big\{ {k ! \over \prod_{i = 1}^{m_{\nu }} w_i !} \non \\
&\quad \times \Big( - \prod_{i = 1}^{m_{\nu }} {p_{i, \nu }}^{w_i} + \Big( \prod_{i = 1}^{m_{\nu }} {p_{i, \nu }}^{w_i} \Big) \sum_{i = 1}^{m_{\nu }} w_i \log p_{i, \nu } \non \\
&\quad + \Big( \prod_{i = 1}^{m_{\nu }} {p_{i, \nu }}^{w_i} \Big) E \Big[ - \log {G( \ta , ( \Z _{{\nu }'} ( \ta ) + \de _{\nu , {\nu }'}^{(N)} ( w_i )_{i = 1}^{m_{\nu }} )_{{\nu }' = 1, \dots , N} ) \over G( \ta , \Z ( \ta ))} \Big] \non \\
&\quad + E \Big[ \int_D \pi ( \p ) \Big[ \Big( \prod_{i = 1}^{m_{\nu }} {p_{i, \nu }}^{w_i} \Big) \prod_{{\nu }' = 1}^{N} \Big\{ {p_{0, {\nu }'}}^{t_{{\nu }'} ( \ta )} \prod_{i = 1}^{m_{{\nu }'}} {p_{i, {\nu }'}}^{Z_{i, {\nu }'} ( \ta )} \Big\} \Big] d\p / G( \ta , \Z ( \ta )) \Big]  \Big) \Big\} \Big] d\ta \non \\
&= \int_{0}^{1} \Big[ \sum_{\nu = 1}^{n} {t_{\nu }}' ( \ta ) \sum_{k = 1}^{\infty } {1 \over k} \sum_{( w_i )_{i = 1}^{m_{\nu }} \in \mathcal{W} _{\nu , k}} \Big\{ {k ! \over \prod_{i = 1}^{m_{\nu }} w_i !} \non \\
&\quad \times E \Big[ L^{\rm{KL}} \Big( {G( \ta , ( \Z _{{\nu }'} ( \ta ) + \de _{\nu , {\nu }'}^{(N)} ( w_i )_{i = 1}^{m_{\nu }} )_{{\nu }' = 1, \dots , N} ) \over G( \ta , \Z ( \ta ))} , \prod_{i = 1}^{m_{\nu }} {p_{i, \nu }}^{w_i} \Big) \Big] \Big\} \Big] d\ta \text{.} \non 
\end{align}
Thus, 
\begin{align}
&R( \p , \gh ^{( \pi )} ) \non \\
&= \int_{0}^{1} \Big\{ \sum_{\nu = 1}^{n} {t_{\nu }}' ( \ta ) \sum_{k = 1}^{\infty } {1 \over k} \sum_{( w_i )_{i = 1}^{m_{\nu }} \in \mathcal{W} _{\nu , k}} {k ! \over \prod_{i = 1}^{m_{\nu }} w_i !} E \Big[ L^{\rm{KL}} \Big( E_{\pi } \Big[ \prod_{i = 1}^{m_{\nu }} {p_{i, \nu }}^{w_i} \Big| \Z ( \ta ) \Big] , \prod_{i = 1}^{m_{\nu }} {p_{i, \nu }}^{w_i} \Big) \Big] \Big\} d\ta \text{,} \non 
\end{align}
which is the desired result. 
\hfill$\Box$

\bigskip

\noindent
{\bf Proof of Corollary \ref{cor:NM}.} \ \ By Theorem \ref{thm:NM}, we have 
\begin{align}
&R( \p , \gh ^{( \pi _{M, \bgat , \a _0 , \a } )} ) - R( \p , \gh ^{( \pi _{\a _0 , \a } )} ) \non \\
&= \int_{0}^{1} \Big\{ \sum_{\nu = 1}^{n} {t_{\nu }}' ( \ta ) \sum_{k = 1}^{\infty } {1 \over k} \sum_{( w_i )_{i = 1}^{m_{\nu }} \in \mathcal{W} _{\nu , k}} {k ! \over \prod_{i = 1}^{m_{\nu }} w_i !} E \Big[ L^{\rm{KL}} \Big( E_{\pi _{M, \bgat , \a _0 , \a }} \Big[ \prod_{i = 1}^{m_{\nu }} {p_{i, \nu }}^{w_i} \Big| \Z ( \ta ) \Big] , \prod_{i = 1}^{m_{\nu }} {p_{i, \nu }}^{w_i} \Big) \non \\
&\quad - L^{\rm{KL}} \Big( E_{\pi _{\a _0 , \a }} \Big[ \prod_{i = 1}^{m_{\nu }} {p_{i, \nu }}^{w_i} \Big| \Z ( \ta ) \Big] , \prod_{i = 1}^{m_{\nu }} {p_{i, \nu }}^{w_i} \Big) \Big] \Big\} d\ta \text{.} \non 
\end{align}
Fix $\ta \in [0, 1]$, $\nu = 1, \dots , n$, and $k \in \mathbb{N}$. 
Then 
\begin{align}
&\sum_{( w_i )_{i = 1}^{m_{\nu }} \in \mathcal{W} _{\nu , k}} {k ! \over \prod_{i = 1}^{m_{\nu }} w_i !} E \Big[ L^{\rm{KL}} \Big( E_{\pi _{M, \bgat , \a _0 , \a }} \Big[ \prod_{i = 1}^{m_{\nu }} {p_{i, \nu }}^{w_i} \Big| \Z ( \ta ) \Big] , \prod_{i = 1}^{m_{\nu }} {p_{i, \nu }}^{w_i} \Big) \non \\
&\quad - L^{\rm{KL}} \Big( E_{\pi _{\a _0 , \a }} \Big[ \prod_{i = 1}^{m_{\nu }} {p_{i, \nu }}^{w_i} \Big| \Z ( \ta ) \Big] , \prod_{i = 1}^{m_{\nu }} {p_{i, \nu }}^{w_i} \Big) \Big] \non \\
&= \sum_{( w_i )_{i = 1}^{m_{\nu }} \in \mathcal{W} _{\nu , k}} {k ! \over \prod_{i = 1}^{m_{\nu }} w_i !} E \Big[ E_{\pi _{M, \bgat , \a _0 , \a }} \Big[ \prod_{i = 1}^{m_{\nu }} {p_{i, \nu }}^{w_i} \Big| \Z ( \ta ) \Big] - E_{\pi _{\a _0 , \a }} \Big[ \prod_{i = 1}^{m_{\nu }} {p_{i, \nu }}^{w_i} \Big| \Z ( \ta ) \Big] \non \\
&\quad - \Big( \prod_{i = 1}^{m_{\nu }} {p_{i, \nu }}^{w_i} \Big) \log \Big\{ E_{\pi _{M, \bgat , \a _0 , \a }} \Big[ \prod_{i = 1}^{m_{\nu }} {p_{i, \nu }}^{w_i} \Big| \Z ( \ta ) \Big] / E_{\pi _{\a _0 , \a }} \Big[ \prod_{i = 1}^{m_{\nu }} {p_{i, \nu }}^{w_i} \Big| \Z ( \ta ) \Big] \Big\} \Big] \text{.} \non 
\end{align}
Note that 
\begin{align}
&\sum_{( w_i )_{i = 1}^{m_{\nu }} \in \mathcal{W} _{\nu , k}} {k ! \over \prod_{i = 1}^{m_{\nu }} w_i !} E \Big[ E_{\pi _{M, \bgat , \a _0 , \a }} \Big[ \prod_{i = 1}^{m_{\nu }} {p_{i, \nu }}^{w_i} \Big| \Z ( \ta ) \Big] - E_{\pi _{\a _0 , \a }} \Big[ \prod_{i = 1}^{m_{\nu }} {p_{i, \nu }}^{w_i} \Big| \Z ( \ta ) \Big] \Big] \non \\
&= E[ E_{\pi _{M, \bgat , \a _0 , \a }} [ {p_{\cdot , \nu }}^k | \Z ( \ta ) ] - E_{\pi _{\a _0 , \a }} [ {p_{\cdot , \nu }}^k | \Z ( \ta ) ] ] \non 
\end{align}
and that for all $( w_i )_{i = 1}^{m_{\nu }} \in \mathcal{W} _{\nu , k}$, 
\begin{align}
&E_{\pi _{M, \bgat , \a _0 , \a }} \Big[ \prod_{i = 1}^{m_{\nu }} {p_{i, \nu }}^{w_i} \Big| \Z ( \ta ) \Big] / E_{\pi _{\a _0 , \a }} \Big[ \prod_{i = 1}^{m_{\nu }} {p_{i, \nu }}^{w_i} \Big| \Z ( \ta ) \Big] \non \\
&= \frac{ \displaystyle \int_{0}^{\infty } \Big\{ \prod_{{\nu }' = 1}^{N} {\Ga ( \gat _{{\nu }'} (u) + t_{{\nu }'} ( \ta ) + a_{0, {\nu }'} ) \over \Ga ( \gat _{{\nu }'} (u) + t_{{\nu }'} ( \ta ) + a_{0, {\nu }'} + Z_{\cdot , {\nu }'} ( \ta ) + a_{\cdot , {\nu }'} + \de _{\nu , {\nu }'}^{(N)} k)} \Big\} dM(u) }{ \displaystyle \int_{0}^{\infty } \Big\{ \prod_{{\nu }' = 1}^{N} {\Ga ( \gat _{{\nu }'} (u) + t_{{\nu }'} ( \ta ) + a_{0, {\nu }'} ) \over \Ga ( \gat _{{\nu }'} (u) + t_{{\nu }'} ( \ta ) + a_{0, {\nu }'} + Z_{\cdot , {\nu }'} ( \ta ) + a_{\cdot , {\nu }'} )} \Big\} dM(u) } \non \\
&\quad / \frac{ \displaystyle \prod_{{\nu }' = 1}^{N} {\Ga ( t_{{\nu }'} ( \ta ) + a_{0, {\nu }'} ) \over \Ga ( t_{{\nu }'} ( \ta ) + a_{0, {\nu }'} + Z_{\cdot , {\nu }'} ( \ta ) + a_{\cdot , {\nu }'} + \de _{\nu , {\nu }'}^{(N)} k)} }{ \displaystyle \prod_{{\nu }' = 1}^{N} {\Ga ( t_{{\nu }'} ( \ta ) + a_{0, {\nu }'} ) \over \Ga ( t_{{\nu }'} ( \ta ) + a_{0, {\nu }'} + Z_{\cdot , {\nu }'} ( \ta ) + a_{\cdot , {\nu }'} )} } \non \\
&= E_{\pi _{M, \bgat , \a _0 , \a }} [ {p_{\cdot , \nu }}^k | \Z ( \ta ) ] / E_{\pi _{\a _0 , \a }} [ {p_{\cdot , \nu }}^k | \Z ( \ta ) ] \text{,} \non 
\end{align}
where $Z_{\cdot , {\nu }'} ( \ta ) = \sum_{i = 1}^{m_{{\nu }'}} Z_{i, {\nu }'} ( \ta )$ for ${\nu }' = 1, \dots , N$. 
It follow that 
\begin{align}
&\sum_{( w_i )_{i = 1}^{m_{\nu }} \in \mathcal{W} _{\nu , k}} {k ! \over \prod_{i = 1}^{m_{\nu }} w_i !} E \Big[ L^{\rm{KL}} \Big( E_{\pi _{M, \bgat , \a _0 , \a }} \Big[ \prod_{i = 1}^{m_{\nu }} {p_{i, \nu }}^{w_i} \Big| \Z ( \ta ) \Big] , \prod_{i = 1}^{m_{\nu }} {p_{i, \nu }}^{w_i} \Big) \non \\
&\quad - L^{\rm{KL}} \Big( E_{\pi _{\a _0 , \a }} \Big[ \prod_{i = 1}^{m_{\nu }} {p_{i, \nu }}^{w_i} \Big| \Z ( \ta ) \Big] , \prod_{i = 1}^{m_{\nu }} {p_{i, \nu }}^{w_i} \Big) \Big] \non \\
&= E[ L^{\rm{KL}} ( E_{\pi _{M, \bgat , \a _0 , \a }} [ {p_{\cdot , \nu }}^k | \Z ( \ta ) ], {p_{\cdot , \nu }}^k ) - L^{\rm{KL}} ( E_{\pi _{\a _0 , \a }} [ {p_{\cdot , \nu }}^k | \Z ( \ta ) ], {p_{\cdot , \nu }}^k ) ] \text{.} \non 
\end{align}
This completes the proof. 
\hfill$\Box$

\section*{Acknowledgments}
I would like to thank Professor Tatsuya Kubokawa for his encouragement. 
Research of the author was supported in part by Grant-in-Aid for Scientific Research (20J10427) from Japan Society for the Promotion of Science.

\end{document}